\documentclass{article}

\usepackage[utf8]{inputenc}
\usepackage{amsmath, amsxtra, amsfonts,amscd,amssymb,amsthm,mathtools,bm,mathrsfs}
\usepackage[left=1in,top=1in,right=1in,bottom=1in,letterpaper]{geometry}
\usepackage{multirow,booktabs,makecell}
\usepackage{graphicx,epstopdf}
\usepackage{hyperref}
\usepackage{algorithm,algorithmic}
\usepackage{color}
\allowdisplaybreaks

\newcommand{\bbR}{\mathbb{R}}
\newcommand{\vx}{\mathbf{x}}
\newcommand{\vp}{\mathbf{p}}

\newcommand{\calC}{\mathcal{C}}
\newcommand{\calP}{\mathcal{P}}

\newcommand{\calK}{\mathcal{Y}}
\newcommand{\Cobs}{\calC_{\obs}}
\newcommand{\Cmetric}{\calC_{\metric}}
\newcommand{\vb}{\bm{\beta}}

\newcommand{\vm}{\mathbf{m}}
\newcommand{\vv}{\mathbf{v}}
\newcommand{\mx}{m^x}
\newcommand{\my}{m^y}
\newcommand{\obs}{b}
\newcommand{\metric}{g}
\newcommand{\rhoinit}{\mu_0}
\newcommand{\rhoend}{\mu_1}
\newcommand{\weight}{\gamma}
\newcommand{\costfuncdyn}{L}
\newcommand{\costevo}{\mathcal{F}_I}
\newcommand{\costterm}{\mathcal{F}_T}

\newcommand{\rhotilde}{\widetilde{\rho}}
\newcommand{\vmtilde}{\widetilde{\vm}}
\newcommand{\mxtilde}{\widetilde{m}^x}
\newcommand{\mytilde}{\widetilde{m}^y}
\newcommand{\phitilde}{\widetilde{\phi}}
\newcommand{\obstilde}{\widetilde{\obs}}
\newcommand{\mettilde}{\widetilde{\metric}}

\DeclareMathOperator*{\diag}{diag}
\DeclareMathOperator*{\proj}{proj}
\DeclareMathOperator*{\argmin}{argmin}

\newcommand{\calL}{\mathcal{L}}
\newcommand{\calU}{\mathcal{U}}
\newcommand{\calD}{\mathcal{D}}
\newcommand{\calR}{\mathcal{R}}
\newcommand{\calA}{\mathcal{A}}
\newcommand{\Drho}{\calD_{\rho}}
\newcommand{\Dm}{\calD_{\vm}}
\newcommand{\gd}{{M}}

\newcommand{\vi}{\mathbf{i}}
\newcommand{\idt}{{i_t}}
\newcommand{\idx}{{i_x}}
\newcommand{\idy}{{i_y}}
\newcommand{\id}{{\idx\idy\idt}}
\newcommand{\idtm}{{\idx,\idy,\idt-\half}}
\newcommand{\idxm}{{\idx-\half,\idy,\idt}}
\newcommand{\idym}{{\idx,\idy-\half,\idt}}
\newcommand{\idtp}{{\idx,\idy,\idt+\half}}
\newcommand{\idxp}{{\idx+\half,\idy,\idt}}
\newcommand{\idyp}{{\idx,\idy+\half,\idt}}
\newcommand{\idtpp}{{\idx,\idy,\idt+1}}
\newcommand{\idxpp}{{\idx+1,\idy,\idt}}
\newcommand{\idypp}{{\idx,\idy+1,\idt}}

\newcommand{\idxy}{{\idx,\idy}}
\newcommand{\nt}{{n_t}}
\newcommand{\nx}{{n_x}}
\newcommand{\ny}{{n_y}}
\newcommand{\dt}{\Delta t}
\newcommand{\dx}{\Delta x}
\newcommand{\dy}{\Delta y}

\newcommand{\calG}{\mathcal{G}}
\newcommand{\Grho}{\calG^{\rho}}
\newcommand{\Gmx}{\calG^{\mx}}
\newcommand{\Gmy}{\calG^{\my}}
\newcommand{\Gm}{\calG^{\vm}}
\newcommand{\Gphi}{\calG^{\phi}}

\newcommand{\rhobar}{\overline{\rho}}
\newcommand{\mxbar}{\overline{m^x}}
\newcommand{\mybar}{\overline{m^y}}
\newcommand{\vmbar}{\overline{\vm}}

\newcommand{\It}{I_t}
\newcommand{\Ix}{I_x}
\newcommand{\Iy}{I_y}
\newcommand{\Dt}{D_t}
\newcommand{\Dx}{D_x}
\newcommand{\Dy}{D_y}

\newcommand{\DGrho}{\calD_{\Grho}}
\newcommand{\DGm}{\calD_{\Gm}}

\newcommand{\ptbrho}{\delta_{\rho}}
\newcommand{\ptbmx}{\delta_{\mx}}
\newcommand{\ptbmy}{\delta_{\my}}
\newcommand{\ptbphi}{\delta_{\phi}}
\newcommand{\ptbvm}{\delta_{\vm}}

\newcommand{\stp}{\tau}
\newcommand{\dd}{\mathrm{d}}
\newcommand{\half}{\frac{1}{2}}
\newcommand{\zero}{\mathbf{0}}

\newtheorem{theorem}{Theorem}[section]
\newtheorem{lemma}[theorem]{Lemma}

\theoremstyle{definition}

\newtheorem{problem}[theorem]{Problem}

\theoremstyle{remark}
\newtheorem{remark}[theorem]{Remark}
\newtheorem{assumption}{Assumption}

\newcommand\numberthis{\addtocounter{equation}{1}\tag{\theequation}}

\definecolor{ao}{rgb}{0.0, 0.5, 0.0}

\title{A Bilevel Optimization Method for Inverse Mean-Field Games\thanks{This work is suppoted in part by NSF DMS-2134168.}}
\date{}
\author{
Jiajia Yu \thanks{J. Yu (jiajia.yu@duke.edu) is with the Department of Mathematics, Duke University. } \qquad 
Quan Xiao \thanks{Q. Xiao (xiaoq5@rpi.edu) is with the Department of Electrical, Computer, and Systems Engineering, Rensselaer Polytechnic Institute. } \qquad 
Tianyi Chen \thanks{T. Chen (chent18@rpi.edu) is with the Department of Electrical, Computer, and Systems Engineering, Rensselaer Polytechnic Institute. } \qquad 
Rongjie Lai \thanks{R Lai (lairj@purdue.edu) is with the Department of Mathematics, Purdue University. }
}
\begin{document}

\maketitle

\begin{abstract}

In this paper, we introduce a bilevel optimization framework for addressing inverse mean-field games, alongside an exploration of numerical methods tailored for this bilevel problem. The primary benefit of our bilevel formulation lies in maintaining the convexity of the objective function and the linearity of constraints in the forward problem. Our paper focuses on inverse mean-field games characterized by unknown obstacles and metrics. We show numerical stability for these two types of inverse problems. More importantly, we, for the first time, establish the identifiability of the inverse mean-field game with unknown obstacles via the solution of the resultant bilevel problem. The bilevel approach enables us to employ an alternating gradient-based optimization algorithm with a provable convergence guarantee. To validate the effectiveness of our methods in solving the inverse problems, we have designed comprehensive numerical experiments, providing empirical evidence of its efficacy. 
\end{abstract}

\section{Introduction}
\label{sec: intro}

Mean-field games study the Nash Equilibrium in a non-cooperative game with infinitely many agents. In the game, each agent aims to minimize a combination of dynamic cost, interaction cost, and terminal cost by controlling its own state trajectory. At the Nash Equilibrium, the agents cannot unilaterally reduce their costs.
The theory is proposed in \cite{lasry2007mean,huang2006large,huang2007large} and has attracted increasing attention since then. 

In most existing works, knowing the cost functions is required to solve mean-field games. However, in practice, these cost functions are not always easy to obtain. In contrast, the state distribution, the strategies of agents, and sometimes the value function at the Nash Equilibrium can be observed. Thus, a natural question arises: Can we learn the cost functions from the Nash Equilibrium? We refer to this as the inverse mean-field game problem, and to the original one as the forward problem.

Unlike the forward problem, relatively few studies focus on inverse mean-field games. 
\cite{kachroo2015inverse} derives two traffic flow models as the solutions of non-viscous mean-field games.
\cite{ding2022mean} reconstructs the underlying metric in the dynamic cost and the kernel in the non-local interaction cost from the possibly noisy observation of agent distribution and strategy. \cite{chow2022numerical} learns the running cost from population density and strategy on a given boundary. 
\cite{liu2022inverse,liu2022inversebd,liu2023simultaneously,ren2023unique} establish the theoretical unique identifiability result for a general class of mean-field games, mean-field game boundary problems and multipopulation mean-field games, where infinite pairs of training data are required in the proof. 
Following \cite{klibanov2023lipschitz}, a series of works \cite{klibanov2023lipschitz2,klibanov2023meana,klibanov2023meanb,klibanov2023holder,imanuvilov2023lipschitz} study the stability and uniqueness of inverse mean-field game through Carleman estimates.

In this paper, we study a typical class of forward problems, the potential mean-field games. Applications like crowd motion \cite{ruthotto2020machine} and generative models \cite{lin2020apac} have the formulations of potential mean-field games. 
In a potential mean-field game, the Nash Equilibrium is a pair of agent distribution $\rho$ and strategy $\vm$ minimizing a cost $\calL$ which consists of the dynamic cost $\costfuncdyn$, the interaction cost $\costevo$ and the terminal cost $\costterm$, under a constraint $\calC(\rhoinit)$ for density and strategy evolution dynamics:
\begin{equation}    (\rho^*,\vm^*):=\argmin_{\rho,\vm\in\calC(\rhoinit)}\calL(\rho,\vm;\costfuncdyn,\costevo,\costterm).
\label{eq: intro forward mfg}
\end{equation}
The inverse problem is to identify $\costfuncdyn,\costevo$ or $\costterm$ given $(\rho^*,\vm^*)$.
Typical choices of $\costfuncdyn,\costevo$ and $\costterm$ make \eqref{eq: intro forward mfg} a (strongly) convex optimization problem with linear constraint.  
Taking $\costfuncdyn$ unknown and $\costevo,\costterm$ known as an example, we thus consider the following bilevel optimization problem
\begin{equation}
\begin{aligned}
    \min_{\costfuncdyn}\quad &\calD((\rho^*,\vm^*),(\rho(\costfuncdyn),\vm(\costfuncdyn))) + \calR(\costfuncdyn)\\
    &\text{s.t. }(\rho(\costfuncdyn),\vm(\costfuncdyn)):=\argmin_{\rho,\vm\in\calC(\rhoinit)}\calL(\rho,\vm;\costfuncdyn,\costevo,\costterm).
\end{aligned}
\label{eq: intro inverse mfg}
\end{equation}
Here $\calD$ is a fidelity term and $\calR$ is a regularity term.
Existing works~\cite{ding2022mean, chow2022numerical} use the nonlinear and nonconvex PDE optimality conditions as constraints. Consequently, achieving a theoretical convergence guarantee is challenging. In contrast, we propose a bilevel formulation for inverse mean-field games, which directly incorporates the forward problem as the constraint. This bilevel formulation maintains the desirable convex-linear structure of the forward problem (the lower-level problem) and enables us to adopt a gradient-based bilevel optimization algorithm \cite{chen2021closing, vicol2022implicit,ji2021bilevel,hong2023two,xiao2023alternating} to address the inverse problem \eqref{eq: intro inverse mfg}. Moreover, leveraging this convex-linear structure, we have developed a convergence result, demonstrating that our algorithm converges to the stationary point of the bilevel problem. 

A common question in inverse mean-field games concerns the stability and unique identifiability of the unknown parameter or function relative to the data. In our setting, we ask whether it is possible to uniquely recover the cost functions from a single pair of observations $(\rho^*,\vm^*)$ and whether the recovered cost function continuously depends on these observations. This setup differs significantly from the theoretical work discussed in~\cite{liu2022inverse, liu2022inversebd}. In those studies, the authors demonstrate that if the interaction and terminal costs are local, holomorphic in $\rho(\vx,t)$, and meet zero admissibility conditions, then it is possible to uniquely recover them from infinitely many observations either throughout the full domain or on its boundary. However, in our case, the cost function for a crowd motion model typically does not satisfy the zero admissibility condition. Moreover, obtaining infinitely many observations is not feasible in practice. In this work, we establish stability results for a general model and unique identifiability results for crowd motion models at a discrete level. Specifically, for a general model, we prove that our model can achieve a close solution to the ground truth with noisy observation, and for the crowd motion model, we prove that only one pair of complete observation 
$(\rho^*,\vm^*)$ is sufficient to uniquely recover the obstacle function, up to a constant. Thus, compared to the requirement of infinitely many observations in \cite{liu2022inverse}, our result is more practical and offers insights into what constitutes an effective observation for accurately recovering the ground truth.
 
\textbf{Contribution:}
We summarize our contributions as follows.
\begin{enumerate}
    \item We propose a bilevel optimization framework for modeling inverse mean-field games.
    \item We study a general model of mean-field games and show that the unknown cost parameters continuously depend on the observation of the Nash Equilibrium.
    \item For the crowd motion model, we prove that up to a constant, the ground truth obstacle function is the unique minimizer to the bilevel optimization problem of the inverse mean-field game.
    \item We apply an alternating gradient-based bilevel optimization algorithm to solve inverse mean-field games and prove the algorithm converges to the stationary point of the bilevel problem.
    \item We implement the algorithm and illustrate the effectiveness of our model and algorithm with comprehensive numerical experiments.
\end{enumerate}

\textbf{Organization:}
The paper is organized as follows.
In section \ref{sec: review}, we briefly review the potential mean-field games and provide two examples of forward mean-field game models whose inverse problem will be addressed in this paper.
In section \ref{sec: model}, we provide the bilevel optimization model for inverse mean-field games and discretize the model. We also state the stability of both models and the unique identifiability of the inverse crowd motion model and prove them in section \ref{sec: proof}. In section \ref{sec: alg}, we present the algorithm to solve our bilevel model for inverse mean-field games and prove the convergence in section \ref{sec: proof}.
In section \ref{sec: num egs}, we demonstrate our model and algorithm with experiments. Finally, we conclude our work in section \ref{sec:con}.

\section{A Review of Potential Mean-Field Games}
\label{sec: review}

In this section, we first review potential mean-field games and their optimality conditions~\cite{lasry2007mean,huang2006large,huang2007large}. Then we present two example problems that we would like to solve in the inverse problem setup.

Consider a problem defined spatially on $\Omega\subset\bbR^d$ and temporally on $[0,1]$. $\rho:\Omega\times[0,1]\to\bbR$ is the state density. $\vv:\Omega\times[0,1]\to\bbR^d$ represents the velocity (control) field of the agents and $\vm:=\rho\vv$ the flux.
A potential mean-field game typically has the following formulation:
\begin{equation}
    \min_{(\rho,\vm)\in\calC(\rhoinit)} \calL(\rho,\vm):=\int_0^1\int_\Omega \rho(\vx,t)\costfuncdyn\left(\vx,\frac{\vm(\vx,t)}{\rho(\vx,t)}\right)\dd \vx\dd t + \int_0^1 \costevo(\rho(\cdot,t))\dd t + \costterm(\rho(\cdot,1))
\label{eq: rew mfg}
\end{equation}
with the constraint set being
\begin{equation}
    \calC(\rhoinit):=\{(\rho,\vm):\partial_t\rho + \nabla\cdot\vm = 0, \rho(\cdot,0) = \rhoinit,\vm\cdot\mathbf{n}=0\text{ for }\vx\in\partial\Omega,\rho(\cdot,\cdot)\geq0\}.
\label{eq: rew cstr set}
\end{equation}
where $\mathbf{n}$ is the normal direction on the boundary $\partial\Omega$. It is clear to see that any pair of $(\rho,\vm)\in\calC(\rhoinit)$ satisfies mass conservation and zero boundary flux condition with the initial density of $\rho$ being $\rhoinit$.
In this objective function, $\costfuncdyn:\Omega\times\bbR^d\to\bbR$ models the dynamic cost, $\costevo:\calP(\Omega)\to\bbR$ the interaction cost and $\costterm:\calP(\Omega)\to\bbR$ the terminal cost.

To derive the optimality condition of \eqref{eq: rew mfg}, we introduce the Lagrangian multiplier $\phi$ and formulate the Lagrangian
\begin{equation}
\begin{aligned}
    \calA(\rho,\vm,\phi):=&\calL(\rho,\vm) - \int_0^1\int_\Omega \phi(\vx,t)\left(\partial_t\rho(\vx,t)+\nabla\cdot\vm(\vx,t)\right)\dd\vx\dd t\\
    =&\calL(\rho,\vm) 
      + \int_0^1\int_\Omega\rho(\vx,t)\partial_t\phi(\vx,t)\dd\vx\dd t
      + \int_0^1\int_\Omega\vm(\vx,t)\cdot\nabla\phi(\vx,t)\dd\vx\dd t\\
      &- \int_\Omega\phi(\vx,1)\rho(\vx,1)\dd\vx + \int_\Omega\phi(\vx,0)\rhoinit(\vx)\dd\vx,
\end{aligned}
\label{eq: rew lagrangian}
\end{equation}
where the second equality is due to integration by part.
The optimal solution solves the saddle point problem 
\begin{equation}
    \min_{\rho\geq0,\vm}\max_{\phi}\calA(\rho,\vm,\phi).
\label{eq: rew saddle pt prob}
\end{equation}
When $L(\vx,\vv)$ is convex in $\vv$,  let the Legendre transformation of $L$ be
\begin{equation}
    H:\Omega\times\bbR^d\to\bbR,(\vx,\vp)\mapsto \sup_{\vv}\left\{-\langle\vp,\vv \rangle - L(\vx,\vv) \right\}.
\end{equation}
Then if $\rho>0$, the optimality condition of \eqref{eq: rew mfg} is 
\begin{equation}
\left\{\begin{aligned}
    &-\partial_t\phi(\vx,t)+H(\vx,\nabla\phi(\vx,t))=\frac{\delta \costevo(\rho)}{\delta\rho}(\vx), 
    \quad \phi(\vx,1)=\frac{\delta \costterm(\rho)}{\delta\rho}(\vx),\\
    &\partial_t\rho(\vx,t)-\nabla\cdot\left(\rho(\vx,t)\partial_{\vp}H(\vx,\nabla\phi(\vx,t))\right)=0,
    \quad \rho(\cdot,0)=\rhoinit.\\
\end{aligned}\right.
\label{eq: rew mfg pde}
\end{equation}
We use this forward-backward PDE system to explore the properties of the inverse problem later.

In this paper, we focus on the following two problems.
\begin{problem}[Crowd motion with obstacle]
\label{eg: rew obs}
A common example comes from crowd motion \cite{ruthotto2020machine}, whose formulation is 
\begin{equation}
\begin{aligned}
    \min_{(\rho,\vm)\in\calC(\rhoinit)}
    \calL(\rho,\vm;\obs):=&\int_0^1\int_\Omega \frac{\|\vm(\vx,t)\|_2^2}{2\rho(\vx,t)}\dd \vx\dd t 
    + \int_0^1 \int_\Omega \rho(\vx,t)\obs(\vx)\dd\vx\dd t \\
    &+\weight_I\int_0^1 \int_\Omega \rho(\vx,t)\log\rho(\vx,t)\dd\vx\dd t
    + \weight_T \int_\Omega \rho(\vx,t)\left( \log\rho(\vx,t)-\log\rhoend(\vx) \right)\dd\vx.
\end{aligned}
\label{eq: rew mfg obs}
\end{equation}
Here the terminal cost is the KL divergence $\costterm(\rho(\cdot,1))=\int_\Omega \rho(\vx,t)\left( \log\rho(\vx,t)-\log\rhoend(\vx) \right)\dd\vx$ which aims to match the terminal density $\rho(\cdot,1)$ to the desired density $\rhoend$. 
The interaction cost contains two parts. The entropy term $\int_\Omega \rho(\vx,t)\log\rho(\vx,t)\dd\vx$ penalizes the aggregation of the density. 
And the obstacle term $\int_\Omega\rho(\vx,t)\obs(\vx)\dd\vx$ penalizes the mass going through the obstacle $\vx$ with larger value of $\obs(\vx)$.
With the same initial density $\rhoinit$, different obstacle functions lead to different Nash Equilibrium.
Assuming that we know everything in the objective function \eqref{eq: rew mfg obs} except the obstacle function $\obs$, we aim to recover $\obs$ from observations of the equilibrium $(\rho,\vm)$.
\end{problem}

\begin{problem}[Non-Euclidean metric]
\label{eg: rew metric}
It is also common to consider mean-field games on spaces with non-Euclidean metrics. If at each $\vx\in\Omega,\Omega\subset\bbR^d$, there is a positive definite matrix $\metric(\vx)\in S_{++}^d$ indicating the metric, then the mean-field game problem takes the form
\begin{equation}
\begin{aligned}
    \min_{(\rho,\vm)\in\calC(\rhoinit)}
    \calL(\rho,\vm;\metric):=&\int_0^1\int_\Omega \frac{\vm(\vx,t)^\top \metric(\vx)\vm(\vx,t)}{2\rho(\vx,t)}\dd \vx\dd t \\
    &+\weight_I\int_0^1 \int_\Omega \rho(\vx,t)\log\rho(\vx,t)\dd\vx\dd t
    + \weight_T \int_\Omega \rho(\vx,t)\left( \log\rho(\vx,t)-\log\rhoend(\vx) \right)\dd\vx.
\end{aligned}
\label{eq: rew mfg metric}
\end{equation}
We also work on solving the metric $\metric$ from the observations of the equilibrium $(\rho,\vm)$, assuming other terms in \eqref{eq: rew mfg metric} are known.

\end{problem}

In summary, we are interested in the mean-field game problem with the objective function
\begin{equation}
\begin{aligned}
    \calL(\rho,\vm;\metric,\obs):=
    &\int_0^1\int_\Omega \frac{\vm(\vx,t)^\top \metric(\vx)\vm(\vx,t)}{2\rho(\vx,t)}\dd \vx\dd t 
    + \int_0^1 \int_\Omega \rho(\vx,t)\obs(\vx)\dd\vx\dd t \\
    &+\weight_I\int_0^1 \int_\Omega \rho(\vx,t)\log\rho(\vx,t)\dd\vx\dd t
    + \weight_T \int_\Omega \rho(\vx,1)\left( \log\rho(\vx,1)-\log\rhoend(\vx) \right)\dd\vx.
\end{aligned}
\label{eq: rew mfg obj}
\end{equation}
We write $\calL(\rho,\vm;\metric)$ when $\obs\equiv0$ and $\calL(\rho,\vm;\obs)$ when $\metric\equiv I_d$.
With $\rho>0$, the optimality condition for the problem 
\begin{equation}
    \min_{(\rho,\vm)\in\calC(\rhoinit)}\calL(\rho,\vm;\metric,\obs),
\label{eq: rew mfg eg}
\end{equation}
is 
\begin{equation}
\left\{\begin{aligned}
    &-\partial_t\phi(\vx,t)+\half(\nabla\phi(\vx,t))^\top(\metric(\vx))^{-1}\nabla\phi(\vx,t)
    =\weight_I(\log\rho(\vx,t)+1)+\obs(\vx), \\
    &\qquad \phi(\vx,1)=\weight_T(\log\rho(\vx,t)-\log\rhoend(\vx)+1),\\
    &\partial_t\rho(\vx,t)-\nabla\cdot\left(\rho(\vx,t)(\metric(\vx))^{-1}\nabla\phi(\vx,t)\right)=0,
    \quad \rho(\cdot,0)=\rhoinit.\\
\end{aligned}\right.
\label{eq: rew mfg pde eg}
\end{equation}

We call the potential mean-field games \eqref{eq: rew mfg}, as well as \eqref{eq: rew mfg obs},\eqref{eq: rew mfg metric}, the forward problem. In this paper, we aim to learn the unknown variables $\obs,\metric$ from one or a set of observations of the Nash Equilibrium $\left\{\left(\rhotilde^n,\vmtilde^n\right)\right\}_{n=1}^N$ that solve the forward problems, and we name this the inverse problem.
Note that the forward problem has a convex objective function and linear constraint, while the optimality condition is nonlinear and nonconvex.
To preserve the nice convex-linear structure of the forward problem, we formulate the inverse mean-field game as a bilevel optimization problem and treat the forward problem as the constraint.

\section{A Bilevel Formulation of Inverse Mean-Field Games}
\label{sec: model}

In this section, we first review the general formulation of a bilevel optimization problem, then provide the bilevel formulation of inverse mean-field games, as well as two concrete inverse problems that we would like to solve in this work. After that, we discretize the model for numerical implementation.

\subsection{Bilevel Formulation}
\label{subsec: blo}
The general formulation of a bilevel optimization problem is
\begin{equation}
\begin{aligned}
    &\min_{\xi\in\Xi} \quad u(\xi):=\calU(\eta^*(\xi);\xi)\\
    &\text{where } \eta^*(\xi)=\argmin_{\eta\in H}\calL(\eta;\xi).
\end{aligned}
\label{eq: general blo}
\end{equation}
Here we consider linear constraint set $H=\{\eta\mid A\eta=c\}$ and convex set $\Xi$, where $A\in\mathbb{R}^{d_c\times d_\eta}, c\in\mathbb{R}^{d_c}$. $d_c<d_\eta$. 
The optimization problem over $\calU$ is referred to as the upper-level problem and that over $\calL$ as the lower-level problem.
We formulate our inverse problems as bilevel optimization problems, with the upper-level objective being a combination of fidelity $\Drho,\Dm$ and regularity $\calR$, and the lower-level problem being the forward problem.
\begin{equation*}
\begin{aligned}
    \min_{\costfuncdyn\in\calC_{\costfuncdyn},\costevo\in\calC_{\costevo}}
    &\quad\calU\left((\rho,\vm),(\rhotilde,\vmtilde) ; \costfuncdyn, \costevo\right)
    :=\left(\Drho(\rho,\rhotilde) + \Dm(\vm,\vmtilde)\right) + \calR(\costfuncdyn,\costevo)\\
    &\text{ s.t. }(\rho,\vm):=\argmin_{(\rho,\vm)\in\calC(\rhoinit)}\calL(\rho,\vm;\costfuncdyn,\costevo).
\end{aligned}
\label{eq: inv mfg general}
\end{equation*}
The dynamic cost $\costfuncdyn$ and interaction cost functional $\costevo$ are the upper-level variables and the density-flux pair $(\rho,\vm)$ is the lower-level variable. For convenience, we choose $\Drho(\rho,\rhotilde) = \frac{1}{2}\int_0^1\int_\Omega (\rho(\vx,t) - \rhotilde(\vx,t))^2\dd\vx\dd t$ and $\Drho(\vm,\vmtilde) = \frac{1}{2}\int_0^1\int_\Omega \|\vm(\vx,t) - \vmtilde(\vx,t)\|_2^2\dd\vx\dd t$

We formulate the inverse problems of problem \ref{eg: rew obs} and \ref{eg: rew metric} as follows.
\begin{problem}[The inverse problem of crowd motion (problem \ref{eg: rew obs})]
\label{eg: blo obs}
Let the regularity be $\calR(\obs)=0$. The inverse problem of \eqref{eq: rew mfg obj} is
\begin{equation}
\begin{aligned}
    \min_{\obs\in\Cobs} 
    &\quad\Drho(\rho,\rhotilde) + \Dm(\vm,\vmtilde)\\
    &\text{s.t. }(\rho,\vm):=\argmin_{(\rho,\vm)\in\calC(\rhoinit)}\calL(\rho,\vm;\obs).
\end{aligned}
\label{eq: inv mfg obs}
\end{equation}
Here $(\rhotilde,\vmtilde)=\argmin_{(\rho,\vm)\in\calC(\rhoinit)}\calL(\rho,\vm;\obstilde)$ are the observed data with ground truth $\obstilde$.
Notice that for any constant $c\in\bbR$, if $(\rhotilde,\vmtilde)$ minimizes $\calL(\rho,\vm;\obstilde)$, then $(\rhotilde,\vmtilde)$ also minimizes $\calL(\rho,\vm;\obstilde+c)$. To remove the ambiguity, we restrict our focus to obstacle functions with zero integral, i.e.
\begin{equation}
\begin{aligned}
    \Cobs:=\left\{\obs:\int_{\Omega}\obs(\vx)\dd\vx=0\right\}.
\end{aligned}
\label{eq: cstr set obs}
\end{equation}
Ideally, we expect $\proj_{\Cobs}(\obstilde)$ to be the unique minimizer of the bilevel problem \eqref{eq: inv mfg obs}. We prove this unique identifiability property for the discretization of \eqref{eq: inv mfg obs} in section \ref{sec: proof}.
\end{problem}

\begin{problem}[The inverse problem of unknown metric (problem \ref{eg: rew metric})]
\label{eg: blo metric}
Similarly, we have the bilevel formulation to recover the metric $\mettilde$ from the data $\displaystyle(\rhotilde,\vmtilde)=\argmin_{(\rho,\vm)\in\calC(\rhoinit)}\calL(\rho,\vm;\mettilde)$.
\begin{equation}
\begin{aligned}
    \min_{\metric\in\Cmetric} 
    &\quad \Drho(\rho,\rhotilde) + \Dm(\vm,\vmtilde) + \calR(\metric)\\
    &\text{s.t. }(\rho,\vm):=\argmin_{(\rho,\vm)\in\calC(\rhoinit)}\calL(\rho,\vm;\metric).
\end{aligned}
\label{eq: inv mfg metric}
\end{equation}
To make sure that $\metric$ induces a metric on $\Omega$, we set the constraint of $\metric$ as
\begin{equation}
\begin{aligned}
    \Cmetric:=\{\metric:\Omega\rightarrow S^{d}_{++}:\metric(\vx)\text{ are positive definite matrices},\forall \vx\in\Omega\}.
\end{aligned}
\label{eq: cstr set metric}
\end{equation}
\end{problem}

For one observation, if the density is zero in an open set, it means almost no players pass the region and it is impossible to obtain the exact information in that region. However multiple observations may complement the missing information, and therefore it is meaningful to consider the following inverse MFG with multiple observations.
\begin{problem}[The inverse problem of unknown metric (problem \ref{eg: rew metric}) with multiple observations]
\label{eg: blo metric multidata}

Suppose that we have multiple observations of the Nash Equilibrium with a given $\mettilde$ from different initial densities $\rhoinit^n,n=1,\cdots,N$, i.e. $\displaystyle(\rhotilde^n,\vmtilde^n)=\argmin_{(\rho,\vm)\in\calC(\rhoinit^n)}\calL(\rho,\vm;\mettilde)$ for $n=1,2,\cdots, N$. Then we can solve the following bilevel optimization problem to recover the true metric.
\begin{equation}
\begin{aligned}
    \min_{\metric\in\Cmetric} 
    &\quad\sum_{n=1}^N \left(\Drho(\rho^n,\rhotilde^n) + \Dm(\vm^n,\vmtilde^n)\right) + \calR(\metric)\\
    &\text{s.t. }
    \{(\rho^n,\vm^n)\}_{n=1}^N:= \argmin_{(\rho_n,\vm_n)\in\calC(\rhoinit^n)}\sum_{n=1}^N\calL(\rho_n,\vm_n;\metric).
\end{aligned}
\label{eq: inv mfg metric multidata}
\end{equation}
The lower-level is equivalent to a concatenation of $N$ forward problems since $(\rho^n,\vm^n)$ are independent.
\end{problem}

\subsection{Discretization}
\label{subsec: disct}
\begin{figure}[h]
    \centering
    \includegraphics[width=8cm]{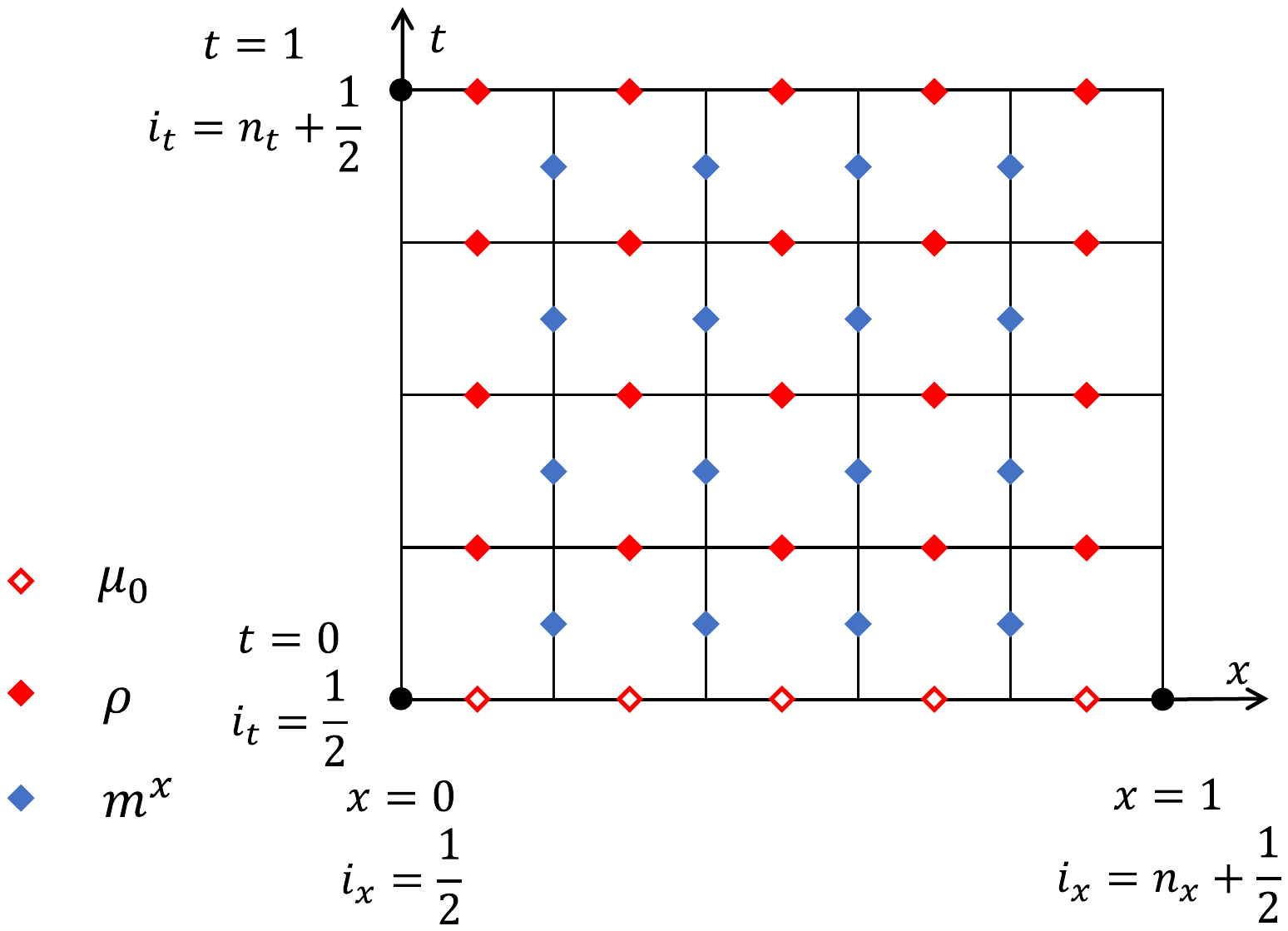}
    \includegraphics[width=8cm]{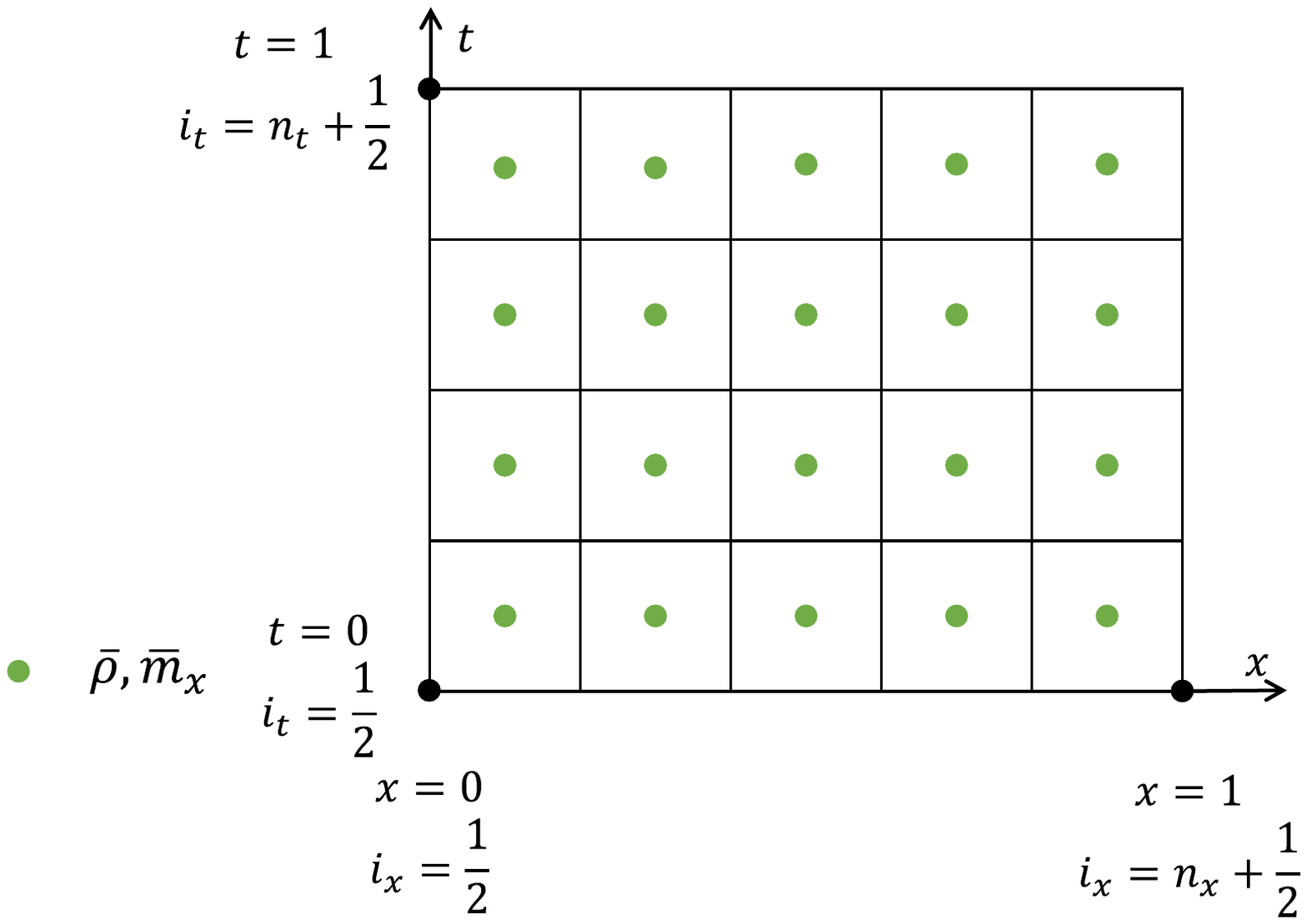}
    \caption{Illustrations of the staggered (left) and central (right) grids.}
\label{fig: grid}
\end{figure}

We conduct numerical experiments on $\bbR^d$ with $d=1,2$. 
Taking $d=2$ as an example, we let $\Omega = [0,1]\times[0,1]$ and the space-time joint domain be $[0,1]^3$, and we write $\vm = (\mx,\my)$.
We follow the discretization in \cite{yu2021fast}, with which the discrete optimizer is consistent with the continuous optimizer under certain regularity conditions.
To be precise, we equally divide $[0,1]$ into $\nx,\ny,\nt$ parts, and each cube is of size $\dx\dy\dt$, with $\dx=\frac{1}{\nx},\dy=\frac{1}{\ny},\dt=\frac{1}{\nt}.$
Let $x_i=(i-\half)\dx,y_i=(i-\half)\dy,t_i=(i-\half)\dt$, and $(f)_{\id}$ approximates function $f$
on points $(x_{\idx},y_{\idy},t_{\idt})$. Similarly, $(f)_{\idxy}\approx f(x_\idx,y_\idy)$.
We define $\Grho,\Gmx$ and $\Gmy$ as the sets of grid point indices on $t$-, $x$- and $y$-staggered grids, respectively, where
\begin{equation*}
\begin{aligned}
    &\Grho:=\{(\idtp):\idx=1,\cdots,\nx,\idy=1,\cdots,\ny,\idt=1,\cdots,\nt\},\\
    &\Gmx:=\{(\idxp):\idx=1,\cdots,\nx-1,\idy=1,\cdots,\ny,\idt=1,\cdots,\nt\},\\
    &\Gmy:=\{(\idyp):\idx=1,\cdots,\nx,\idy=1,\cdots,\ny-1,\idt=1,\cdots,\nt\}.
\end{aligned}
\end{equation*}
Then we approximate the function $\rho,\mx$ and $\my$ on $t$-, $x$- and $y$-staggered grids by $\rho_{\Grho},\mx_{\Gmx}$ and $\my_{\Gmy}$, respectively,
i.e.
$\rho_{\Grho}:=\left\{(\rho)_\vi\right\}_{\vi\in\Grho}\in\bbR^{\nx\ny\nt}$, 
$\mx_{\Gmx}=\left\{(\mx)_\vi\right\}_{\vi\in\Gmx}\in\bbR^{(\nx-1)\ny\nt}$ and 
$\my_{\Gmy}=\left\{(\my)_\vi\right\}_{\vi\in\Gmy}\in\bbR^{\nx(\ny-1)\nt}$. 
We denote $\Gm:=\Gmx\times\Gmy$ as the concatenation of $\Gmx,\Gmy$ and $\vm_{\Gm}:=\{\mx_{\Gmx},\my_{\Gmy}\}$ as the concatenation of $\mx_{\Gmx},\my_{\Gmy}$.
We will omit the under scripts of grids wherever there is no ambiguity according to context. 
The left part of Figure \ref{fig: grid} illustrates the staggered grids and the corresponding $\rho,\mx$ for $d=1$.

We define the inner products on the staggered grids as
\begin{equation*}
\begin{aligned}
    &\langle\rho_1,\rho_2\rangle_{\Grho} :=\dx\dy\dt\sum_{\vi\in\Grho} (\rho_1)_\vi(\rho_2)_\vi,\\
    &\langle\mx_1,\mx_2\rangle_{\Gmx}:=\dx\dy\dt\sum_{\vi\in\Gmx}(\mx_1)_\vi(\mx_2)_{\vi},\\
    &\langle\my_1,\my_2\rangle_{\Gmy}:=\dx\dy\dt\sum_{\vi\in\Gmy}(\my_1)_\vi(\my_2)_{\vi},
\end{aligned}
\end{equation*}
and denote their induced norm as $\|\cdot\|_{\Grho},\|\cdot\|_{\Gmx}$ and $\|\cdot\|_{\Gmy}$.
Based on these, we approximate the discrepancy between lower-level minimizer and observed data $\Drho,\Dm$ by the sum of element-wise differences on grids $\DGrho,\DGm$, 
\begin{equation}
\begin{aligned}
    \DGrho(\rho,\rhotilde):=&\half\|\rho-\rhotilde\|^2_{\Grho}\\
    \DGm(\vm,\vmtilde):=&\half\|\mx-\mxtilde\|^2_{\Gmx} + \half\|\my-\mytilde\|^2_{\Gmy}.
\end{aligned}
\label{eq: disct upper level obj}
\end{equation}

To compute the objective function, we consider the central grid (see the right plot in Figure \ref{fig: grid})
$$\Gphi:=\{(\idx,\idy,\idt):\idx=1,\cdots,\nx,\idy=1,\cdots,\ny,\idt=1,\cdots,\nt\}.$$
We define the inner product and induced norm on the central grid similarly and denote them as $\langle\cdot,\cdot\rangle_{\Gphi}$ and $\|\cdot\|_{\Gphi}$.
With the interpolation operators, $\rhobar=\It(\rho;\rhoinit)$, $\mxbar=\Ix(\mx),\mybar=\Iy(\my)$ meet on the central grid points:
\begin{eqnarray*}
&(\rhobar)_\id=(\It(\rho;\rhoinit))_\id&:=\begin{cases}
\half\left((\rhoinit)_\idxy+(\rho)_\idtp\right), & \idt=1,\\
\half\left((\rho)_\idtm+(\rho)_\idtp\right), & \idt=2,\cdots,\nt.
\end{cases}\\
&(\mxbar)_\id=(\Ix(\mx))_\id&:=\begin{cases}
\half(\mx)_\idxp, & \idx=1,\\
\half\left((\mx)_\idxm+(\mx)_\idxp\right), & \idx=2,\cdots,\nx-1,\\
\half(\mx)_\idxm, & \idx=\nx.
\end{cases}\\
&(\mybar)_\id=(\Iy(\my))_\id&:=\begin{cases}
\half(\my)_\idyp, & \idy=1,\\
\half\left((\my)_\idym+(\my)_\idyp\right), & \idy=2,\cdots,\ny-1,\\
\half(\my)_\idym, & \idy=\ny.
\end{cases}
\end{eqnarray*}
Here, the definition of $\mxbar$ on $\idx=1,\nx$ and $\mybar$ on $\idy=1,\ny$ are consistent with the zero-flux boundary condition in the continuous setting. 
The objective functions of the forward problem can therefore be approximated by
\begin{equation}
\begin{aligned}
    \calL_\calG(\rho,\vm;\metric,\obs)
    &:=\dx\dy\dt\sum_{\vi\in\Gphi}
    \left(\frac{ (\vmbar)_\vi^\top (\metric)_{\idxy} (\vmbar)_{\vi} }{2\left(\rhobar\right)_\vi}
    +\weight_I(\rhobar)_\vi\log((\rhobar)_\vi)\right)
    +\dx\dy\dt\sum_{\vi\in\Grho}(\rho)_\vi(b)_\idxy\\
    &+\weight_T\dx\dy\sum_{\idx=1}^\nx\sum_{\idy=1}^\ny(\rho)_{\idx,\idy,\nt+\half}\left(\log(\rho)_{\idx,\idy,\nt+\half}-\log(\rhoend)_{\idxy}\right).
\end{aligned}
\label{eq: disct forward mfg obj}
\end{equation}
where $\vm=\{\mx,\my\}$, $(\vmbar)_\id^\top:=((\mxbar)_{\id},(\mybar)_{\id})$ and the subscript of $\calL_\calG$ indicates the cost is defined on the discrete space.
Similar to the continuous notation, we write $\calL_\calG(\rho,\vm;\metric)$ when $\obs=\zero$ and $\calL_\calG(\rho,\vm;\obs)$ when $\metric\equiv1 ~ (d=1)$ or $\metric\equiv I_2 ~ (d=2)$.

With this discretization,  $\calL_\calG(\rho,\vm;\metric,\obs)$ preserves the following properties on $(\rho,\vm)$ from the continuous setting.

\begin{lemma}
\label{lem: convexity}
For $\calL_\calG(\rho,\vm;\metric,\obs)$ defined on $(\rho_{\Grho},\mx_{\Gmx},\my_{\Gmy})\in\bbR^{\nx\ny\nt}\times\bbR^{(\nx-1)\ny\nt}\times\bbR^{\nx(\ny-1)\nt}$ with $\displaystyle\min_{\vi\in\Grho}(\rho)_{\vi}>0$, the following statements hold:
\begin{enumerate}
    \item If $\weight_I,\weight_T\geq0$ and $\metric_{\idx,\idy}$ is positive definite for all $\idx,\idy$, 
    then $\calL_\calG(\rho,\vm;\metric,\obs)$ is convex in $\rho,\vm$.
    \item In addition to 1, if we restrict the domain to 
    $\rho$ with $\displaystyle\min_{\vi\in\Grho}(\rho)_{\vi}\geq\underline{c}_\rho>0$, 
    $\mx$ with $\displaystyle\max_{\vi\in\Gmx}(|\mx_\vi|)\leq\overline{c}_m$, 
    and $\my$ with $\displaystyle\max_{\vi\in\Gmy}(|\my_\vi|)\leq\overline{c}_m$, 
    then $\calL_\calG(\rho,\vm;\metric,\obs)$ is Lipschitz smooth in $\rho,\vm$.
    \item In addition to 1,2, if we further restrict the domain to 
    $\rho\in\bbR^{\nx\ny\nt}$ with $\displaystyle\max_{\vi\in\Grho}(\rho)_{\vi}\leq\overline{c}_\rho$, 
    then for any $\weight_I,\weight_T>0$, $\calL_\calG(\rho,\vm;\metric,\obs)$ is strongly convex in $\rho,\vm$. 
\end{enumerate}
\end{lemma}
We postpone the proof of the lemma in section \ref{sec: proof} for better readability.
The Lipschitz smoothness and strong convexity of the lower-level objective are important to guarantee the convergence of our alternating gradient algorithm, as detailed later in section \ref{sec: alg}.

Following the nature of the staggered grid, we choose a central difference scheme to approximate the differential operators.
\begin{eqnarray*}
&(\Dt(\rho;\rhoinit))_\id&:=\begin{cases}
\frac{1}{\dt}\left((\rho)_\idtp-(\rhoinit)_\idxy\right), & \idt=1,\\
\frac{1}{\dt}\left((\rho)_\idtp-(\rho)_\idtm\right), & \idt=2,\cdots,\nt.
\end{cases}\\
&(\Dx(\mx))_\id&:=\begin{cases}
\frac{1}{\dx}(\mx)_\idxp, & \idx=1,\\
\frac{1}{\dx}\left((\mx)_\idxp-(\mx)_\idxm\right), & \idx=2,\cdots,\nx-1,\\
-\frac{1}{\dx}(\mx)_\idxm, & \idx=\nx.
\end{cases}\\
&(\Dy(\my))_\id&:=\begin{cases}
\frac{1}{\dy}(\my)_\idyp, & \idy=1,\\
\frac{1}{\dy}\left((\my)_\idyp-(\my)_\idym\right), & \idy=2,\cdots,\ny-1,\\
-\frac{1}{\dy}(\my)_\idym, & \idy=\ny.
\end{cases}
\end{eqnarray*}
Again, the definitions of $\Dx,\Dy$ on $\idx=1,\nx,\idy=1,\ny$, respectively, are consistent with the zero-flux boundary condition.
The discrete constraint set is
\begin{equation}
\begin{aligned}
    \calC_{\calG}(\rhoinit):=\left\{ (\rho,\vm):\Dt(\rho;\rhoinit)+\Dx(\mx)+\Dy(\my)=\zero \right\}.
\end{aligned}
\label{eq: disct cstr set cteqn}
\end{equation}

Based on the above notations, we restate the inverse problems \ref{eg: blo obs} and \ref{eg: blo metric} in the discretized space. We intentionally write down the problems for more general cases with multiple pairs of training data as they will reduce to the case with a single pair data by chosing $N=1$. 

\begin{problem}[The discretization of the inverse crowd motion problem \ref{eg: blo obs}]
\label{eg: disct blo obs}
The discretization of \eqref{eq: inv mfg obs} has the following formulation
\begin{equation}
\begin{aligned}
    \min_{\obs\in\calC_{\calG,\obs}} 
    &\quad\sum_{n=1}^N \left(\DGrho(\rho^n,\rhotilde^n) + \DGm(\vm^n,\vmtilde^n)\right)\\
    &\text{s.t. }(\rho^n,\vm^n):=\argmin_{(\rho,\vm)\in\calC_{\calG}(\rhoinit^n)}\calL_{\calG}(\rho,\vm;\obs), n=1,2,\cdots,N,
\end{aligned}
\label{eq: disct inv mfg obs}
\end{equation}
where $(\rhotilde^n,\vmtilde^n)=\argmin_{(\rho,\vm)\in\calC_{\calG}(\rhoinit^n)}\calL_\calG(\rho,\vm;\obstilde)$ are the observed data and 
\begin{equation}
\begin{aligned}
    \calC_{\calG,\obs}:=\left\{\obs:\sum_{\idx,\idy}^{\nx,\ny}(\obs)_{\idxy}=0\right\}.
\end{aligned}
\label{eq: disct cstr set obs}
\end{equation}

\end{problem}

\begin{problem}[The discretization of the inverse metric problem \ref{eg: blo metric}]
\label{eg: disct blo metric}
Similarly, given the data $\displaystyle(\rhotilde^n,\vmtilde^n)=\argmin_{(\rho,\vm)\in\calC(\rhoinit^n)}\calL(\rho,\vm;\mettilde)$, we implement algorithms to solve
\begin{equation}
\begin{aligned}
    \min_{\metric\in\calC_{\calG,\metric}} 
    &\quad\sum_{n=1}^N \left(\DGrho(\rho^n,\rhotilde^n) + \DGm(\vm^n,\vmtilde^n)\right) + \calR_{\calG}(\metric)\\
    &\text{s.t. }(\rho^n,\vm^n):=\argmin_{(\rho,\vm)\in\calC_\calG(\rhoinit^n)}\calL_{\calG}(\rho,\vm;\metric),n=1,2,\cdots,N,
\end{aligned}
\label{eq: disct inv mfg metric}
\end{equation}
with the constraint of $\metric$ being
\begin{equation}
\begin{aligned}
    \calC_{\calG,\metric}:=\{\metric:(\metric)_{\idxy}\in\bbR^{d\times d}\text{ are positive definite matrices}, \idx=1,\cdots,\nx,\idy=1,\cdots,\ny\}.
\end{aligned}
\label{eq: disct cstr set metric}
\end{equation}
\end{problem}

\subsection{Regularity and Unique Identifiability of the Inverse Problems}

At the end of this section, we state the regularity of the inverse problem \ref{eg: disct blo obs} and \ref{eg: disct blo metric} and the unique identifiability of the inverse crowd motion problem \ref{eg: disct blo obs}.

The regularity and unique identifiability of the inverse problem rely on the KKT system of the discretized forward problem
\begin{equation}
\begin{aligned}
    \min_{(\rho,\vm)\in\calC(\rhoinit)}\calL_{\calG}(\rho,\vm;\metric,\obs).
\end{aligned}
\label{eq: disct mfg eg}
\end{equation}
To write the KKT system in a concise way, we define the adjoint operators of $\Ix,\Iy,\It$ for any $\phi=\phi_{\Gphi}$ on the central grid as
\begin{eqnarray*}
&(\It^*(\phi))_\idtp&:=\begin{cases}
\half\left((\phi)_\id+(\phi)_\idtpp\right), & \idt=1,\cdots,\nt-1,\\
\half(\phi)_\id, & \idt=\nt.
\end{cases}\\
&(\Ix^*(\phi))_\idxp&:=
\half\left((\phi)_\id+(\phi)_\idxpp\right), \idx=1,\cdots,\nx-1.\\
&(\Iy^*(\phi))_\idyp&:=
\half\left((\phi)_\id+(\phi)_\idypp\right), \idy=1,\cdots,\ny-1.
\end{eqnarray*}
and the adjoint operators of $\Dx,\Dy,\Dt$ as
\begin{eqnarray*}
&(\Dt^*(\phi))_\idtp&:=\begin{cases}
-\frac{1}{\dt}\left((\phi)_\idtpp-(\phi)_\id\right), & \idt=1,\cdots,\nt-1\\
\frac{1}{\dt}(\phi)_\id, & \idt=\nt.
\end{cases}\\
&(\Dx^*(\phi))_\idxp&:=
\frac{1}{\dx}\left((\phi)_\idxpp-(\phi)_\id\right),  \idx=1,\cdots,\nx-1\\
&(\Dy^*(\phi))_\idyp&:=
\frac{1}{\dy}\left((\my)_\idypp-(\my)_\id\right),  \idy=1,\cdots,\ny-1.
\end{eqnarray*}
The adjoint relation in the discretized space holds based on the definitions. To be precise, for the interpolation operators, we have
\begin{equation*}
\begin{aligned}
\langle \It(\rho;\rhoinit),\phi \rangle_{\Gphi}
= &\langle \rho,\It^*(\phi) \rangle_{\Grho} 
+\half\sum_{\idx=1}^\nx\sum_{\idy=1}^\ny(\rhoinit)_\idxy(\phi)_{\idx,\idy,1}.\\
\langle \Ix(\mx),\phi \rangle_{\Gphi}
= &\langle \mx,\Ix^*(\phi), \rangle_{\Gmx}.\\
\langle \Iy(\my),\phi \rangle_{\Gmy}
= &\langle \my,\Iy^*(\phi),\rangle_{\Gmy}.\\
\end{aligned}
\end{equation*}
and for differential operators, we have
\begin{equation*}
\begin{aligned}
\langle \Dt(\rho;\rhoinit),\phi \rangle_{\Gphi}
= &\langle \rho,\Dt^*(\phi) \rangle_{\Grho} 
-\frac{1}{\dt}\sum_{\idx=1}^\nx\sum_{\idy=1}^\ny(\rhoinit)_\idxy(\phi)_{\idx,\idy,1}.\\
\langle \Dx(\mx),\phi \rangle_{\Gphi}
= &\langle \mx,\Dx^*(\phi) \rangle_{\Gmx}.\\
\langle \Dy(\my),\phi \rangle_{\Gphi}
= &\langle \my,\Dy^*(\phi) \rangle_{\Gmy}.\\
\end{aligned}
\end{equation*}

With the adjoint operators, we define the $\calK$ operators as following to describe the optimality condition for the forward problem,
\begin{equation}\left\{
\begin{aligned}
&\vi\in\Grho,\idt=1,\cdots,\nt-1,\\
&(\calK_{\rho}(\rho,\vm,\phi;\metric,\obs))_\vi
    :=-\left( \Dt^*(\phi) \right)_\vi 
    +\left( \It^*\left( -\frac{(\vmbar)^\top \metric \vmbar}{2\rhobar^2}+\weight_I(\log(\rhobar)+1) \right) \right)_\vi
    +\obs_\idxy,\\ 
&\vi\in\Grho,\it=\nt,\\
&(\calK_{\rho}(\rho,\vm,\phi;\metric,\obs))_\vi
    :=-\left( \Dt^*(\phi) \right)_\vi
    +\left( \It^*\left( -\frac{(\vmbar)^\top \metric \vmbar}{2\rhobar^2}+\weight_I(\log(\rhobar)+1) \right) \right)_\vi
    +\obs_\idxy\\
    &\qquad\qquad +\frac{\weight_T}{\dt}(\log(\rho_\vi)-\log((\rhoend)_{\idxy})+1),\\
&\vi\in\Gmx,\quad(\calK_{\mx}(\rho,\vm,\phi;\metric,\obs))_\vi:=
    -\left( \Dx^*(\phi) \right)_\vi + \left( \Ix^*\left(\frac{\metric_{xx}\mxbar+\metric_{xy}\mybar}{\rhobar}\right) \right)_\vi,\\
&\vi\in\Gmy,\quad(\calK_{\my}(\rho,\vm,\phi;\metric,\obs))_\vi:=
    -\left( \Dy^*(\phi) \right)_\vi + \left( \Iy^*\left(\frac{\metric_{xy}\mxbar+\metric_{yy}\mybar}{\rhobar}\right) \right)_\vi,\\    
&\vi\in\Gphi,\quad(\calK_{\phi}(\rho,\vm,\phi;\metric,\obs))_\vi:= 
    \left(\Dt(\rho;\rhoinit)+\Dx(\mx)+\Dy(\my)\right)_\vi.
\end{aligned}\right.
\label{eq: definition of calK}
\end{equation}
$\calK_\rho,\calK_\mx,\calK_\my,\calK_\phi$ are obtained by taking gradients on the Lagrangian of the forward problem \eqref{eq: disct mfg eg}.
By viewing $\rho,\vm,\phi,\obs$ and $\calK_\rho,\calK_\mx,\calK_\my,\calK_\phi$ as long vectors and denoting $\calK:=(\calK_{\rho},\calK_{\mx},\calK_{\my},\calK_{\phi})^\top$, we define a function $\calK:\bbR^{d_l}\times\bbR^{d_u}\to\bbR^{d_l}$ with $d_u=(\frac{d(d+1)}{2}+1)\nx\ny$ corresponding to the dimension of $(\metric,\obs)$ and $d_l=\nx\ny\nt + (\nx-1)\ny\nt + \nx(\ny-1)\nt + \nx\ny\nt$ to the dimension of $\rho,\mx,\my,\phi$.
Since the constraint is linear, the optimizer of \eqref{eq: disct mfg eg} satisfies the KKT condition. The formal statement is the following.
\begin{lemma}
    If $(\rhotilde,\vmtilde)\in\calC(\rhoinit)$ is a minimizer of $\calL_\calG(\rho,\vm;\mettilde,\obstilde)$, and $\displaystyle\min_{\vi\in\Grho}\left\{\rhotilde_\vi\right\}>0$, 
    then there exists $\phitilde\in\bbR^{\nx\ny\nt}$ such that
    \begin{equation}
    \calK(\rhotilde,\vmtilde,\phitilde;\mettilde,\obstilde)=\zero.
    \label{eq: kkt obs}
    \end{equation} 
\label{lem: kkt -> zero function}
\end{lemma}

With the discrete PDE description of the Nash Equilibrium, we state the regularity result for the inverse problem \ref{eg: disct blo obs} and \ref{eg: disct blo metric}.
\begin{theorem}[Regularity]
\label{thm: differentiability}
Assume that $(\rhotilde,\vmtilde)$ is the Nash Equilibrium given the metric $\mettilde$, obstacle function $\obstilde$ and $\weight_I>0,\weight_T>0$, i.e. \eqref{eq: disct mfg eg} holds,
and that $\displaystyle\min_{\vi\in\Grho}\rhotilde_{\vi}>0$, 
then there exists $r_u>0$ and a radius $r_u$ open ball $B_{r_u}(\mettilde,\obstilde)$ centered at $(\mettilde,\obstilde)$, and a mapping $\mathscr{T}$ defined on $B_{r_u}(\mettilde,\obstilde)$ satisfying the following 
\begin{itemize}
    \item For any $(\metric,\obs)\in B_{r_u}(\mettilde,\obstilde)$, there exist a unique $(\rho,\vm,\phi)=\mathscr{T}(\metric,\obs)\in B_{r_l}(\rhotilde,\vmtilde,\phitilde)$, a radius $r_l$ open ball centered at $(\rhotilde,\vmtilde,\phitilde)$, such that $(\rho,\vm,\phi)$ solves the forward problem with $\calL_{\calG}(\rho,\vm;\metric,\obs)$.
    \item $\mathscr{T}(\mettilde,\obstilde)=(\rhotilde,\vmtilde,\phitilde)$, $\mathscr{T}$ is of class $C^1$ and 
    \begin{equation}
        D \mathscr{T}(\metric,\obs)=-\left(D_{\rho,\vm,\phi}\calK(\mathscr{T}(\metric,\obs);\metric,\obs))\right)^{-1} \left( D_{\obs}\calK(\mathscr{T}(\metric,\obs);\metric,\obs) \right), \text{ for all } (\metric,\obs)\in B_{r_{u}}(\mettilde,\obstilde).
    \label{eq: D_g,b rho,m,phi}
    \end{equation}   
\end{itemize}

\end{theorem}

In addition, we have the unique identifiability of the inverse crowd motion problem because the lower-level objective has a simple dependence on the obstacle $\obs$. To be concrete, by solving the inverse crowd motion problem \ref{eg: disct blo obs}, we uniquely recover the ground truth obstacle $\obstilde$ up to a constant from only one good observation of the Nash Equilibrium.

\begin{theorem}[Unique identifiability]
\label{thm: unique identifiablity}
Assume that $(\rhotilde,\vmtilde)$ is the Nash Equilibrium given the obstacle function $\obstilde$, i.e. 
    \begin{equation}
        (\rhotilde,\vmtilde):=\argmin_{(\rho,\vm)\in\calC_\calG(\rhoinit)}\calL_\calG(\rho,\vm;\obstilde),
    \label{eq: disct mfg forward obs prob}
    \end{equation}
and that $\displaystyle\min_{\vi\in\Grho}\rhotilde_{\vi}>0$, 
then any minimizer $\obs$ of the bilevel minimization problem 
\begin{equation}
\begin{aligned}
    \min_{\obs} 
    &\quad \DGrho(\rho,\rhotilde) + \DGm(\vm,\vmtilde)\\
    &\text{s.t. }(\rho,\vm):=\argmin_{(\rho,\vm)\in\calC_{\calG}(\rhoinit)}\calL_{\calG}(\rho,\vm;\obs),
\end{aligned}
\label{eq: inv mfg obs one data}
\end{equation}
has the form $\obs=\obstilde+c$ where $c\in\bbR$ is a constant.
This implies that $\obstilde$ is the unique minimizer of the bilevel minimization problem \eqref{eq: inv mfg obs one data} up to a constant.

\end{theorem}

The proofs are postponed to section \ref{sec: proof}. We close this section with some remarks on the theorems.

\begin{remark}[Numerical stability]
\label{rem: robustness reasoning}
    While the unique identifiability Theorem \ref{thm: unique identifiablity} holds without the entropy term and the regularity Theorem \ref{thm: differentiability}, we emphasize that the entropy term and regularity theorem are meaningful for studying the numerical stability of the inverse problem.
    In fact, the entropy term guarantees the strong convexity of the objective function and thus the uniqueness of the forward problem. And it is important for the regularity Theorem \ref{thm: differentiability} to hold. 
    The regularity argument states the differentiability of the forward optimizer with respect to the metric $\metric$ and the obstacle $\obs$ and reveals the rate of change. 
    According to Theorem \ref{thm: unique identifiablity}, if the smallest singular value of $D\mathscr{T}(\metric,\obs)$ is large, then a small perturbation to $(\rhotilde,\vmtilde)$ can still give a reasonable approximation of the ground truth $\mettilde,\obstilde$.     
    It is also worth noting that when $\min_\vi\rhotilde_\vi$ is close to 0, the condition number of the Jacobian matrix $D_{\rho,\vm,\phi}\calK(\rhotilde,\vmtilde,\phitilde;\obs)$ in \eqref{eq: D_g,b rho,m,phi} can be extremely large. Therefore the Jacobian matrix $D_\obs(\rho,\vm,\phi)$ is close to singular, and the observation error may cause a failure to recover the ground truth obstacle.
\end{remark}

\begin{remark}[Unique identifiability in the function space]
\label{rem: cts unique identifiablity}
Theorem \ref{thm: unique identifiablity} establishes the unique identifiability of the obstacle $\obs_{\calG}\in\bbR^{\nx\ny}$ in the discretized finite-dimensional space.
To prove the parallel result for the obstacle function $\obs:\Omega\to\bbR$ in the infinite-dimension space, it is subtle to choose the function space for $\obs,\rho,\vm$, and $\phi$.
The function space is expected to be large enough to guarantee the existence of the lower-level optimizers $\rho^*(\obs),\vm^*(\obs)$ for different $\obs$, and to guarantee the existence of the bilevel problem optimizer $\obs^*$. 
Meanwhile, the functions in the space require enough regularity for $\rho^*(\obs),\vm^*(\obs)$ to be differentiable with respect to $\obs$.
This is out of the scope of this paper. We refer interested readers to \cite{liu2022inverse,liu2022inversebd,ren2023unique} for efforts in studying the unique identifiability in the infinite-dimensional space, where infinity many pairs of training data are required. 
\end{remark}

\begin{remark}[Unique identifiability of the unknown metric]
\label{rem: recover metric}
To establish the local unique identifiability of the metric as a corollary of the stability Theorem \ref{thm: differentiability}, we need $D_{\metric}\calK(\rhotilde,\vmtilde,\phitilde;\metric)$ to have full rank. 
However, for 1D metric, the rank of $D_{\metric}\calK(\rhotilde,\vmtilde,\phitilde;\metric)$ depends on the data $\rhotilde,\vmtilde,\phitilde$, which is different from $D_{\obs}\calK(\rhotilde,\vmtilde,\phitilde;\obs)$ being a constant. 
Therefore, we may not uniquely recover the metric from the data. 
Besides the degenerated rank, while uniquely identifying $\metric$ requires the knowledge of $\phitilde$, we do not have $\phitilde$ in our problem setting and this can also cause non-uniqueness of the inverse problem. 
By experiments in \cite{ding2022mean}, the lack of information on $\phitilde$ can be overcome by giving partial true information on the metric and incorporating regularity terms in the upper-level objective.
For 2D metric, if we view $\metric_{xx},\metric_{xy},\metric_{yy}$ as independent variables, then $D_{\metric}\calK(\rhotilde,\vmtilde,\phitilde;\metric)$ is not a full-rank matrix and theoretically there is no hope to uniquely recover the ground truth metric. 
If the metric $\metric_{\vi}\in S^2_{++}$ has intrinsic structures such that the number of variables to determine the metric is $\nx\ny$ instead of $3\nx\ny$, numerically we recover the ground truth with a low error as shown by the numerical experiment in section \ref{subsec: num 2d met}.
The numerical experiment in section \ref{sssec: multidata} also shows that another way to resolve the ambiguity is to have multiple observations for more complete information in the region.
\end{remark}

\section{Alternating Gradient Method}
\label{sec: alg}

In this section, we present the alternating gradient method (AGM) to solve the general bilevel optimization problem \eqref{eq: general blo}, as well as two inverse mean-field game problems \ref{eg: disct blo obs} and \ref{eg: disct blo metric}. 

\subsection{Preliminary on AGM for Bilevel Optimization}
The idea of the AGM is iteratively conducting gradient descent on the lower-level variable and the upper-level variable.
To illustrate our algorithm, we first consider the following unconstrained bilevel   problem 
\begin{equation}
\begin{aligned}
    &\min_{\xi\in\bbR^{d_u}} \quad u(\xi):=\calU(\eta^*(\xi);\xi)\\
    &\text{where } \eta^*(\xi)=\argmin_{\eta\in\bbR^{d_l}}\calL(\eta;\xi).
\end{aligned}
\label{eq: general blo uncstr}
\end{equation}
The computation of the lower-level gradient is straightforward. 
To obtain the upper-level gradient, we assume that $\calU,\calL$ are differentiable and denote the gradient operator with respect to their first and second entries as $\nabla_\eta,\nabla_\xi$.
If for any given $\xi$, there exists a unique $\eta^*(\xi)$ solving the lower-level optimization problem and the function mapping $\xi$ to its corresponding minimizer $\eta^*(\xi)$ is differentiable, then by chain rule, we have
\begin{equation}
    \nabla u(\xi) = \nabla_{\xi}\eta^*(\xi)^\top\nabla_\eta \calU(\eta^*(\xi);\xi) 
                    +\nabla_\xi\calU(\eta^*(\xi);\xi),
\label{eq: ul grad chain rule}
\end{equation}
with $\nabla_{\xi}\eta^*(\xi)=(\partial_{\xi_1}\eta^*(\xi),\ldots,\partial_{\xi_{d_u}}\eta^*(\xi))\in\bbR^{d_l\times d_u}$ being the Jacobian matrix of $\eta^*$.
We clarify that here $\nabla_\xi\calU(\eta^*(\xi);\xi)$ is the gradient of $\calU$ with respect to its second entry evaluated at $(\eta^*(\xi);\xi)$ without considering the dependence of $\eta^*$ on $\xi$.
Therefore $\nabla_\eta \calU(\eta^*(\xi);\xi)$ and $\nabla_\xi\calU(\eta^*(\xi);\xi)$ in \eqref{eq: ul grad chain rule} are easy to compute.

When the exact lower-level solution $\eta^*(\xi)$ is unavailable, the upper-level gradient $\nabla u(\xi)$ is inaccessible. 
However, we can approximate $\eta^*(\xi)$ and therefore approximate $\nabla_\xi u(\xi)$.
Specifically, for $\xi^{k_u}$ at the $k_u$-th iteration, we run $K_l$-step gradient descent of lower-level with stepsize $\stp_l$ to approximate $\eta^*(\xi^{k_u})$, i.e.
\begin{equation}\left\{
\begin{aligned}
    &\eta^{k_u,1} = \eta^{k_u};\\
    &\eta^{k_u,k_l+1}=\eta^{k_u,k_l}-\stp_l\nabla_\eta \calL\left(\eta^{k_u,k_l};\xi^{k_u}\right),k_l=1,\cdots,K_l;\\
    &\eta^{k_u+1} = \eta^{k_u,K_l+1}.    
\end{aligned}\right.
\label{eq: ll update uncstr}
\end{equation}
It is easy to see that $\eta^{k_u,k_l+1}=\eta^{k_u,k_l+1}(\xi^{k_u})(k_l=1,\cdots,K_l)$ and $\eta^{k_u+1}=\eta^{k_u+1}(\xi^{k_u})$ are functions of $\xi^{k_u}$. We drop the dependence for notation conciseness and estimate the upper-level gradient $\nabla u(\xi^{k_u})$ by
\begin{equation}
\begin{aligned}
     \widehat{\nabla}u(\xi^{k_u})
    :=\left(\nabla_{\xi^{k_u}} \eta^{k_u+1}\right)^\top \nabla_{\eta}\calU\left(\eta^{k_u+1};\xi^{k_u}\right)+
    \nabla_{\xi} \calU\left(\eta^{k_u+1};\xi^{k_u}\right),  
\end{aligned}
\label{eq: ul grad unroll uncstr}
\end{equation}
where the $\eta^{k_u}$ is a lower-level estimator of the lower-level optimizer $\eta^*(\xi^{k_u})$, 
and $(\nabla_{\xi^{k_u}} \eta^{k_u+1})_{ij} = \partial_{\xi^{k_u}_j} \eta^{k_u+1}_i$ estimates $\nabla_{\xi}\eta^*(\xi^{k_u})$ by unrolling the lower-level iterates through the chain rule.
With the estimator in \eqref{eq: ul grad unroll uncstr}, we then update the upper-level variable by gradient descent with stepsize $\stp_u$, i.e.
\begin{equation}
    \xi^{k_u+1} = \xi^{k_u} - \stp_u \widehat{\nabla}u(\xi^{k_u}).
\label{eq: ul update uncstr}
\end{equation}

We summarize the algorithm in Algorithm \ref{alg: AGM general uncstr}
\begin{algorithm}[h]
\caption{General AGM for unconstrained bilevel optimization problem \eqref{eq: general blo uncstr}}
\begin{algorithmic}
\STATE{Initialization:  $\xi^1,\eta^1$, stepsizes $\{\stp_u,\stp_l\}$}\\
\FOR{$k_u=1,2,\cdots,K_u$}
    \STATE Initialize lower-level update by $\eta^{k_u,1}=\eta^{k_u}$.
    \FOR{$k_l=1,2\cdots,K_l$}
        \STATE lower-level gradient descent
        \begin{equation}
            \eta^{k_u,k_l+1}= \eta^{k_u,k_l}-\stp_l\nabla_\eta \mathcal{L}(\eta^{k_u,k_l};\xi^{k_u}).
        \end{equation}
        
    \ENDFOR
    \STATE Let the lower-level estimator be $\eta^{k_u+1} = \eta^{k_u,K_l+1}$ and compute $\widehat{\nabla}u(\xi^{k_u})$ by \eqref{eq: ul grad unroll uncstr}.
    \STATE Conduct upper-level gradient descent
    \begin{equation}
        \xi^{k_u+1}=\xi^{k_u}-\stp_u \widehat{\nabla}u(\xi^{k_u}).
    \end{equation}
\ENDFOR
\end{algorithmic}
\label{alg: AGM general uncstr}
\end{algorithm}

\begin{remark}[Error of unrolled differentiation]
\label{rem: unroll differentiation}
\eqref{eq: ul grad chain rule} gives the exact value of the upper-level gradient.  To obtain the unknown $\nabla_{\xi}\eta^*(\xi)$ in \eqref{eq: ul grad chain rule}, we refer to the first-order optimality condition from the lower-level problem $\nabla_\eta\calL(\eta^*(\xi);\xi)=\zero$.
We view $\nabla_\eta\calL(\eta^*(\xi);\xi)$ as a vector-valued function of $\xi$, and its Jacobian matrix gives
\begin{equation}
    \nabla_{\xi}\eta^*(\xi)^\top  \nabla_{\eta\eta}\calL(\eta^*(\xi);\xi)
    +\nabla_{\xi\eta}\calL(\eta^*(\xi);\xi)=\zero,
\label{eq: ll opt}
\end{equation}
where $(\nabla_{\xi\eta}\calL)_{ij}(\eta,\xi)=\partial_{\xi_i}\partial_{\eta_j}\calL(\eta,\xi)$ and $(\nabla_{\eta\eta}\calL)_{ij}(\eta,\xi)=\partial_{\eta_i}\partial_{\eta_j}\calL(\eta,\xi)$ are blocks of the Hessian matrix of $\calL$.
Therefore
\begin{equation}
    \nabla_{\xi}\eta^*(\xi)^\top=
    -\nabla_{\xi\eta}\calL(\eta^*(\xi);\xi)
    \left(\nabla_{\eta\eta}\calL(\eta^*(\xi);\xi)\right)^{-1}.
\label{eq: optimizer grad}
\end{equation}
Plugging \eqref{eq: optimizer grad} into \eqref{eq: ul grad chain rule} gives the upper-level gradient
\begin{equation}
    \nabla u(\xi) = \widehat{\nabla}_{\xi}\calU(\eta^*(\xi);\xi),
\label{eq: ul grad}
\end{equation}
where
\begin{equation}
    \widehat{\nabla}_{\xi}\calU(\eta;\xi)=\nabla_\xi\calU(\eta;\xi)
                    -\nabla_{\xi\eta}\calL(\eta;\xi)
                    \left(\nabla_{\eta\eta}\calL(\eta;\xi)\right)^{-1}
                    \nabla_\eta\calU(\eta;\xi).
\label{eq: ul grad est uncstr}
\end{equation}
The gradient estimator \eqref{eq: ul grad unroll uncstr} approximate the true gradient by approximating $\eta^*$ by $\eta^{k_u+1}$ and approximating $\left(\nabla_{\eta\eta}\calL(\eta;\xi)\right)^{-1}$ by unrolling differentiation.
A key to the convergence of the AGM algorithm is to control the error of unrolling differentiation. For unconstrained problems, \cite{grazzi2020iteration,ji2021bilevel} proved that under sufficient smoothness assumptions, the errors of the approximations decrease as $k_l$ increases.
In Lemma \ref{lem: lm-sub} of this paper, we study and prove the error can also be bounded for linear equality constrained lower-level problems. 

\end{remark}

\subsection{AGM for Inverse Mean-Field Games}

Building upon Algorithm \ref{alg: AGM general uncstr} for unconstrained bilevel optimization problems \eqref{eq: general blo uncstr}, we propose Algorithm \ref{alg: AGM general} to solve the constrained bilevel optimization problem \eqref{eq: general blo} and its special cases in inverse mean-field game problems \ref{eg: blo obs} and \ref{eg: blo metric}. 

\begin{algorithm}[h]
\caption{General AGM for \eqref{eq: general blo}}
\begin{algorithmic}
\STATE{Initialization:  $\xi^1,\eta^1$, stepsizes $\{\stp_u,\stp_l\}$}\\
\FOR{$k_u=1,2,\cdots,K_u$}
    \STATE Initialize lower-level update by $\eta^{k_u,1}=\eta^{k_u}$.
    \FOR{$k_l=1,2\cdots,K_l$}
        \STATE lower-level gradient descent
        \begin{equation}
            \eta^{k_u,k_l+1}= \proj_H\left(\eta^{k_u,k_l}-\stp_l\nabla_\eta \mathcal{L}(\eta^{k_u,k_l};\xi^{k_u})\right).
        \end{equation}
        
    \ENDFOR
    \STATE Let the lower-level estimator be $\eta^{k_u+1} = \eta^{k_u,K_l+1}$ and compute $\widehat{\nabla}u(\xi^{k_u})$ by \eqref{eq: ul grad unroll uncstr}.
    \STATE Conduct upper-level projected gradient descent
    \begin{equation}
        \xi^{k_u+1}=\proj_{\Xi}\left(\xi^{k_u}-\stp_u \widehat{\nabla}u(\xi^{k_u}) \right)
    \end{equation}
\ENDFOR
\end{algorithmic}
\label{alg: AGM general}
\end{algorithm}

Algorithm \ref{alg: AGM general} applies the projected gradient descent to estimate the lower-level optimizer and to update the upper-level optimizer at each iteration. 
Precisely, by denoting the matrix form of the constraint \eqref{eq: disct cstr set cteqn} as $A\eta=c$, the projection to $H=\{\eta\mid A\eta=c\}$ is
$$\proj_H(\eta)=(I-A^\dagger A)\eta + \eta_0,$$
where 
$A^\dagger$ is the Moore-Penrose inverse and $\eta_0$ is a fixed solution to $A\eta=c$.
The projection operator is invariant to the lower-level objective and the number of iterations.
As discussed in \cite{yu2021fast}, the main cost of the lower-level projected gradient descent is to compute the inverse of the discretized Laplacian operator $(AA^\top)^{-1}$, which can be solved efficiently using the fast cosine transform. We refer to section 3.2 in \cite{yu2021fast} for all detailed discussions.
Since each step in projected gradient descent is explicit, it is possible to unroll the differentiation to estimate the upper-level gradient and thus conduct AGM for the constrained bilevel optimization problem.
It is worth emphasizing that the proximal gradient solver for the lower-level problem \cite{yu2021fast} makes it easy and efficient to unroll the differentiation and estimate the upper-level gradient. 
This is not the case for other popular lower-level solvers, for example, augmented Lagrangian \cite{benamou2014augmented,benamou2015augmented}and primal-dual \cite{papadakis2014optimal,papadakis2015optimal} because the implicit updates in these methods make unrolling the differentiation prohibitively complicated and expensive. 
Meanwhile, although it is widely acknowledged in unconstrained bilevel optimization \cite{grazzi2020iteration,ji2021bilevel} that the error arising from unrolling differentiation is controllable, rigorously adapting this approach to incorporate lower-level linear constraints is, to the best of our knowledge, unexplored. Lemma \ref{lem: lm-sub} in this paper investigates the error of this approximation, indicating that the gradient estimation error can be effectively bounded by the accuracy of the lower-level solution. 

The complexity of resolving the upper-level constraint is similar to a single-level optimization problem.
In our cases, for the inverse crowd motion problem \ref{eg: disct blo obs}, the upper-level constraint set $\Xi=\calC_{\calG,\obs}$ as defined in \eqref{eq: disct cstr set obs} is the set of matrices of size $n_x\times n_y$ with entry sum zero. And the projection is simply  $\proj_{\calC_{\calG,\obs}}(\obs)=\tilde{\obs}$, where $(\tilde{\obs})_{\idxy}=(\obs)_{\idxy}-\frac{1}{\nx\ny}\sum_{\idxy}^{\nx,\ny}(\obs)_{\idxy}$. 
And for the inverse metric problem \ref{eg: disct blo metric}, $\Xi=\calC_{\calG,\metric}$, where $\calC_{\calG,\metric}$ is defined in \eqref{eq: disct cstr set metric}. We compute the projection $\tilde{\metric}:=\proj_{\calC_{\calG,\metric}}(\metric)$ pointwisely. To be specific, for $(\metric)_{\idxy}$,  we first compute its eigenvalue decomposition $(\metric)_{\idxy}=Q\Lambda Q^{-1}$ where $\Lambda=\diag(\lambda_1,\lambda_2)$ and let $(\tilde{\metric})_\idxy:=Q\tilde{\Lambda}Q^{-1}$ where $\tilde{\Lambda}=\diag(\max(\lambda_1,\epsilon),\max(\lambda_2,\epsilon))$ with a pre-selected small positive value $\epsilon$. 

Different from our bilevel formulation and AGM algorithm, \cite{chow2022numerical,ding2022mean} treat the forward MFG PDE system as the constraint of their optimization problem and apply primal-dual algorithm \cite{chambolle2011first} to solve it. However, the nonlinear and nonconvex constraint makes it challenging to prove the algorithm convergence. 
On the contrary, our bilevel formulation takes advantage of the convex-linear structure of the forward MFG and we establish the following convergence theorem of our Algorithm \ref{alg: AGM general}. 

If the upper-level and lower-level objective functions satisfy the following regularity assumptions,
\begin{assumption}\label{as1}
Assume that $\mathcal{U},\nabla \mathcal{U},\nabla \mathcal{L},\nabla^2 \mathcal{L}$ is Lipschitz continuous with $\ell_{u,0},\ell_{u,1},\ell_{l,1},$
$\ell_{l,2}$, respectively.  
\end{assumption}
\begin{assumption}\label{as2}
For any fixed $\xi$, assume that $\mathcal{L}(\eta;\xi)$ is $\mu_l$-strongly convex with respect to $\eta$. 
\end{assumption}
\begin{assumption}\label{ash}
$\Xi$ is a linear constraint set $\Xi=\{\xi\mid B\xi=e\}$, and $H$ and $\Xi$ are nonempty. 
\end{assumption}
then we have the following theorem.
\begin{theorem}
Under Assumption \ref{as1}--\ref{ash}, 
let $\stp_l\leq\frac{1}{2\ell_{l,1}}, K_l={\cal O}(\log K_u)$ and $\stp_u={\cal O}(1)$, then the iterates of Algorithm \ref{alg: AGM general} satisfy
\begin{align}
    \frac{1}{K_u}\sum_{k_u=1}^{K_u}\|\xi^{k_u}-\proj_{\Xi}(\xi^{k_u}-\nabla u(\xi^{k_u}))\|^2={\cal O}\left(\frac{1}{K_u}\right)
\end{align}
where ${\mathcal O}$ omits the $log$ dependency. 
\label{thm-conv}
\end{theorem}
Let us define $\epsilon$ stationary point as $\|\xi-\proj_{\Xi}(\xi-\nabla u(\xi))\|^2\leq\epsilon$,
then Theorem \ref{thm-conv} states that Algorithm \ref{alg: AGM general} achieves $\epsilon$ stationary point by ${\mathcal O}(\epsilon^{-1})$ iterations. This matches the iteration complexity of the single-level projected gradient descent method. We postpone the proof in section \ref{sec: proof}.

Lemma \ref{lem: convexity} states that when $\rho,\vm$ are bounded, and when the entropy in the objective function is non-zero ($\lambda>0$), then our inverse problems \ref{eg: disct blo obs} and \ref{eg: disct blo metric} satisfy assumptions \ref{as1} and \ref{as2}.
Moreover, since the upper-level constraint set of the inverse crowd motion problem \ref{eg: disct blo obs} is linear, assumption \ref{ash} is satisfied and Theorem \ref{thm-conv} guarantees the algorithm convergence when solving problem \ref{eg: disct blo obs}. 
For the inverse metric problem \ref{eg: disct blo metric} where the upper-level constraint set is a convex cone, the convergence of the algorithm can be established similarly. However, the convergence rate is possibly different. We leave the study of the convergence rate for general upper-level constraints set to future research.

At the end of this section, we discuss how to unroll differentiation in practice. 
\begin{remark}[Unroll differentiation in practice]
\label{rem: comp unroll diff}
Recall that in our problem, the lower-level variable $\eta=(\rho_{\calG^{\rho}},\vm_{\calG^{\vm}})$ and the upper-level variable $\xi=(\metric_{\calG},\obs_{\calG})$ are of size $\mathcal{O}(d^2\nt\nx\ny)$.
To obtain the upper-level gradient estimator \eqref{eq: ul grad unroll uncstr}, the computation of $\nabla_{\xi} \calU\left(\eta^{k_u+1};\xi^{k_u}\right)$ is straightforward. 
But it is not practicable to directly formulate $\nabla_{\xi^{k_u}} \eta^{k_u+1}$ since the size of the Jacobian matrix is $\mathcal{O}(d\nt\nx\ny) \times \mathcal{O}(d^2\nt\nx\ny)$ and the sparsity structure of the Jacobian matrix is not straightforward. 
Denote the gradient descent mapping $\gd(\eta;\xi):=\eta-\stp_l\nabla_\eta \calL\left(\eta;\xi\right)$. Then the Jacobian of $\gd$, $\nabla\gd=(\nabla_\eta\gd, \nabla_\xi\gd)=(I-\stp\nabla_{\eta\eta}\calL, -\stp\nabla_{\eta\xi}\calL)$ is sparse because the number of non-zero entries of $\nabla_{\eta\eta}\calL$ and $\nabla_{\eta\xi}\calL$ is $\mathcal{O}(d\nt\nx\ny)$.
In practice, we avoid formulating the matrix $\nabla_{\xi^{k_u}} \eta^{k_u+1}$ by chain rule and the sparsity structure of $\nabla\gd$.
Specifically, let $P:=I-A^\dagger A$ be the projection matrix, $\nabla_{\eta^{k_u,k_l}}\calU\left(\eta^{k_u+1};\xi^{k_u}\right)$ be the gradient of $\calU\left(\eta^{k_u+1};\xi^{k_u}\right)$ with respect to $\eta^{k_u,k_l}$,
and $\nabla_{\xi^{k_u}}\eta^{k_u,k_l}$ be the Jacobian of $\eta^{k_u,k_l}$ with respect to $\xi^{k_u}$,
then we have the following relation by back-propagation
\begin{equation}\left\{
\begin{aligned}
    &\nabla_{\eta^{k_u,K_l+1}}\calU\left(\eta^{k_u+1};\xi^{k_u}\right) 
    = \nabla_\eta\calU\left(\eta^{k_u+1};\xi^{k_u}\right),\\ 
    &\nabla_{\eta^{k_u,k_l}}\calU\left(\eta^{k_u+1};\xi^{k_u}\right)
    = \left( \nabla_\eta\gd(\eta^{k_u,k_l};\xi^{k_u}) \right)^\top P \nabla_{\eta^{k_u,k_l+1}}\calU\left(\eta^{k_u};\xi^{k_u}\right),
     \quad k_l=1,\cdots,K_l.
\end{aligned}\right.
\label{eq: ul grad unroll llv}
\end{equation}
Consequently, the upper-level gradient estimator is
\begin{equation}
\begin{aligned}
     \widehat{\nabla}u(\xi^{k_u})
    =&\left(\nabla_{\xi^{k_u}} \eta^{k_u,K_l+1}\right)^\top \nabla_{\eta^{k_u,K_l+1}}\calU\left(\eta^{k_u+1};\xi^{k_u}\right)
    + \nabla_{\xi} \calU\left(\eta^{k_u+1};\xi^{k_u}\right)\\
    =&\left( \nabla_{\eta}\gd(\eta^{k_u,K_l};\xi^{k_u}) \nabla_{\xi^{k_u}}\eta^{k_u,K_l} \right)^\top P\nabla_{\eta^{k_u,K_l+1}}\calU\left(\eta^{k_u+1};\xi^{k_u}\right)\\
    & +\left( \nabla_{\xi}\gd(\eta^{k_u,K_l};\xi^{k_u}) \right)^\top P\nabla_{\eta^{k_u,K_l+1}}\calU\left(\eta^{k_u+1};\xi^{k_u}\right) + \nabla_{\xi} \calU\left(\eta^{k_u+1};\xi^{k_u}\right)\\
    \stackrel{\text{by }\eqref{eq: ul grad unroll llv}}{=}
    &\left( \nabla_{\xi^{k_u}}\eta^{k_u,K_l} \right)^\top \nabla_{\eta^{k_u,K_l}}\calU\left(\eta^{k_u+1};\xi^{k_u}\right) \\
    & +\left( \nabla_{\xi}\gd(\eta^{k_u,K_l};\xi^{k_u}) \right)^\top P\nabla_{\eta^{k_u,K_l+1}}\calU\left(\eta^{k_u+1};\xi^{k_u}\right) + \nabla_{\xi} \calU\left(\eta^{k_u+1};\xi^{k_u}\right)\\
    =&\left( \nabla_{\eta}\gd(\eta^{k_u,K_l-1};\xi^{k_u}) \nabla_{\xi^{k_u}}\eta^{k_u,K_l-1} \right)^\top P\nabla_{\eta^{k_u,K_l}}\calU\left(\eta^{k_u+1};\xi^{k_u}\right) \\
    & +\left( \nabla_{\xi}\gd(\eta^{k_u,K_l-1};\xi^{k_u}) \right)^\top P\nabla_{\eta^{k_u,K_l}}\calU\left(\eta^{k_u+1};\xi^{k_u}\right)\\
    & +\left( \nabla_{\xi}\gd(\eta^{k_u,K_l};\xi^{k_u}) \right)^\top P\nabla_{\eta^{k_u,K_l+1}}\calU\left(\eta^{k_u+1};\xi^{k_u}\right) + \nabla_{\xi} \calU\left(\eta^{k_u+1};\xi^{k_u}\right)\\
    \stackrel{\text{by }\eqref{eq: ul grad unroll llv}}{=}
    & \left( \nabla_{\xi^{k_u}}\eta^{k_u,K_l-1} \right)^\top \nabla_{\eta^{k_u,K_l-1}}\left(\eta^{k_u+1};\xi^{k_u}\right) \\
    & + \sum_{i=K_l-1}^{K_l}
    \left( \nabla_{\xi}\gd(\eta^{k_u,i};\xi^{k_u}) \right)^\top P\nabla_{\eta^{k_u,i+1}}\calU\left(\eta^{k_u+1};\xi^{k_u}\right)
    + \nabla_{\xi} \calU\left(\eta^{k_u+1};\xi^{k_u}\right)\\
    =& \cdots\\
    \stackrel{(a)}{=}
    & \sum_{i=1}^{K_l}    
    \left( \nabla_{\xi}\gd(\eta^{k_u,i};\xi^{k_u}) \right)^\top P\nabla_{\eta^{k_u,i+1}}\calU\left(\eta^{k_u+1};\xi^{k_u}\right)
    + \nabla_{\xi} \calU\left(\eta^{k_u+1};\xi^{k_u}\right)
\end{aligned}
\label{eq: ul grad unroll ulv}
\end{equation}
where $(a)$ is because that $\eta^{k_u,1}$ is independent of $\xi^{k_u}$. In this way, each term in the estimator can be computed by sparse matrix and vector multiplication.
\end{remark}

\section{Proofs of Main Theorems}
\label{sec: proof}
In this section, we provide the proofs of main theorems. 
Theorem \ref{thm: differentiability} shows that the observations of the Nash Equilibrium continuously depend on the unknown parameters.
Theorem \ref{thm: unique identifiablity} states that with only one good observation of the Nash Equilibrium, we can uniquely recover the obstacle $\obs$ up to a constant by solving the bilevel problem \eqref{eq: disct inv mfg obs}.
This illustrates the effectiveness of our model. 
Lemma \ref{lem: convexity} and Theorem \ref{thm-conv} together guarantee that Algorithm \ref{alg: AGM general} converges to a stationary point to the bilevel problem \eqref{eq: disct inv mfg obs} if the forward problem has enough regularity.
This illustrates the effectiveness of our algorithm.

\subsection{Proof of Theorem \ref{thm: differentiability} and \ref{thm: unique identifiablity}}

Recall that $\calK(\rho,\vm,\phi;\metric,\obs)=\zero$ gives the optimality condition.
Denote the Jacobian matrix of $\calK$ as $\nabla\calK=\left( (\nabla_{\rho,\vm,\phi}\calK)_{d_l\times d_l},(\nabla_{\metric,\obs}\calK)_{d_l\times d_u} \right)$.
The proof of the regularity Theorem \ref{thm: differentiability} is based on the implicit function theorem and the key is to show that the matrix $\nabla_{\rho,\vm,\phi}\calK$ is invertible at a good observation $(\rhotilde,\vmtilde,\phitilde;\mettilde,\obstilde)$.
\begin{lemma}
\label{lem: invertible}
    If $\weight_I>0$,$\weight_T>0$ and $\min_{\vi\in\Grho}\{\rhotilde_\vi\}>0$, then $\nabla_{\rho,\vm,\phi}\calK(\rhotilde,\vmtilde,\phitilde;\mettilde,\obstilde)$ is invertible.
\end{lemma}

\begin{proof}
To prove the lemma is equivalent to showing that 
\begin{equation}
    \nabla_{\rho,\vm,\phi}\calK(\rhotilde,\vmtilde,\phitilde;\mettilde,\obstilde)(\ptbrho,\ptbvm,\ptbphi)=\zero,
\label{eq: hessian*v=0}
\end{equation}
if and only if $(\ptbrho,\ptbvm,\ptbphi)=\zero$. Here $\ptbvm:=\{\ptbmx,\ptbmy\}$.
By definition, 
\begin{equation}
    \nabla_{\rho,\vm,\phi}\calK(\rhotilde,\vmtilde,\phitilde;\mettilde,\obstilde)(\ptbrho,\ptbvm,\ptbphi)
    =\lim_{\epsilon\to0}\frac{1}{\epsilon}\left( 
    \calK(\rhotilde+\epsilon\ptbrho,\vmtilde+\epsilon\ptbvm,\phitilde+\epsilon\ptbphi;\mettilde,\obstilde) - \calK(\rhotilde,\vmtilde,\phitilde;\mettilde,\obstilde)
    \right)
\end{equation}
Therefore \eqref{eq: hessian*v=0} is equivalent to
\begin{equation}
\allowdisplaybreaks\left\{
\begin{aligned}
&\vi\in\Grho,\idt=1,\cdots,\nt-1,\\
&   -(\Dt^*(\ptbphi))_\vi 
    + \Bigg(\It^*\Bigg( 
        -\frac{\metric_{xx}\overline{\mxtilde}+\metric_{xy}\overline{\mytilde}}{\overline{\rhotilde}^2}\overline{\ptbmx} 
        -\frac{\metric_{xy}\overline{\mxtilde}+\metric_{yy}\overline{\mytilde}}{\overline{\rhotilde}^2}\overline{\ptbmy} \\
& \qquad\qquad \qquad\qquad \qquad\qquad \qquad\qquad \qquad
        + \frac{(\overline{\vmtilde})^{\top}\metric\overline{\vmtilde}}{\overline{\rhotilde}^3}\overline{\ptbrho} 
        + \frac{\weight_I}{\overline{\rhotilde}}\overline{\ptbrho} \Bigg)\Bigg)_\vi=0,\\ 
&\vi\in\Grho,\it=\nt,\\
&   -(\Dt^*(\ptbphi))_\vi 
    + \Bigg(\It^*\Bigg( 
        -\frac{\metric_{xx}\overline{\mxtilde}+\metric_{xy}\overline{\mytilde}}{\overline{\rhotilde}^2}\overline{\ptbmx} 
        -\frac{\metric_{xy}\overline{\mxtilde}+\metric_{yy}\overline{\mytilde}}{\overline{\rhotilde}^2}\overline{\ptbmy} \\
& \qquad\qquad \qquad\qquad \qquad\qquad \qquad\qquad \qquad        
        +\frac{(\overline{\vmtilde})^{\top}\metric\overline{\vmtilde}}{\overline{\rhotilde}^3}\overline{\ptbrho} 
        +\frac{\weight_I}{\overline{\rhotilde}}\overline{\ptbrho} \Bigg)\Bigg)_\vi  
    + \frac{\weight_T}{\dt(\rhotilde)_\vi}(\ptbrho)_{\vi}=0,\\
&\vi\in\Gmx,\quad
    -(\Dx^*(\ptbphi))_\vi 
    + \left(\Ix^*\left( \frac{\metric_{xx}}{\overline{\rhotilde}}\overline{\ptbmx}
        +\frac{\metric_{xy}}{\overline{\rhotilde}}\overline{\ptbmy}
        - \frac{\metric_{xx}\overline{\mxtilde}+\metric_{xy}\overline{\mytilde}}{\overline{\rhotilde}^2}\overline{\ptbrho} \right)\right)_\vi=0,\\
&\vi\in\Gmy,\quad
    -(\Dy^*(\ptbphi))_\vi 
    + \left(\Iy^*\left( \frac{\metric_{xy}}{\overline{\rhotilde}}\overline{\ptbmx}
        +\frac{\metric_{yy}}{\overline{\rhotilde}}\overline{\ptbmy} 
        - \frac{\metric_{xy}\overline{\mxtilde}+\metric_{yy}\overline{\mytilde}}{\overline{\rhotilde}^2}\overline{\ptbrho} \right)\right)_\vi=0,\\    
&\vi\in\Gphi,\quad
    (\Dt(\ptbrho;\zero)+\Dx(\ptbmx)+\Dy(\ptbmy))_\vi=0.
\end{aligned}\right.
\label{eq: hessian*v=0, long}
\end{equation}
Note that $\rhotilde,\vmtilde,\phitilde$ are viewed as constants with respect to $(\ptbrho,\ptbvm,\ptbphi)$ in the system. It clear that the system \ref{eq: hessian*v=0, long} is linear in $(\ptbrho,\ptbvm,\ptbphi)$ and therefore \eqref{eq: hessian*v=0} holds if $(\ptbrho,\ptbvm,\ptbphi)=\zero$.
If both $(\ptbrho,\ptbvm,\ptbphi)$ and $(\ptbrho',\ptbvm',\ptbphi')$ satisfy \eqref{eq: hessian*v=0}, then by plugging them into \eqref{eq: hessian*v=0, long} and subtracting, we have
\begin{equation}\left\{
\begin{aligned}
&\vi\in\Grho,\idt=1,\cdots,\nt-1,\\
&   -(\Dt^*(\ptbphi-\ptbphi'))_\vi 
    + \Bigg(\It^*\Bigg( 
        -\frac{\metric_{xx}\overline{\mxtilde}+\metric_{xy}\overline{\mytilde}}{\overline{\rhotilde}^2}(\overline{\ptbmx}-\overline{\ptbmx'})\\
& \qquad\qquad\qquad\qquad\qquad\qquad        
        -\frac{\metric_{xy}\overline{\mxtilde}+\metric_{yy}\overline{\mytilde}}{\overline{\rhotilde}^2}(\overline{\ptbmy}-\overline{\ptbmy'})  \\
& \qquad\qquad\qquad\qquad\qquad\qquad
        + \frac{(\overline{\vmtilde})^\top\metric\overline{\vmtilde}}{\overline{\rhotilde}^3}(\overline{\ptbrho}-\overline{\ptbrho'}) 
        + \frac{\weight_I}{\overline{\rhotilde}}(\overline{\ptbrho}-\overline{\ptbrho'}) \Bigg)\Bigg)_\vi=0,\\ 
&\vi\in\Grho,\it=\nt,\\
&   -(\Dt^*(\ptbphi-\ptbphi'))_\vi 
        + \Bigg(\It^*\Bigg( -\frac{\metric_{xx}\overline{\mxtilde}+\metric_{xy}\overline{\mytilde}}{\overline{\rhotilde}^2}(\overline{\ptbmx}-\overline{\ptbmx'})\\
& \qquad\qquad\qquad\qquad\qquad\qquad        
        -\frac{\metric_{xy}\overline{\mxtilde}+\metric_{yy}\overline{\mytilde}}{\overline{\rhotilde}^2}(\overline{\ptbmy}-\overline{\ptbmy'})  \\
& \qquad\qquad\qquad\qquad\qquad\qquad
        +\frac{(\overline{\vmtilde})^\top\metric\overline{\vmtilde}}{\overline{\rhotilde}^3}(\overline{\ptbrho}-\overline{\ptbrho'}) 
        +\frac{\weight_I}{\overline{\rhotilde}}(\overline{\ptbrho}-\overline{\ptbrho'}) \Bigg)\Bigg)_\vi  
    + \frac{\weight_T}{\dt(\rhotilde)_\vi}(\ptbrho-\ptbrho')_{\vi}=0,\\
\end{aligned}\right.
\label{eq: hessian*v=0, rho}
\end{equation}
\begin{equation}\left\{
\begin{aligned}
&\vi\in\Gmx,\quad
    -(\Dx^*(\ptbphi-\ptbphi'))_\vi 
    + \Bigg(\Ix^*\Bigg( \frac{\metric_{xx}}{\overline{\rhotilde}}(\overline{\ptbmx}-\overline{\ptbmx'}) 
        + \frac{\metric_{xy}}{\overline{\rhotilde}}(\overline{\ptbmy}-\overline{\ptbmy'})\\
&\qquad\qquad\qquad\qquad\qquad\qquad\qquad\qquad
        - \frac{\metric_{xx}\overline{\mxtilde}+\metric_{xy}\overline{\mytilde}}{\overline{\rhotilde}^2}(\overline{\ptbrho}-\overline{\ptbrho'}) \Bigg)\Bigg)_\vi=0,\\
&\vi\in\Gmy,\quad
    -(\Dy^*(\ptbphi-\ptbphi'))_\vi 
    + \Bigg(\Iy^*\Bigg( \frac{\metric_{xy}}{\overline{\rhotilde}}(\overline{\ptbmx}-\overline{\ptbmx'}) 
        + \frac{\metric_{yy}}{\overline{\rhotilde}}(\overline{\ptbmy}-\overline{\ptbmy'})\\
&\qquad\qquad\qquad\qquad\qquad\qquad\qquad\qquad
        - \frac{\metric_{xy}\overline{\mxtilde}+\metric_{yy}\overline{\mytilde}}{\overline{\rhotilde}^2}(\overline{\ptbrho}-\overline{\ptbrho'}) \Bigg)\Bigg)_\vi=0,\\    
\end{aligned}\right.
\label{eq: hessian*v=0, m}
\end{equation}
and
\begin{equation}
\vi\in\Gphi,\quad
    (\Dt(\ptbrho-\ptbrho';\zero)+\Dx(\ptbmx-\ptbmx')+\Dy(\ptbmy-\ptbmy'))_\vi=0.
\label{eq: hessian*v=0, phi}
\end{equation}
Pointwisely multiplying \eqref{eq: hessian*v=0, rho} with $(\ptbrho-\ptbrho')$ and summing over $\Grho$ gives us
\begin{equation}
\begin{aligned}
    & -\left\langle \ptbrho-\ptbrho',\Dt^*(\ptbphi-\ptbphi') \right\rangle_{\Grho}\\
    & -\left\langle \ptbrho-\ptbrho',\It^*\left( \frac{\metric_{xx}\overline{\mxtilde}+\metric_{xy}\overline{\mytilde}}{\overline{\rhotilde}^2}(\overline{\ptbmx}-\overline{\ptbmx'}) \right) \right\rangle_{\Grho}
      -\left\langle \ptbrho-\ptbrho',\It^*\left( \frac{\metric_{xy}\overline{\mxtilde}+\metric_{yy}\overline{\mytilde}}{\overline{\rhotilde}^2}(\overline{\ptbmy}-\overline{\ptbmy'})\right) \right\rangle_{\Grho}\\
    & +\left\langle \ptbrho-\ptbrho',\It^*\left( \frac{(\overline{\vmtilde})^\top\metric\overline{\vmtilde}}{\overline{\rhotilde}^3}(\overline{\ptbrho}-\overline{\ptbrho'})\right)  \right\rangle_{\Grho}
    +\left\langle \ptbrho-\ptbrho',\It^*\left(\frac{\weight_I}{\overline{\rhotilde}}(\overline{\ptbrho}-\overline{\ptbrho'})\right) \right\rangle_{\Grho}\\
    & +\dx\dy\sum_{\nx=1}^\nx\sum_{\ny=1}^\ny \frac{\weight_T}{(\rhotilde)_{\idx,\idy,\nt}} (\ptbrho-\ptbrho')_{\idx,\idy,\nt}^2=0.
\end{aligned}
\label{eq: hessian*v=0, in prod, rho}
\end{equation}
Similarly \eqref{eq: hessian*v=0, m} and \eqref{eq: hessian*v=0, phi} imply
\begin{equation}
\begin{aligned}
    & -\left\langle \ptbmx-\ptbmx',\Dx^*(\ptbphi-\ptbphi') \right\rangle_{\Gmx} \\
    & +\left\langle \ptbmx-\ptbmx', \Ix^*\left( \frac{\metric_{xx}}{\overline{\rhotilde}}(\overline{\ptbmx}-\overline{\ptbmx'})
    +  \frac{\metric_{xy}}{\overline{\rhotilde}}(\overline{\ptbmy}-\overline{\ptbmy'})
    \right) \right\rangle_{\Gmx} \\
    & -\left\langle \ptbmx-\ptbmx', \Ix^*\left( \frac{\metric_{xx}\overline{\mxtilde}+\metric_{xy}\overline{\mytilde}}{\overline{\rhotilde}^2}(\overline{\ptbrho}-\overline{\ptbrho'}) \right) \right\rangle_{\Gmx} =0,
\end{aligned}
\label{eq: hessian*v=0, in prod, mx}
\end{equation}
and 
\begin{equation}
\begin{aligned}
    & -\left\langle \ptbmy-\ptbmy',\Dy^*(\ptbphi-\ptbphi') \right\rangle_{\Gmy} \\
    & +\left\langle \ptbmy-\ptbmy', \Iy^*\left( \frac{\metric_{xy}}{\overline{\rhotilde}}(\overline{\ptbmx}-\overline{\ptbmx'})
    +\frac{\metric_{yy}}{\overline{\rhotilde}}(\overline{\ptbmy}-\overline{\ptbmy'})    
    \right) \right\rangle_{\Gmy} \\
    &-\left\langle \ptbmy-\ptbmy', \Iy^*\left( \frac{\metric_{xy}\overline{\mxtilde}+\metric_{yy}\overline{\mytilde}}{\overline{\rhotilde}^2}(\overline{\ptbrho}-\overline{\ptbrho'}) \right) \right\rangle_{\Gmy} =0,
\end{aligned}
\label{eq: hessian*v=0, in prod, my}
\end{equation}
and
\begin{equation}
    \left\langle \ptbphi-\ptbphi',\Dt(\ptbrho-\ptbrho';\zero) \right\rangle_{\Gphi} 
    +\left\langle \ptbphi-\ptbphi',\Dx(\ptbmx-\ptbmx') \right\rangle_{\Gphi} 
    +\left\langle \ptbphi-\ptbphi',\Dy(\ptbmy-\ptbmy') \right\rangle_{\Gphi}=0.
\label{eq: hessian*v=0, in prod, phi}
\end{equation}
Next, we add \eqref{eq: hessian*v=0, in prod, rho}-\eqref{eq: hessian*v=0, in prod, phi} and combine terms with the same components in groups. 
The first group is 
\begin{equation}
\begin{aligned}
    & -\left\langle \ptbrho-\ptbrho',\Dt^*(\ptbphi-\ptbphi') \right\rangle_{\Grho}
    -\left\langle \ptbmx-\ptbmx',\Dx^*(\ptbphi-\ptbphi') \right\rangle_{\Gmx} 
    -\left\langle \ptbmy-\ptbmy',\Dy^*(\ptbphi-\ptbphi') \right\rangle_{\Gmy} \\
    &\left\langle \ptbphi-\ptbphi',\Dt(\ptbrho-\ptbrho';\zero) \right\rangle_{\Gphi} 
    +\left\langle \ptbphi-\ptbphi',\Dx(\ptbmx-\ptbmx') \right\rangle_{\Gphi} 
    +\left\langle \ptbphi-\ptbphi',\Dy(\ptbmy-\ptbmy') \right\rangle_{\Gphi},  
\label{eq: group 1}
\end{aligned}
\end{equation}
and by the adjoint relation between $\Dt,\Dt^*$, this group sums to 0.
The second group consists of
\begin{equation}
\begin{aligned}
    &\left\langle \ptbrho-\ptbrho',\It^*\left(\frac{\weight_I}{\overline{\rhotilde}}(\overline{\ptbrho}-\overline{\ptbrho'})\right) \right\rangle_{\Grho}
    &+\dx\dy\sum_{\nx=1}^\nx\sum_{\ny=1}^\ny \frac{\weight_T}{(\rhotilde)_{\idx,\idy,\nt}} (\ptbrho-\ptbrho')_{\idx,\idy,\nt}^2.    
\label{eq: group 2}
\end{aligned}
\end{equation}
and the sum is equal to $\displaystyle\weight_I\left\| \frac{\overline{\ptbrho}-\overline{\ptbrho'}}{\overline{\rhotilde}^{1/2}} \right\|_{\Gphi}^2
    +\dx\dy\sum_{\nx=1}^\nx\sum_{\ny=1}^\ny \frac{\weight_T}{(\rhotilde)_{\idx,\idy,\nt}} (\ptbrho-\ptbrho')_{\idx,\idy,\nt}^2$.
The rest terms form the last group and sum to
\begin{equation}
    \sum_{\vi\in\calG^{\phi}}\frac{1}{\overline{\rhotilde}^3_{\vi}}
    \left\| \left( (\overline{\ptbrho}-\overline{\ptbrho'})\overline{\vmtilde} - \overline{\rhotilde}(\overline{\ptbvm}-\overline{\ptbvm'}) \right)_{\vi} \right\|^2_{\metric_{\vi}},
\label{eq: group 3}
\end{equation}
where $\|\vv_\vi\|^2_{\metric_{\vi}}=(\vv)_{\vi}^\top\metric_{\vi}\vv_{\vi}$.
Overall, adding \eqref{eq: hessian*v=0, in prod, rho}-\eqref{eq: hessian*v=0, in prod, phi} gives 
\begin{equation}
\begin{aligned}
    &\weight_I\left\| \frac{\overline{\ptbrho}-\overline{\ptbrho'}}{\overline{\rhotilde}^{1/2}} \right\|_{\Gphi}^2
    +\dx\dy\sum_{\nx=1}^\nx\sum_{\ny=1}^\ny \frac{\weight_T}{(\rhotilde)_{\idx,\idy,\nt}} (\ptbrho-\ptbrho')_{\idx,\idy,\nt}^2\\
    &+\sum_{\vi\in\calG^{\phi}}\frac{1}{\overline{\rhotilde}^3_{\vi}}
    \left\| \left( (\overline{\ptbrho}-\overline{\ptbrho'})\overline{\vmtilde} - \overline{\rhotilde}(\overline{\ptbvm}-\overline{\ptbvm'}) \right)_{\vi} \right\|^2_{\metric_{\vi}}=0.
\end{aligned}
\label{eq: positive integral=0}
\end{equation}
We conclude that each term in \eqref{eq: positive integral=0} is zero since they are non-negative and sum to zero. 
Combining $(\ptbrho)_{\idx,\idy,\nt}=(\ptbrho')_{\idx,\idy,\nt}$ and $\overline{\ptbrho}=\overline{\ptbrho'}$ gives $\ptbrho=\ptbrho'$. 
Consequently, $\overline{\ptbmx}=\overline{\ptbmx'}$ and $\overline{\ptbmy}=\overline{\ptbmy'}$. Because $\Ix,\Iy$ are full rank linear operators, $\ptbmx=\ptbmx'$ and $\ptbmy=\ptbmy'$.
Based on $\ptbrho=\ptbrho',\ptbvm=\ptbvm'$, \eqref{eq: hessian*v=0, rho} and \eqref{eq: hessian*v=0, m} lead to $\ptbphi=\ptbphi'$.
Therefore \eqref{eq: hessian*v=0} has unique solution $(\rho,\vm,\phi)=\zero$, i.e. $\nabla_{\rho,\vm,\phi}\calK(\rhotilde,\vmtilde,\phitilde;\obstilde)$ is invertible.
    
\end{proof}

With Lemma \ref{lem: invertible}, we apply implicit function theorem to $\calK$ at $(\rhotilde,\vmtilde,\phitilde;\mettilde,\obstilde)$ and then the regularity Theorem \ref{thm: differentiability} is true.

Next, we prove the unique identifiability Theorem \ref{thm: unique identifiablity} for inverse obstacle problem \ref{eg: disct blo obs}.

\begin{proof}[Proof of Theorem \ref{thm: unique identifiablity}]
    Since the upper-level objective is non-negative and equals 0 when $\obs=\obstilde$, any minimizer $\obs$ of the bilevel minimization problem satisfies
    \begin{equation}
        (\rhotilde,\vmtilde)=\argmin_{(\rho,\vm)\in\calC_\calG(\rhoinit)}\calL_\calG(\rho,\vm;\obs),
    \end{equation}
    and by Lemma \ref{lem: kkt -> zero function}, there exists $\phi$ such that $\calK(\rhotilde,\vmtilde,\phi;\obs)=\zero$.
    Assume that $\obs'$ is a minimizer, $\obs'\neq\obstilde$, and
    $$ \calK(\rhotilde,\vmtilde,\phitilde;\obstilde)= \calK(\rhotilde,\vmtilde,\phi';\obs') =\zero, $$    
    then 
    $$ \calK(\rhotilde,\vmtilde,\phitilde;\obstilde)- \calK(\rhotilde,\vmtilde,\phi';\obs') =\zero, $$
    which is equivalent to
    \begin{equation}\left\{
    \begin{aligned}
    &\vi\in\Grho,
    -\left( \Dt^*(\phi'-\phitilde) \right)_\vi 
        +((\obs')_\idxy-(\obstilde)_\idxy)=0,\\ 
    &\vi\in\Gmx,
        \left( \Dx^*(\phi'-\phitilde) \right)_\vi=0,\\
    &\vi\in\Gmy,
        \left( \Dy^*(\phi'-\phitilde) \right)_\vi =0.
    \end{aligned}\right.
    \label{eq: Y(b)-Y(b')}
    \end{equation} 
    The equation on $\Grho$ gives $(\phi'-\phitilde)_\id = (\nt-\idt+1)(\obs'-\obstilde)_\idxy$. Plugging in equations on $\Gmx,\Gmy$, we have $(\obs'-\obstilde)_\idxy=c$ where $c$ is a constant for different $\idx,\idy$.

\end{proof}

\subsection{Proof of Theorem \ref{thm-conv} and Lemma \ref{lem: convexity}}

In this section, we provide the nonasymptotic analysis for AGM on general constrained bilevel optimization \eqref{eq: general blo}. We follow conventional notations in bilevel optimization by using commas to separate lower-level and upper-level variables, i.e., $\calL(\eta,\xi)=\calL(\eta;\xi), \calU(\eta,\xi)=\calU(\eta;\xi)$ 

Recall that the lower level constraint is $H=\{\eta\mid A\eta=c\}$. 
Denote the singular value decomposition of $A$ as $A=U\Sigma V^\top $, where $$\Sigma=\left[\begin{array}{ll}
\Sigma_1 & 0 \\
0 & 0 
\end{array}\right]\in\mathbb{R}^{d_c\times d_\eta},$$ 
$U=[U_1~ U_2]$, $V=[V_1~ V_2]$, $U\in\mathbb{R}^{d_c\times d_c},V\in\mathbb{R}^{d_\eta\times d_\eta}$ are orthogonal matrix and $U_1\in\mathbb{R}^{d_c\times r}, V_1\in\mathbb{R}^{d_\eta\times r}$ are the submatrix corresponds to full rank diagonal submatrix $\Sigma_1\in\mathbb{R}^{r\times r}$. 
Then $V_2$ is the orthogonal basis of $\operatorname{Ker}(A):=\{\eta\mid A\eta=0\}$. 
Let $\eta_0\in H$ be a feasible lower-level solution, then the lower-level update is equivalent to
\begin{equation}
    \eta^{k_u,1}=\eta^{k_u}; ~~\eta^{k_u,k_l+1}=V_2V_2^\top \left(\eta^{k_u,k_l}-\tau_l\nabla_\eta \mathcal{L}(\eta^{k_u,k_l},\xi^{k_u})\right)+\eta_0;~~ \eta^{k_u+1}=\eta^{k_u,K_l+1}.
\label{eq: ll update for cvg}
\end{equation}
With $\eta^{k_u+1}$ approximating $\eta^*(\xi^{k_u})$, we approximate the lower-level gradient with
\begin{equation}
    \widehat{\nabla}u(\xi^{k_u}):=\nabla_{\xi} \calU(\eta^{k_u+1},\xi^{k_u})+(\nabla_{\xi^{k_u}}\eta^{k_u+1})^\top \nabla_{\eta}\calU(\eta^{k_u+1},\xi^{k_u})
\label{UL-unroll}    
\end{equation}
and $\nabla_{\xi^{k_u}} \eta^{k_u+1}$ is obtained by unrolling the lower-level iterates 
\begin{equation}\left\{
\begin{aligned}
    \nabla_{\xi^{k_u}} \eta^{k_u,1}&=\zero,\\  
     \nabla_{\xi^{k_u}} \eta^{k_u,k_l+1}
     &=V_2V_2^\top \nabla_{\xi^{k_u}} \eta^{k_u,k_l}-\tau_l V_2 V_2^\top\left(\nabla_{\eta\xi}\mathcal{L}(\eta^{k_u,k_l},\xi^{k_u})+\nabla_{\eta\eta}\mathcal{L}(\eta^{k_u,k_l},\xi^{k_u})\nabla_{\xi^{k_u}} \eta^{k_u,k_l}\right)\\
    &=V_2V_2^\top\left(I-\tau_l \nabla_{\eta\eta}\mathcal{L}(\eta^{k_u,k_l},\xi^{k_u})\right)\nabla_{\xi^{k_u}} \eta^{k_u,k_l}-\tau_l V_2 V_2^\top\nabla_{\eta\xi}\mathcal{L}(\eta^{k_u,k_l},\xi^{k_u}),k_l=1,\cdots,K_l. 
\end{aligned}\right.
\label{chain-rule-grad}
\end{equation}

To prove the convergence, we first present the regularity of the lower-level optimizer established in \cite{xiao2023alternating}. To be self-contained, we also provide its proof.

\begin{lemma}[The regularity of lower-level optimizer]
Under Assumption \ref{as1}--\ref{as2}, $\eta^*(\xi)$ is differentiable with respect to $\xi$ with the following gradient
\begin{align*}
    \nabla \eta^*(\xi)=-V_2(V_2^\top \nabla_{\eta\eta}\mathcal{L}(\eta^*(\xi),\xi)V_2)^{-1}V_2^\top\nabla_{\eta\xi}\mathcal{L}(\eta^*(\xi),\xi).
\end{align*}
where $V_2$ is the orthogonal basis of $\operatorname{Ker}(A)$.
Therefore, $\eta^*(\xi)$ is $L_\eta$-Lipschitz continuous and
and $L_{\eta\xi}$ smooth with 
\begin{align*}
    L_\eta:=\frac{\ell_{l,1}}{\mu_l}={\cal O}(\kappa),\qquad L_{\eta\xi}:=\frac{\ell_{l,2}(1+\frac{\ell_{l,1}}{\mu_l})^2}{\mu_l}={\cal O}(\kappa^3).
\end{align*}
\label{lem: y-smooth}
\end{lemma}

\begin{proof}\allowdisplaybreaks
First, we prove the differentiability and compute the Jacobian matrix. 
We choose a fixed $\eta_0$ satisfying $A\eta_0=c$. Using the aforementioned SVD of $A$, the constraint set $H= \{\eta_0+ V_2 z~|~ z\in \mathbb{R}^{d_\eta - r} \}$. 
Letting $\calL_z(z,\xi):=\calL(\eta_0+V_2z,\xi)$ and $z^*(\xi)=\arg\min_{z}\calL_z(z,\xi)$, we have $\eta^*(\xi)=\eta_0+V_2z^*(\xi)$. 
By optimality condition, $z^*(\xi)$ satisfies 
\begin{equation}
    \nabla_z \calL_z(z^*(\xi),\xi)=V_2^\top \nabla_\eta \mathcal{L}(\eta_0+V_2z^*(\xi),\xi)=0,
\label{eq: opt of z}
\end{equation}
Since
\begin{equation}
    \nabla_{zz}\calL_z(z^*(\xi),\xi)=V_2^\top \nabla_{\eta\eta} \mathcal{L}(\eta_0+V_2z^*(\xi),\xi) V_2
\end{equation}
and by strong convexity of $\calL$ with respect to $\eta$, $\nabla_{zz}\calL_z(z^*(\xi),\xi)$ is invertible. By implicit function theorem, $z^*(\xi)$ is differentiable with respect to $\xi$.
As a consequence, $\eta^*(\xi)$ is differentiable with respect to $\xi$.
Taking the gradient with respect to $\xi$ on both sides of \eqref{eq: opt of z} gives us
\begin{align*}
    0
    &=\nabla_{\xi\eta}\mathcal{L}(\eta_0+V_2z^*(\xi),\xi)V_2+\left(\nabla_\xi z^*(\xi)^\top V_2^\top\right)\nabla_{\eta\eta}\mathcal{L}(\eta_0+V_2z^*(\xi),\xi)V_2\\
    &=\nabla_{\xi\eta}\mathcal{L}(\eta_0+V_2z^*(\xi),\xi)V_2+\nabla_\xi z^*(\xi)^\top V_2^\top\nabla_{\eta\eta}\mathcal{L}(\eta_0+V_2z^*(\xi),\xi)V_2.
\end{align*}
Then, we have (cf. $\nabla_{\eta\eta}\mathcal{L}(\eta^*(\xi),\xi)=\nabla_{\eta\eta}\mathcal{L}(\eta_0+V_2z^*(\xi),\xi)$)
\begin{equation}
  \nabla z^*(\xi)=-\left(V_2^\top\nabla_{\eta\eta}\mathcal{L}(\eta^*(\xi),\xi)V_2\right)^{-1}V_2^\top\nabla_{\eta\xi}\mathcal{L}(\eta^*(\xi),\xi)  
\end{equation}
and as a result,
\begin{align*}
    \nabla \eta^*(\xi)&=V_2\nabla z^*(\xi)\\
    &=-V_2\left(V_2^\top\nabla_{\eta\eta}\mathcal{L}(\eta^*(\xi),\xi)V_2\right)^{-1}V_2^\top\nabla_{\eta\xi}\mathcal{L}(\eta^*(\xi),\xi).
\end{align*}

Next, utilizing the fact that $V_2$ is the orthogonal matrix, we know $\mu_l I\preceq V_2^\top\nabla_{\eta\eta}\calL(\eta,\xi)V_2$. Therefore, we have for any $\xi,\eta$, 
\begin{align}\label{bound_v}
    V_2\left(V_2^\top\nabla_{\eta\eta}\mathcal{L}(\eta,\xi)V_2\right)^{-1}V_2^\top\preceq\frac{1}{\mu_l} I.
\end{align}
As a result, $\nabla \eta^*(\xi)$ is bounded by
\begin{align*}
    \|\nabla \eta^*(\xi)\|&\leq \|V_2\left(V_2^\top\nabla_{\eta\eta}\mathcal{L}(\eta^*(\xi),\xi)V_2\right)^{-1}V_2^\top\|\|\nabla_{\eta\xi}\mathcal{L}(\eta^*(\xi),\xi)\|\leq \frac{\ell_{l,1}}{\mu_l}=L_\eta
\end{align*}
which implies $\eta^*(\xi)$ is $L_\eta$ Lipschitz continuous. 

Finally, we aim to prove the smoothness of $\eta^*(\xi)$. For any $\xi_1$ and $\xi_2$, we have
\begin{align*}
    &\quad\|\nabla \eta^*(\xi_1)-\nabla \eta^*(\xi_2)\|\\
    &= \|V_2\left(V_2^\top\nabla_{\eta\eta}\mathcal{L}(\eta^*(\xi_1),\xi_1)V_2\right)^{-1}V_2^\top\nabla_{\eta\xi}\mathcal{L}(\eta^*(\xi_1),\xi_1)\\
    &\quad-V_2\left(V_2^\top\nabla_{\eta\eta}\mathcal{L}(\eta^*(\xi_2),\xi_2)V_2\right)^{-1}V_2^\top\nabla_{\eta\xi}\mathcal{L}(\eta^*(\xi_2),\xi_2)\|\\
    &\leq\|V_2B_1^{-1}V_2^\top\|\|\nabla_{\eta\xi}\mathcal{L}(\eta^*(\xi_1),\xi_1)-\nabla_{\eta\xi}\mathcal{L}(\eta^*(\xi_2),\xi_2))\|\\
    &\quad+\|V_2(B_1^{-1}-B_2^{-1})V_2^\top\|\|\nabla_{\eta\xi}\mathcal{L}(\eta^*(\xi_2),\xi_2)\|\\
    &\stackrel{(a)}{\leq} \frac{1}{\mu_l}\|\nabla_{\eta\xi}\mathcal{L}(\eta^*(\xi_1),\xi_1)-\nabla_{\eta\xi}\mathcal{L}(\eta^*(\xi_2),\xi_2)\|\\
    &\quad +\frac{\ell_{l,1}}{\mu_l^2}\|\nabla_{\eta\eta}\mathcal{L}(\eta^*(\xi_1),\xi_1)-\nabla_{\eta\eta}\mathcal{L}(\eta^*(\xi_2),\xi_2)\|\\
    &\stackrel{(b)}{\leq}\frac{\ell_{l,2}(1+\frac{\ell_{l,1}}{\mu_l})^2}{\mu_l}\|\xi_1-\xi_2\|\numberthis\label{refer} 
\end{align*}
where $B_1=V_2^\top\nabla_{\eta\eta}\mathcal{L}(\eta^*(\xi_1),\xi_1)V_2$ and $B_2=V_2^\top\nabla_{\eta\eta}\mathcal{L}(\eta^*(\xi_2),\xi_2)V_2$, (a) comes from \eqref{bound_v} and the following fact:
\begin{align*}
    &\quad V_2\left(B_1^{-1}-B_2^{-1}\right)V_2^\top\\
    &=V_2B_1^{-1}\left(B_2-B_1\right)B_2^{-1}V_2^\top\\
    &=V_2B_1^{-1}\left(\left(V_2^\top\nabla_{\eta\eta}\mathcal{L}(\eta^*(\xi_2),\xi_2)V_2\right)-\left(V_2^\top\nabla_{\eta\eta}\mathcal{L}(\eta^*(\xi_1),\xi_1)V_2\right)\right)B_2^{-1}V_2^\top\\
    &=V_2B_1^{-1}V_2^\top \left(\nabla_{\eta\eta}\mathcal{L}(\eta^*(\xi_2),\xi_2)-\nabla_{\eta\eta}\mathcal{L}(\eta^*(\xi_1),\xi_1)\right)V_2B_2^{-1}V_2^\top
\end{align*}
and (b) comes from 
\begin{align*}
    \|\nabla^2 \mathcal{L}(\eta^*(\xi_1),\xi_1)-\nabla^2 \mathcal{L}(\eta^*(\xi_2),\xi_2)\|&\leq \ell_{l,2}\left[\|\xi_1-\xi_2\|+\|\eta^*(\xi_1)-\eta^*(\xi_2)\|\right]\\
    &\leq \ell_{l,2}\left(1+\frac{\ell_{l,1}}{\mu_l}\right)\|\xi_1-\xi_2\|.
\end{align*}
\end{proof}

In Algorithm \ref{alg: AGM general}, we approximate $\nabla\eta^*(\xi)$ by unrolling the differentiation. 
The following lemma investigates the error of this approximation in constrained bilevel problems for the first time, indicating that the gradient estimation error can be effectively bounded by the accuracy of the lower-level solution. 

\begin{lemma}[Error of unrolling differentiation]
Suppose that Assumption \ref{as1}--\ref{ash} hold and choose $\tau_l\leq\frac{1}{2\ell_{l,1}}$, the error of implicit gradient estimator can be bounded by
\begin{align*}
    \|\nabla\eta^*(\xi^{k_u})-\nabla_{\xi^{k_u}} \eta^{k_u+1}\|^2\leq 2\left(1-\tau_l\mu_l\right)^{2K_l+2}+2C_{K_l}C_l^2\|\eta^*(\xi^{k_u})-\eta^{k_u}\|^2
\end{align*}
where $C_l^2:=\left(1+\frac{\ell_{l,1}}{\mu_l}\right)\ell_{l,2}^2\left(\frac{2}{\mu_l^2}+\frac{3}{2\ell_{l,1}^2}\right)$ and $C_{K_l}$ is the upper bound of $K_l(1-\tau_l\mu_l)^{K_l-1}$ and is finite. 
\label{lem: lm-sub}
\end{lemma}
\begin{proof}\allowdisplaybreaks
According to \eqref{chain-rule-grad}, we know that $\nabla_{\xi^{k_u}}\eta^{k_u,1}=0$ and
\begin{align*}
    \nabla_{\xi^{k_u}} \eta^{k_u,k_l+1}
    &=V_2V_2^\top\left(I-\tau_l \nabla_{\eta\eta}\mathcal{L}(\eta^{k_u,k_l},\xi^{k_u})\right)\nabla_{\xi^{k_u}} \eta^{k_u,k_l}-\tau_l V_2 V_2^\top\nabla_{\eta\xi}\mathcal{L}(\eta^{k_u,k_l},\xi^{k_u}).
\end{align*}
For any given $\xi^{k_u}$, we can define an auxiliary sequence $\{w^{k_l}\}_{k_l=0}^\infty$ and $w^*:=\lim_{K_l\rightarrow\infty}w^{K_l}$, where $w^1=0$ and 
\begin{align}
w^{k_l+1}=V_2V_2^\top\left(I-\tau_l \nabla_{\eta\eta}\mathcal{L}(\eta^{*}(\xi^{k_u}),\xi^{k_u})\right)w^{k_l}-\tau_l V_2 V_2^\top\nabla_{\eta\xi}\mathcal{L}(\eta^{*}(\xi^{k_u}),\xi^{k_u}).\label{chain-rule-w}
\end{align}
We can see that \eqref{chain-rule-grad} and \eqref{chain-rule-w} only differ in $\eta^{*}(\xi^{k_u})$ and $\eta^{k_u,k_l}$. For the sequence $w^{k_l}$, we can calculate the explicit form of $w^{K_l+1}$ as
\begin{align*}
w^{K_l+1}
&=\sum_{s=0}^{K_l}\left(V_2V_2^\top-\tau_l V_2V_2^\top\nabla_{\eta\eta}\mathcal{L}(\eta^*(\xi^{k_u}),\xi^{k_u})\right)^s\left(-\tau_l V_2V_2^\top\nabla_{\eta\xi}\mathcal{L}(\eta^*(\xi^{k_u}),\xi^{k_u})\right)\\
&=\sum_{s=0}^{K_l}\left(V_2V_2^\top-\tau_l V_2V_2^\top\nabla_{\eta\eta}\mathcal{L}(\eta^*(\xi^{k_u}),\xi^{k_u})V_2V_2^\top\right)^s\left(-\tau_l \nabla_{\eta\xi}\mathcal{L}(\eta^*(\xi^{k_u}),\xi^{k_u})\right)\\
&=\sum_{s=0}^{K_l}\left(V_2\left(I-\tau_l V_2^\top\nabla_{\eta\eta}\mathcal{L}(\eta^*(\xi^{k_u}),\xi^{k_u})V_2\right)V_2^\top\right)^s\left(-\tau_l \nabla_{\eta\xi}\mathcal{L}(\eta^*(\xi^{k_u}),\xi^{k_u})\right)\\
&=\sum_{s=0}^{K_l}V_2\left(I-\tau_l V_2^\top\nabla_{\eta\eta}\mathcal{L}(\eta^*(\xi^{k_u}),\xi^{k_u})V_2\right)^sV_2^\top\left(-\tau_l \nabla_{\eta\xi}\mathcal{L}(\eta^*(\xi^{k_u}),\xi^{k_u})\right)\\
&=V_2\left(\sum_{s=0}^{K_l}\left(I-\tau_l V_2^\top\nabla_{\eta\eta}\mathcal{L}(\eta^*(\xi^{k_u}),\xi^{k_u})V_2\right)^s\right)V_2^\top\left(-\tau_l \nabla_{\eta\xi}\mathcal{L}(\eta^*(\xi^{k_u}),\xi^{k_u})\right)\numberthis\label{w^K+1}
\end{align*}
where the first equality comes from unrolling \eqref{chain-rule-w}, the second and the fourth equality are due to $(V_2V_2^\top)^s=V_2V_2^\top$. 
Let $D:=I-\tau_l V_2^\top\nabla_{\eta\eta}\mathcal{L}(\eta^*(\xi^{k_u}),\xi^{k_u})V_2$. When $\tau_l<\frac{2}{\ell_{l,1}}$, the operator norm of $D$ satisfies $\|D\|<1$, the limit $\sum_{s=0}^{+\infty}D^s:=\lim_{K_l\to+\infty}\sum_{s=0}^{K_l}D^s=(I-D)^{-1}=(\tau_l V_2^\top\nabla_{\eta\eta}\mathcal{L}(\eta^*(\xi^{k_u}),\xi^{k_u})V_2)^{-1}.$
Therefore, the limit point of $w^{k_l}$ is equal to $\nabla\eta^*(\xi^{k_u})$ since
\begin{align*}
w^*:=\lim_{K_l\rightarrow\infty}w^{K_l}&=V_2\left(\sum_{s=0}^{\infty}\left(I-\tau_l V_2^\top\nabla_{\eta\eta}\mathcal{L}(\eta^*(\xi^{k_u}),\xi^{k_u})V_2\right)^s\right)V_2^\top\left(-\tau_l \nabla_{\eta\xi}\mathcal{L}(\eta^*(\xi^{k_u}),\xi^{k_u})\right)\\
&=V_2\left(\tau_l V_2^\top\nabla_{\eta\eta}\mathcal{L}(\eta^*(\xi^{k_u}),\xi^{k_u})V_2\right)^{-1}V_2^\top\left(-\tau_l \nabla_{\eta\xi}\mathcal{L}(\eta^*(\xi^{k_u}),\xi^{k_u})\right)\\
&=-V_2\left(V_2^\top\nabla_{\eta\eta}\mathcal{L}(\eta^*(\xi^{k_u}),\xi^{k_u})V_2\right)^{-1}V_2^\top\nabla_{\eta\xi}\mathcal{L}(\eta^*(\xi^{k_u}),\xi^{k_u})=\nabla\eta^*(\xi^{k_u})\numberthis\label{w*-form}
\end{align*}
Moreover, the error by finite-step approximation can be bounded by
\begin{align*}
    \left\|(I-D)^{-1}-\sum_{s=0}^{K_l} D^s\right\|=\left\|\sum_{s=K_l+1}^\infty D^s\right\|\leq\sum_{s=K_l+1}^\infty \|D\|^s=\frac{\|D\|^{K_l+1}}{1-\|D\|}
\end{align*}
Since $(1-\tau_l\ell_{l,1})I\preceq D=I-\tau_l V_2^\top\nabla_{\eta\eta}\mathcal{L}(\eta^*(\xi^{k_u}),\xi^{k_u})V_2\preceq (1-\tau_l\mu_{l}) I$ and according to \eqref{w*-form} and \eqref{w^K+1}, we know that if $\tau_l\leq\frac{1}{2\ell_{l,1}}$, 
\begin{align}
    \|w^{K_l+1}-\nabla\eta^*(\xi^{k_u})\|\leq\tau_l\ell_{l,1}\frac{\left(1-\tau_l\mu_l\right)^{K_l+1}}{1-\tau_l\ell_{l,1}}\leq \left(1-\tau_l\mu_l\right)^{K_l+1} \label{ws-q}.
\end{align}

Next, we aim to bound the distance between $\nabla_{\xi^{k_u}} \eta^{k_u,k_l}$ and the auxiliary sequence $w^{k_l}$. For any $k_l$, according to \eqref{chain-rule-grad} and \eqref{chain-rule-w}, we have
\begin{align*}
    \|\nabla_{\xi^{k_u}} \eta^{k_u,k_l+1}-w^{k_l+1}\|^2
    &=\|\left(V_2 V_2^\top-\tau_l V_2 V_2^\top\nabla_{\eta\eta}\mathcal{L}(\eta^{k_u,k_l},\xi^{k_u})\right)\nabla_{\xi^{k_u}} \eta^{k_u,k_l}-\tau_l V_2 V_2^\top\nabla_{\eta\xi}\mathcal{L}(\eta^{k_u,k_l},\xi^{k_u})\\
    &\quad-\left(V_2 V_2^\top-\tau_l V_2 V_2^\top\nabla_{\eta\eta}\mathcal{L}(\eta^{*}(\xi^{k_u}),\xi^{k_u})\right)w^{k_l}+\tau_l V_2 V_2^\top\nabla_{\eta\xi}\mathcal{L}(\eta^{*}(\xi^{k_u}),\xi^{k_u})\|^2\\
    &\leq \left(1+\gamma\right)\|\left(V_2V_2^\top-\tau_l V_2 V_2^\top\nabla_{\eta\eta}\mathcal{L}(\eta^{k_u,k_l},\xi^{k_u})\right)(\nabla_{\xi^{k_u}} \eta^{k_u,k_l}-w^{k_l})\|^2+\\
    &\quad +2\left(1+\frac{1}{\gamma}\right)\|\tau_l V_2V_2^\top(\nabla_{\eta\eta}\mathcal{L}(\eta^{k_u,k_l},\xi^{k_u})-\nabla_{\eta\eta}\mathcal{L}\left(\eta^*(\xi^{k_u}),\xi^{k_u}\right)w^{k_l}\|^2\\ 
    &\quad +2\left(1+\frac{1}{\gamma}\right)\|\tau_l V_2V_2^\top(\nabla_{\eta\xi}\mathcal{L}(\eta^{k_u,k_l},\xi^{k_u})-\nabla_{\eta\xi}\mathcal{L}\left(\eta^*(\xi^{k_u}),\xi^{k_u}\right)\|^2\\
    &\leq \left(1+\gamma\right)\left(1-\tau_l\mu_l\right)^2\|\nabla_{\xi^{k_u}} \eta^{k_u,k_l}-w^{k_l}\|^2\\
    &\quad +\left(1+\frac{1}{\gamma}\right)\tau_l^2\ell_{l,2}^2\|\eta^*(\xi^{k_u})-\eta^{k_u,k_l}\|^2\left(2+2\|w^{k_l}\|^2\right)\numberthis\label{grad_aux}
\end{align*}
where the first inequality is derived from $\|a+b+c\|_2^2\leq(1+\gamma)\|a\|_2^2+(2+\frac{2}{\gamma})\|b\|_2^2+(2+\frac{2}{\gamma})\|c\|_2^2$ 
and the second inequality is due to  Assumption \ref{as1}--\ref{as2}. 
On the one hand, $\|w^{k_l}\|\leq\|\nabla\eta^*(\xi^{k_u})\|+\|w^{k_l}-\nabla\eta^*(\xi^{k_u})\|\leq\ell_{l,1}\left(\frac{1}{\mu_l}+\tau_l\right)$ is bounded according to Lemma \ref{lem: y-smooth} and \eqref{ws-q}. Thus, if $\tau_l\leq\frac{1}{2\ell_{l,1}}$ and letting $\gamma=\tau_l\mu_{l}$,  \eqref{grad_aux} becomes
\begin{align}
    \|\nabla_{\xi^{k_u}} \eta^{k_u,k_l+1}-w^{k+1}\|^2
    &\leq\left(1-\tau_l\mu_l\right)\|\nabla_{\xi^{k_u}} \eta^{k_u,k_l}-w^{k_l}\|^2 +\left(1+\frac{1}{\tau_l\mu_{l}}\right)\tau_l^2\ell_{l,2}^2\left(\frac{4\ell_{l,1}^2}{\mu_l^2}+3\right)\|\eta^*(\xi^{k_u})-\eta^{k_u,k_l}\|^2\nonumber\\
    &\leq \left(1-\tau_l\mu_l\right)\|\nabla_{\xi^{k_u}} \eta^{k_u,k_l}-w^{k_l}\|^2+C_l^2\|\eta^*(\xi^{k_u})-\eta^{k_u,k_l}\|^2\label{grad_aux1}
\end{align}
where $C_l^2:=\left(1+\frac{\ell_{l,1}}{\mu_l}\right)\ell_{l,2}^2\left(\frac{2}{\mu_l^2}+\frac{3}{2\ell_{l,1}^2}\right)$. 

On the other hand, we know that projected gradient descent is a contraction according to \cite{xiao2023alternating}, i.e.  
\begin{align}
    \|\eta^{k_u,k_l+1}-\eta^*(\xi^{k_u})\|^2\leq\left(1-\tau_l\mu_l\right)\|\eta^{k_u,k_l}-\eta^*(\xi^{k_u})\|^2\label{y^*}
\end{align}
for $0\leq\tau_l\leq\frac{1}{\ell_{l,1}}$. By induction, we have
\begin{align}
    \|\eta^{k_u,k_l+1}-\eta^*(\xi^{k_u})\|^2\leq \left(1-\tau_l\mu_l\right)^{k_l}\|\eta^{k_u}-\eta^*(\xi^{k_u})\|^2
    \label{lower-contr}
\end{align}
Then \eqref{grad_aux1} becomes
\begin{align}
    \|\nabla_{\xi^{k_u}} \eta^{k_u,k_l+1}-w^{k_l+1}\|^2&\leq\left(1-\tau_l\mu_l\right)\|\nabla_{\xi^{k_u}} \eta^{k_u,k_l}-w^{k_l}\|^2+C_l^2\left(1-\tau_l\mu_l\right)^{k_l-1}\|\eta^{k_u}-\eta^*(\xi^{k_u})\|^2.
\end{align}
Then by induction and $w^1=\nabla_{\xi^{k_u}} \eta^{k_u,1}=0,\eta^{k_u+1}=\eta^{k_u,K_l+1}$, we obtain that
\begin{align}
    \|\nabla_{\xi^{k_u}} \eta^{k_u+1}-w^{K_l+1}\|^2&\leq K_l(1-\tau_l\mu_l)^{K_l-1}C_l^2\|\eta^*(\xi^{k_u})-\eta^{k_u}\|^2.
\label{grad_aux3}
\end{align}

Combining \eqref{grad_aux3} with \eqref{ws-q} and setting $\tau_l\leq\frac{1}{2\ell_{l,1}}$, we know that 
\begin{align}
    \|\nabla_{\xi^{k_u}} \eta^{k_u+1}-\nabla\eta^*(\xi^{k_u})\|^2\leq 2\left(1-\tau_l\mu_l\right)^{2K_l+2}+2K_l(1-\tau_l\mu_l)^{K_l-1}C_l^2\|\eta^*(\xi^{k_u})-\eta^{k_u}\|^2. \label{gradient-diff}
\end{align}
Then given $\tau_l$ and let $f(K_l)=K_l(1-\tau_l\mu_l)^{K_l-1}$, we know $\log(f(K_l))=\log K_l+(K_l-1)\log (1-\tau_l\mu_l)$. Taking the gradient of $\log(f(K_l))$, we get
$1/K_l+\log (1-\tau_l\mu_l)$. As $\log (1-\tau_l\mu_l)<0$, we know $\log(f(K_l))$ first increases and then decreases and thus, $\log(f(K_l))$ and $f(K_l)$ have a finite upper bound. Let us denote the upper bound of $K_l(1-\tau_l\mu_l)^{K_l-1}$ as $C_{K_l}={\cal O}(1)$. Then \eqref{gradient-diff} becomes 
\begin{align}
    \|\nabla_{\xi^{k_u}} \eta^{k_u+1}-\nabla\eta^*(\xi^{k_u})\|^2\leq 2\left(1-\tau_l\mu_l\right)^{2K_l+2}+2C_{K_l}C_l^2\|\eta^*(\xi^{k_u})-\eta^{k_u}\|^2. 
\end{align}
which yields the conclusion. 
\end{proof}

Besides, we have the lower-level contraction and error. 

\begin{lemma}[Lower-level error]
Suppose that Assumption \ref{as1}--\ref{ash} hold and $\tau_l\leq\frac{1}{\ell_{l,1}}$, then for any $\gamma>0$, we have
\begin{subequations}
\begin{align}
    &\|\eta^{k_u+1}-\eta^*(\xi^{k_u})\|^2\leq\left(1-\tau_l\mu_l\right)^{K_l}\|\eta^{k_u}-\eta^*(\xi^{k_u})\|^2\label{27a}\\
    &\|\eta^{k_u+1}-\eta^*(\xi^{k_u+1})\|^2\leq (1+\gamma)\|\eta^{k_u+1}-\eta^*(\xi^{k_u})\|^2+L_\eta^2\left(1+\frac{1}{\gamma}\right)\|\xi^{k_u}-\operatorname{Proj}_{\Xi}(\xi^{k_u}- \tau_u\widehat\nabla u(\xi^{k_u}))\|^2\label{27b}
\end{align}
\end{subequations}
\label{lm-ll}
\end{lemma}

\begin{proof}
\eqref{27a} comes from \eqref{lower-contr} when setting $k_l=K_l$. Moreover, 
\begin{align*}
    \|\eta^{k_u+1}-\eta^*(\xi^{k_u+1})\|^2&=\|\eta^{k_u+1}-\eta^*(\xi^{k_u})+\eta^*(\xi^{k_u})-\eta^*(\xi^{k_u+1})\|^2\\
    &\stackrel{(a)}{\leq}(1+\gamma)\|\eta^{k_u+1}-\eta^*(\xi^{k_u})\|^2+\left(1+\frac{1}{\gamma}\right)\|\eta^*(\xi^{k_u})-\eta^*(\xi^{k_u+1})\|^2\\
    &\stackrel{(b)}{\leq}(1+\gamma)\|\eta^{k_u+1}-\eta^*(\xi^{k_u})\|^2+L_\eta^2\left(1+\frac{1}{\gamma}\right)\|\xi^{k_u}-\xi^{k_u+1}\|^2\\
    &= (1+\gamma)\|\eta^{k_u+1}-\eta^*(\xi^{k_u})\|^2+L_\eta^2\left(1+\frac{1}{\gamma}\right)\|\xi^{k_u}-\operatorname{Proj}_{\Xi}(\xi^{k_u}-\tau_u \widehat\nabla u(\xi^{k_u}))\|^2
\end{align*}
where (a) is due to $\|a+b\|_2^2\leq (1+\gamma)\|a\|_2^2+(1+\frac{1}{\gamma})\|b\|_2^2$ for any $\gamma>0$, 
and (b) comes from the Lipschitz continuity of $\eta^*(\xi)$ in Lemma \ref{lem: y-smooth}. 
\end{proof}

\begin{lemma}[Upper-level error] Under Suppose that Assumption \ref{as1}--\ref{ash} hold and $\tau_l\leq\frac{1}{2\ell_{l,1}}$, then it holds that
\begin{align*}
u(\xi^{k_u+1})-u(\xi^{k_u}) &\leq -\frac{\tau_u}{2}\|\xi^{k_u}-\operatorname{Proj}_{\Xi}(\xi^{k_u}-\nabla u(\xi^{k_u}))\|^2-\left(\frac{1}{2\tau_u}-\frac{L_u}{2}\right) \|\xi^{k_u}-\operatorname{Proj}_{\Xi}(\xi^{k_u}-\widehat\nabla u(\xi^{k_u}))\|^2\\
    &\quad+\tau_u\left(\ell_{u,1}\left(1+L_\eta\right)+2\ell_{u,0}C_{K_l}C_l^2\right)^2\|\eta^*(\xi^{k_u})-\eta^{k_u}\|^2+2\tau_u\left(1-\tau_l\mu_l\right)^{4K_l+4}
\end{align*}
\label{lm:UL}
\end{lemma}

\begin{proof}\allowdisplaybreaks
According to Lemma \ref{lem: y-smooth}, we know $u(\xi)=\mathcal{U}(\eta^*(\xi),\xi)$ is Lipschitz smooth and
\begin{align*}
\nabla u(\xi)=\nabla_\xi\mathcal{U}(\eta^*(\xi),\xi)+\nabla_\xi^\top \eta^*(\xi)\nabla_\eta\mathcal{U}(\eta^*(\xi),\xi)
\end{align*}
and for any $\xi_1,\xi_2$, we have
\begin{align*}
\|\nabla u(\xi_1)-\nabla u(\xi_2)\|&= \|\nabla_\xi\mathcal{U}(\eta^*(\xi_1),\xi_1)+\nabla_\xi^\top \eta^*(\xi_1)\nabla_\eta\mathcal{U}(\eta^*(\xi_1),\xi_1)-\nabla_\xi\mathcal{U}(\eta^*(\xi_2),\xi_2)-\nabla_\xi^\top \eta^*(\xi_2)\nabla_\eta\mathcal{U}(\eta^*(\xi_2),\xi_2)\|\\
&\leq \|\nabla_\xi\mathcal{U}(\eta^*(\xi_1),\xi_1)-\nabla_\xi\mathcal{U}(\eta^*(\xi_2),\xi_2)\|+\|\nabla_\xi^\top \eta^*(\xi_1)\|\|\nabla_\eta\mathcal{U}(\eta^*(\xi_1),\xi_1)-\nabla_\eta\mathcal{U}(\eta^*(\xi_2),\xi_2)\|\\
&\quad+\|\nabla_\eta\mathcal{U}(\eta^*(\xi_2),\xi_2)\|\|\nabla_\xi \eta^*(\xi_1)-\nabla_\xi \eta^*(\xi_2)\|\\
&\leq \ell_{u,1}(\|\eta^*(\xi_1)-\eta^*(\xi_2)\|+\|\xi_1-\xi_2\|)+L_\eta \ell_{u,1}(\|\eta^*(\xi_1)-\eta^*(\xi_2)\|+\|\xi_1-\xi_2\|)+\ell_{u,0}L_{\eta\xi}\|\xi_1-\xi_2\|\\
&\leq (\ell_{u,1}(1+L_\eta)^2+\ell_{u,0}L_{\eta\xi})\|\xi_1-\xi_2\|.
\end{align*}
By denoting the smoothness constant of $u(\xi)$ as $L_u:=\ell_{u,1}(1+L\eta)^2+\ell_{u,0}L_{\eta\xi}$, we have the following expansion 
\begin{align*}
    u(\xi^{k_u+1})
    &\leq u(\xi^{k_u})+\langle \nabla u(\xi^{k_u}),\xi^{k_u+1}-\xi^{k_u}\rangle+\frac{L_u}{2}\|\xi^{k_u+1}-\xi^{k_u}\|^2\\
    &= u(\xi^{k_u})-\langle \nabla u(\xi^{k_u}),\xi^{k_u}-\operatorname{Proj}_{\Xi}(\xi^{k_u}-\tau_u\widehat\nabla u(\xi^{k_u}))\rangle+\frac{L_u}{2}\|\xi^{k_u}-\operatorname{Proj}_{\Xi}(\xi^{k_u}-\tau_u\widehat\nabla u(\xi^{k_u}))\|^2\\
    &\stackrel{(a)}{=} u(\xi^{k_u})-\frac{1}{\tau_u}\langle \xi^{k_u}-\operatorname{Proj}_{\Xi}(\xi^{k_u}-\tau_u\nabla u(\xi^{k_u})),\xi^{k_u}-\operatorname{Proj}_{\Xi}(\xi^{k_u}-\tau_u\widehat\nabla u(\xi^{k_u}))\rangle\\
    &\quad +\frac{L_u}{2}\|\xi^{k_u}-\operatorname{Proj}_{\Xi}(\xi^{k_u}-\tau_u\widehat\nabla u(\xi^{k_u}))\|^2\\
    &\stackrel{(b)}{=} u(\xi^{k_u})-\frac{1}{2\tau_u}\|\xi^{k_u}-\operatorname{Proj}_{\Xi}(\xi^{k_u}-\tau_u\nabla u(\xi^{k_u}))\|^2\\
    &\quad +\frac{1}{2\tau_u}\|\operatorname{Proj}_{\Xi}(\xi^{k_u}-\tau_u\nabla u(\xi^{k_u}))-\operatorname{Proj}_{\Xi}(\xi^{k_u}-\tau_u\widehat\nabla u(\xi^{k_u}))\|^2\\
    &\quad -\left(\frac{1}{2\tau_u}-\frac{L_u}{2}\right) \|\xi^{k_u}-\operatorname{Proj}_{\Xi}(\xi^{k_u}-\tau_u\widehat\nabla u(\xi^{k_u}))\|^2\\
    &\stackrel{(c)}{\leq} u(\xi^{k_u})-\frac{\tau_u}{2}\|\xi^{k_u}-\operatorname{Proj}_{\Xi}(\xi^{k_u}-\nabla u(\xi^{k_u}))\|^2+\frac{\tau_u}{2}\|\nabla u(\xi^{k_u})-\widehat\nabla u(\xi^{k_u})\|^2\\
    &\quad -\left(\frac{1}{2\tau_u}-\frac{L_u}{2}\right) \|\xi^{k_u}-\operatorname{Proj}_{\Xi}(\xi^{k_u}-\tau_u\widehat\nabla u(\xi^{k_u}))\|^2\numberthis\label{upper}
\end{align*}
where (a) comes from $\xi^{k_u}=\operatorname{Proj}_{\Xi}(\xi^{k_u})$ and the fact that $\operatorname{Proj}_{\Xi}$ onto a linear equality constraint set is a linear operator, (b) is derived from $2a^\top b= \|a\|^2+\|b\|^2-\|a-b\|^2$ and (c) is because $\operatorname{Proj}_{\Xi}$ is a linear operator and $\|\operatorname{Proj}(A)-\operatorname{Proj}(B)\|\leq\|A-B\|$. 
Besides, we can decompose the gradient bias term as follows 
\begin{align*}
    \|\nabla u(\xi^{k_u})-\widehat\nabla u(\xi^{k_u})\|&=\|\nabla_\xi \mathcal{U}(\eta^*(\xi^{k_u}),\xi^{k_u})-\nabla\eta^*(\xi^{k_u})^\top\nabla_\eta \mathcal{U}(\eta^*(\xi^{k_u}),\xi^{k_u})\\
    &\quad -\nabla_\xi \mathcal{U}(\eta^{k_u+1},\xi^{k_u})+\nabla_{\xi^{k_u}}^\top \eta^{k_u+1}\nabla_\eta \mathcal{U}(\eta^{k_u+1},\xi^{k_u})\|\\
    &\leq\|\nabla_\xi \mathcal{U}(\eta^*(\xi^{k_u}),\xi^{k_u})-\nabla_\xi \mathcal{U}(\eta^{k_u+1},\xi^{k_u})\|\\
    &\quad +\|\nabla\eta^*(\xi^{k_u})\|\|\nabla_\eta \mathcal{U}(\eta^*(\xi^{k_u}),\xi^{k_u})-\nabla_\eta \mathcal{U}(\eta^{k_u+1},\xi^{k_u})\|\\
    &\quad +\|\nabla_\eta\mathcal{U}(\eta^{k_u+1},\xi^{k_u})\|\|\nabla\eta^*(\xi^{k_u})-\nabla_{\xi^{k_u}} \eta^{k_u+1}\|\\
    &\leq \ell_{u,1}\left(1+L_\eta\right)\|\eta^*(\xi^{k_u})-\eta^{k_u+1}\|+\ell_{u,0}\|\nabla\eta^*(\xi^{k_u})-\nabla_{\xi^{k_u}} \eta^{k_u+1}\|\\
    &\stackrel{(a)}{\leq} \left(\ell_{u,1}\left(1+L_\eta\right)+2\ell_{u,0}C_{K_l}C_l^2\right)\|\eta^*(\xi^{k_u})-\eta^{k_u}\|+2\left(1-\tau_l\mu_l\right)^{2K_l+2}\numberthis\label{grad_bias}
\end{align*}
where (a) comes from lower-level contraction \eqref{lower-contr}. 
Thus, plugging \eqref{grad_bias} to \eqref{upper}, we get that
\begin{align*}
    u(\xi^{k_u+1})-u(\xi^{k_u})&\leq-\frac{\tau_u}{2}\|\xi^{k_u}-\operatorname{Proj}_{\Xi}(\xi^{k_u}-\nabla u(\xi^{k_u}))\|^2-\left(\frac{1}{2\tau_u}-\frac{L_u}{2}\right) \|\xi^{k_u}-\operatorname{Proj}_{\Xi}(\xi^{k_u}-\widehat\nabla u(\xi^{k_u}))\|^2\\
    &\quad+\frac{\tau_u}{2}\left(\left(\ell_{u,1}\left(1+L_\eta\right)+2\ell_{u,0}C_{K_l}C_l^2\right)\|\eta^*(\xi^{k_u})-\eta^{k_u}\|+2\left(1-\tau_l\mu_l\right)^{2K+2}\right)^2\\
    &\leq -\frac{\tau_u}{2}\|\xi^{k_u}-\operatorname{Proj}_{\Xi}(\xi^{k_u}-\nabla u(\xi^{k_u}))\|^2-\left(\frac{1}{2\tau_u}-\frac{L_u}{2}\right) \|\xi^{k_u}-\operatorname{Proj}_{\Xi}(\xi^{k_u}-\widehat\nabla u(\xi^{k_u}))\|^2\\
    &\quad+\tau_u\left(\ell_{u,1}\left(1+L_\eta\right)+2\ell_{u,0}C_{K_l}C_l^2\right)^2\|\eta^*(\xi^{k_u})-\eta^{k_u}\|^2+2\tau_u\left(1-\tau_l\mu_l\right)^{4K_l+4}
\end{align*}
\end{proof}

With Lemma \ref{lem: y-smooth}, \ref{lem: lm-sub}, \ref{lm-ll} and \ref{lm:UL}, we restate the convergence Theorem \ref{thm-conv} in a more formal way and prove the theorem as follows.
\begin{theorem}
Under Assumption \ref{as1}--\ref{ash}, 
let $\tau_l\leq\frac{1}{2\ell_{l,1}}, K_l={\cal O}(\log K_u)$ and $\tau_u={\cal O}(1)$ satisfies
\begin{align*}
&\tau_u\leq\min\left\{\frac{1}{2L_u(1+2L_\eta)}, \frac{\tau_l L_u\mu_l}{L_\eta\left((\ell_{u,1}(1+L_\eta)+2\ell_{u,0}C_{K_l}C_l^2)^2+4L_u^2\right)}\right\}, 
\end{align*}
then the iterates of Algorithm \ref{alg: AGM general} satisfy
\begin{align}
    \frac{1}{K_u}\sum_{k_u=1}^{K_u}\|\xi^{k_u}-\operatorname{Proj}_{\Xi}(\xi^{k_u}-\nabla u(\xi^{k_u}))\|^2={\cal O}\left(\frac{1}{K_u}\right)
\end{align}
where ${\mathcal O}$ omits the $log$ dependency. 
\end{theorem}
\begin{proof}
We can define Lyapunov function as $$\mathbb{V}^{k_u}\:= u(\xi^{k_u})+\frac{ L_u}{L_\eta}\|\eta^*(\xi^{k_u})-\eta^{k_u}\|^2$$ 

On the one hand, plugging \eqref{27a} to \eqref{27b}, we get 
\begin{align}
    \|\eta^{k_u+1}-\eta^*(\xi^{k_u+1})\|^2\leq (1+\gamma)\left(1-\tau_l\mu_l\right)\|\eta^{k_u}-\eta^*(\xi^{k_u})\|^2+L_\eta^2\left(1+\frac{1}{\gamma}\right)\|\xi^{k_u}-\operatorname{Proj}_{\Xi}(\xi^{k_u}- \tau_u\widehat\nabla u(\xi^{k_u}))\|^2. \label{eq:LLgeneral}
\end{align}

On the other hand, according to Lemma \ref{lm:UL} and \eqref{eq:LLgeneral}, it holds that 
\begin{align*}
\mathbb{V}^{k_u+1}-\mathbb{V}^{k_u}&\leq -\frac{\tau_u}{2}\|\xi^{k_u}-\operatorname{Proj}_{\Xi}(\xi^{k_u}-\nabla u(\xi^{k_u}))\|^2-\left(\frac{1}{2\tau_u}-\frac{L_u}{2}\right) \|\xi^{k_u}-\operatorname{Proj}_{\Xi}(\xi^{k_u}-\widehat\nabla u(\xi^{k_u}))\|^2\\
    &\quad+\tau_u\left(\ell_{u,1}\left(1+L_\eta\right)+2\ell_{u,0}C_{K_l}C_l^2\right)^2\|\eta^*(\xi^{k_u})-\eta^{k_u}\|^2+2\tau_u\left(1-\tau_l\mu_l\right)^{4K+4}\\
    &\quad+\frac{ L_u}{L_\eta}\left[(1+\gamma)\left(1-\tau_l\mu_l\right)-1\right]\|\eta^{k_u}-\eta^*(\xi^{k_u})\|^2+L_u L_\eta\left(1+\frac{1}{\gamma}\right)\|\xi^{k_u}-\operatorname{Proj}_{\Xi}(\xi^{k_u}- \tau_u\widehat\nabla u(\xi^{k_u}))\|^2\\
    &\stackrel{(a)}{\leq} -\frac{\tau_u}{2}\|\xi^{k_u}-\operatorname{Proj}_{\Xi}(\xi^{k_u}-\nabla u(\xi^{k_u}))\|^2-\left(\frac{1}{4\tau_u}-\frac{L_u}{2}-L_u L_\eta\right) \|\xi^{k_u}-\operatorname{Proj}_{\Xi}(\xi^{k_u}-\tau_u\widehat\nabla u(\xi^{k_u}))\|^2\\    &\quad-\left(\frac{\tau_l L_u\mu_l}{L_\eta}-\tau_u\left((\ell_{u,1}\left(1+L_\eta\right)+2\ell_{u,0}C_{K_l}C_l^2)^2+4L_u^2\right)\right)\|\eta^*(\xi^{k_u})-\eta^{k_u}\|^2+2\tau_u\left(1-\tau_l\mu_l\right)^{4K_l+4}\\
    &\stackrel{(b)}{\leq} -\frac{\tau_u}{2}\|\xi^{k_u}-\operatorname{Proj}_{\Xi}(\xi^{k_u}-\nabla u(\xi^{k_u}))\|^2+2\tau_u\left(1-\tau_l\mu_l\right)^{4K_l+4} \numberthis\label{V-difference}
\end{align*}
where (a) is earned by setting $\gamma=4L_u L_\eta\tau_u$ and (b) comes from the conditions 
\begin{align}
&\frac{1}{4\tau_u}-\frac{L_u}{2}-L_u L_\eta\geq 0, \quad \text{ and } \quad \frac{\tau_l L_u\mu_l}{L_\eta}-\tau_u\left((\ell_{u,1}\left(1+L_\eta\right)+2\ell_{u,0}C_{K_l}C_l^2)^2+4L_u^2\right)\geq 0. \label{condi-nonnegative}
\end{align} 
The sufficient conditions for \eqref{condi-nonnegative} are 
\begin{align*}
&\tau_u\leq\min\left\{\frac{1}{2L_u(1+2L_\eta)}, \frac{\tau_l L_u\mu_l}{L_\eta\left((\ell_{u,1}(1+L_\eta)+2\ell_{u,0}C_{K_l}C_l^2)^2+4L_u^2\right)}\right\}.  
\end{align*} 
Rearranging terms and telescoping \eqref{V-difference} yield 
\begin{align*}
\frac{1}{K_u}\sum_{k_u=1}^{K_u} \left\|\xi^{k_u}-\operatorname{Proj}_{\Xi}(\xi^{k_u}-\nabla u(\xi^{k_u}))\right\|^2 &\leq \frac{2(\mathbb{V}^{1}-\mathbb{V}^{K_u+1})}{\tau_u K_u}+4(1-\tau_l\mu_l)^{4K_l+4}\\
&\leq \frac{2(\mathbb{V}^{1}-\inf_{\xi} u(\xi))}{\tau_u K_u}+4(1-\tau_l\mu_l)^{4K_l+4}. 
\end{align*}
Then by choosing $K_l={\cal O}(\log(K_u))$, the convergence rate of Algorithm \ref{alg: AGM general} is ${\cal O}\left(\frac{1}{K_u}\right)$. 
\end{proof}

The above theorem guarantees the Algorithm \ref{alg: AGM general} converges to an $\epsilon$ stationary point given that assumptions \ref{as1}--\ref{ash} are satisfied.
And Lemma \ref{lem: convexity} states the convexity and Lipschitz smoothness of the lower-level objective functions $\calL_{\calG}(\rho,\vm;\metric,\obs)$ in \eqref{eq: disct forward mfg obj} and shows that our problem setting satisfies the assumptions.

Since the interpolation operators $\Ix,\Iy,\It$  are linear and positive definite, to prove Lemma \ref{lem: convexity}, it is sufficient to prove the (strong) convexity and the Lipschitz smoothness of $L_{G,\weight}:\bbR^+\times\bbR^d\to\bbR,(\alpha,\vb)\mapsto\frac{\vb^\top G\vb}{2\alpha}+\weight\alpha\log(\alpha)$.

\begin{lemma}
    Let $G$ be a $d\times d$ symmetric positive definite matrix and  $L_{G,\weight}:\bbR^+\times\bbR^d\to\bbR,(\alpha,\vb)\mapsto\frac{\vb^\top G\vb}{2\alpha}+\weight\alpha\log(\alpha)$.
    For any $\weight\geq0$, $L_{G,\weight}$ is convex in $\bbR^+\times\bbR^d$ 
    and Lipschitz smooth in $\{\alpha\in\bbR:\alpha\geq\underline{c}_\rho>0\}\times\{\vb\in\bbR^d:\|\vb\|\leq\overline{c}_m\}$ ($\overline{c}_m>0$). 
    And for any $\weight>0$, $L_{G,\weight}$ is strongly convex in $\{\alpha\in\bbR:\overline{c}_{\rho}\geq\alpha\geq\underline{c}_\rho>0\}\times\bbR^d$.
\label{lem: dyn cvx, entropy, strongly cvx}
\end{lemma}

\begin{proof}
    Since $G$ is symmetric and positive definite, we write the singular value decomposition of $G$ as $G=U\Sigma_G U^\top$, with $UU^\top=U^\top U=I$, $\Sigma_G=\diag(\sigma_{G,d},\sigma_{G,d-1},\cdots,\sigma_{G,1})$, $ (\sigma_{G,d}\geq\sigma_{G,d-1}\geq\cdots\geq\sigma_{G,1})$. 
    And $\sigma_{G,i},i=1,\cdots,d$ are the singular values of $G$.
    Denote $\Sigma_G^{\half}:=\diag(\sqrt{\sigma_{G,d}},\sqrt{\sigma_{G,d-1}},\cdots,\sqrt{\sigma_{G,1}})$ and $S=\Sigma_G^{\half}U^\top$. Then $G=S^\top S$ and the singular values of $S$ are $\sigma_{S,i}=\sqrt{\sigma_{G,i}}$.
    
    Obviously, $L_{G,\weight}$ is twice differentiable in $\bbR^+\times\bbR^d$    and 
    \begin{equation}
    \begin{aligned}
        \nabla^2L_{G,\weight}(\alpha,\vb) &= \frac{1}{\alpha^3}
    \left[ \begin{matrix} \vb^\top G\vb & -\alpha \vb^\top G^\top\\
        -\alpha G\vb & \alpha^2 G\end{matrix} \right]\\
    &=\frac{1}{\alpha^3}  
    \left[ \begin{matrix} 1 & \\
        & S^\top \end{matrix} \right]
    \left[ \begin{matrix} (S\vb)^\top S\vb +\weight\alpha^2 & -\alpha (S\vb)^\top \\
        -\alpha S\vb & \alpha^2 I\end{matrix} \right]
    \left[ \begin{matrix} 1 & \\
        & S \end{matrix} \right]\\
    &=\frac{1}{\alpha^3}  
    \left[ \begin{matrix} 1 & \\
        & S^\top \end{matrix} \right]
    \nabla^2L_{I,\weight}(\alpha,S\vb)
    \left[ \begin{matrix} 1 & \\
        & S \end{matrix} \right].        
    \end{aligned}
    \end{equation}
    We denote the minimal and maximal singular values of $\nabla^2L_{G,\weight}(\alpha,\vb)$ as $\sigma_{\min}^{G,\weight}(\alpha,\vb)$ and $\sigma_{\max}^{G,\weight}(\alpha,\vb)$.
    Then we have
    \begin{equation}
        \sigma_{\min}^{G,\weight}(\alpha,\vb)\geq \frac{\min(1,\sigma_{G,1})}{\alpha^3}\sigma_{\min}^{I,\weight}(\alpha,S\vb),\quad
        \sigma_{\max}^{G,\weight}(\alpha,\vb)\leq \frac{\max(1,\sigma_{G,d})}{\alpha^3}\sigma_{\max}^{I,\weight}(\alpha,S\vb).
    \end{equation}
    
    By computation, the eigenvalues $\lambda^{I,\weight}(\alpha,\vb)$ of 
    $\nabla^2L_{I,\weight}(\alpha,\vb)$
    satisfy 
    \begin{equation}
        \left(\lambda^2 -\left(\|\vb\|^2+(\weight+1)\alpha^2\right)\lambda + \weight\alpha^4\right)(\lambda-\alpha^2)^{d-1}=0.
    \end{equation}
    Therefore $\lambda^{I,\weight}(\alpha,\vb)\geq0$ and 
    \begin{equation}\left\{
    \begin{aligned}
        &\sigma_{\min}^{I,\weight}(\alpha,\vb) \geq 
        \frac{\weight\alpha^4}{\|\vb\|^2+(\weight+1)\alpha^2}\geq0,\\ 
        &\sigma_{\max}^{I,\weight}(\alpha,\vb) \leq \|\vb\|^2+(\weight+1)\alpha^2
    \end{aligned}\right.
    \end{equation}
    
    For $\weight\geq0$, 
    $\sigma_{\min}^{G,\weight}(\alpha,\vb)\geq 0$ 
    hold for any$\alpha>0,\vb\in\bbR^d$, which implies $L_\weight$ is convex. 
    For $\weight\geq0,\alpha\geq\underline{c}_{\rho},\|\vb\|\leq\overline{c}_m$, 
    $\sigma_{\max}^{G,\weight}(\alpha,\vb)\leq \max(1,\sigma_{G,d}) \left(\frac{\sigma_{G,d}\overline{c}^2_m + (\weight+1)\underline{c}_{\rho}^2}{\underline{c}^3_{\rho}}\right)$ 
    hold for any $\alpha\geq\underline{c}_{\rho},\|\vb\|\leq\overline{c}_m$, which implies $L_{G,\weight}$ is Lipschitz smooth.
    And for $\weight>0,\overline{c}_{\rho}\geq\alpha\geq\underline{c}_{\rho}$, 
    $\sigma_{\min}^{G,\weight}(\alpha,\vb) \geq \min(1,\sigma_{G,1})\min\left( \frac{\weight\overline{c}_\rho}{\sigma_{G,d}\overline{c}_m^2+(\weight+1)\overline{c}_{\rho}^2}, \frac{\weight\underline{c}_\rho}{\sigma_{G,d}\overline{c}_m^2+(\weight+1)\underline{c}_{\rho}^2} \right)$.
\end{proof}

\section{Numerical Experiments}
\label{sec: num egs}

\subsection{Experiment Settings}

This section presents several numerical experiments to illustrate the effectiveness of our model and algorithm.
We generate the data by solving the forward problem using the projected gradient descent algorithm proposed in \cite{yu2021fast} based on the FISTA algorithm~\cite{beck2009fast}. 
In each experiment, we report the relative error versus the number of iterations for recovering the obstacle and the metric.
The relative error for recovering the obstacle is
\begin{equation}
    \sqrt{ \frac{\sum_{\idx=1}^\nx\sum_{\idy=1}^\ny ((\obs^{K_u})_{\idx,\idy}-(\obstilde)_{\idx,\idy})^2}{\sum_{\idx=1}^\nx\sum_{\idy=1}^\ny (\obstilde)_{\idx,\idy}^2 }},
\end{equation}
and the relative errors for the metrics are
\begin{equation}
    \text{ (1D) } \sqrt{ \frac{\sum_{\idx=1}^\nx ((\metric^{K_u})_{\idx}-(\widetilde{\metric})_{\idx})^2}{\sum_{\idx=1}^\nx (\widetilde{\metric})_{\idx}^2 }}, \qquad 
\text{ (2D) }\sqrt{ \frac{\sum_{\idx=1}^\nx\sum_{\idy=1}^\ny \|(\metric^{K_u})_{\idx,\idy}-(\widetilde{\metric})_{\idx,\idy}\|_F^2}{\sum_{\idx=1}^\nx\sum_{\idy=1}^\ny \|(\widetilde{\metric})_{\idx,\idy}\|_F^2 }},
\end{equation}
where $\obs^{K_u},\metric^{K_u}$ are the numerical results after $K_u$ upper-level updates and $\obstilde,\widetilde{\metric}$ are the ground truth.
We implement all of our numerical experiments in Matlab on a PC with an Intel(R) i7-8550U 1.80GHz CPU  and 16 GB memory.

\subsection{Theoretical arguments verification}

\subsubsection{Algorithm convergence and obstacle unique identifiability}
The first experiment aims to numerically verify the stability Theorem \ref{thm: differentiability}, the unique identifiability Theorem \ref{thm: unique identifiablity} and convergence analysis in Theorem \ref{thm-conv} of the bilevel algorithm with lower and upper-level constraints.

We discretize the space with $\nt=16,\nx=\ny=64$. Denote $p_g(x,y;\mu_x,\mu_y,\sigma_x,\sigma_y)$ as the probability density function of Gaussian distribution with mean $(\mu_x,\mu_y)$ and covariance matrix $\diag(\sigma_x^2,\sigma_y^2)$.
We feed the model with one pair of observations, i.e. $N=1$, with $\rhoinit=p_g(\cdot,\cdot;-0.25,0,0.08,0.08)$, $\rhoend=p_g(\cdot,\cdot;0.25,0,0.08,0.08)$ and $\weight_I=0.1,\weight_T=5$.
We choose the obstacle function as $\obs(x,y)=\weight_b p_g(x,y;0,0,0.08,0.1)$.
With different values of $\weight_b$, the agents avoid the center of the obstacle to different degrees. Higher values of $\weight_b$ lead to lower density values at $(x,y)=(0,0)$.
According to remark \ref{rem: robustness reasoning}, low-density values in the data result in difficulties in accurately reconstructing the obstacle.

\begin{figure}[h!]
\centering
\includegraphics[width=3.9cm]{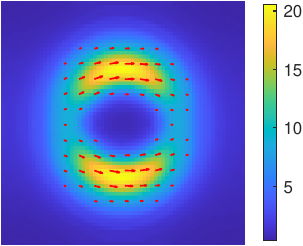}
\includegraphics[width=3.9cm]{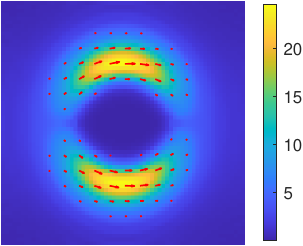}
\includegraphics[width=3.9cm]{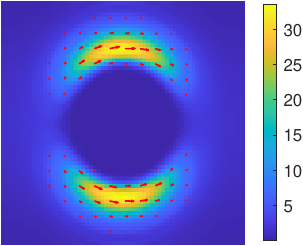}
\includegraphics[width=3.9cm]{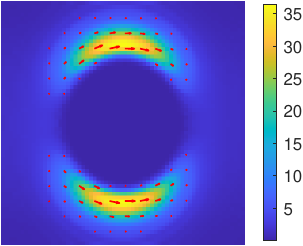}
\\
\includegraphics[width=3.9cm]{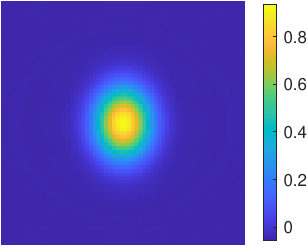}
\includegraphics[width=3.9cm]{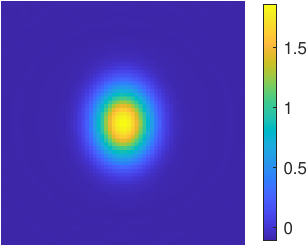}
\includegraphics[width=3.9cm]{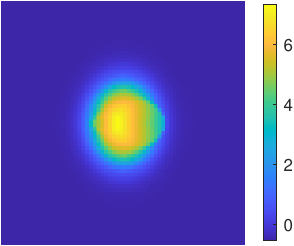}
\includegraphics[width=3.9cm]{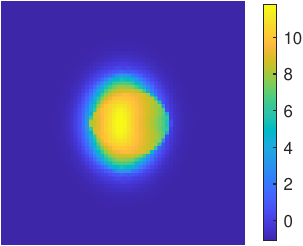}
\\
\includegraphics[width=3.9cm]{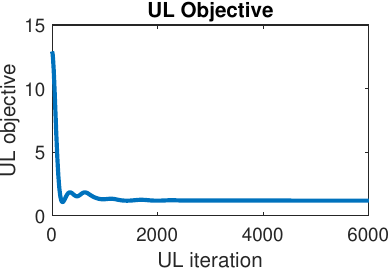}
\includegraphics[width=3.9cm]{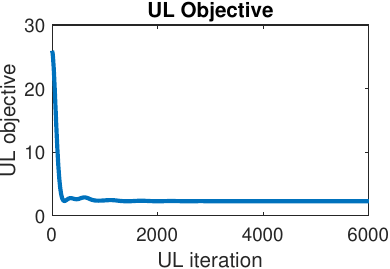}
\includegraphics[width=3.9cm]{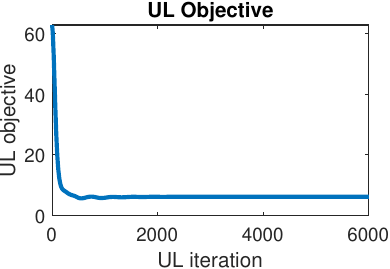}
\includegraphics[width=3.9cm]{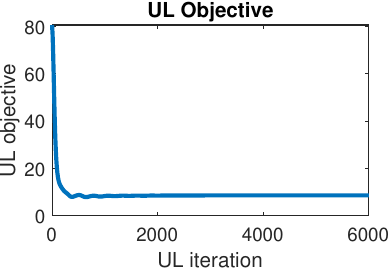}
\\
\includegraphics[width=3.9cm]{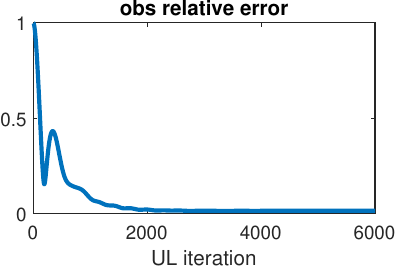}
\includegraphics[width=3.9cm]{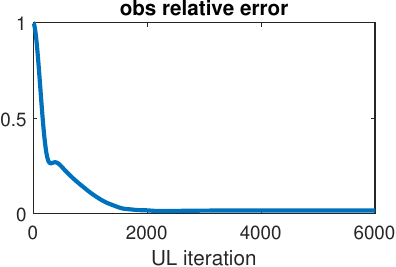}
\includegraphics[width=3.9cm]{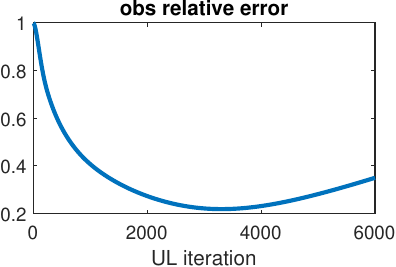}
\includegraphics[width=3.9cm]{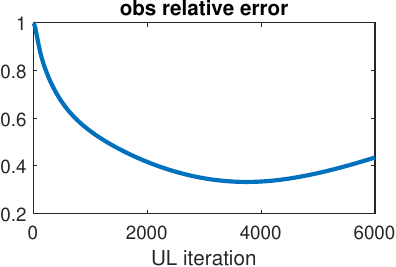}
\caption{Convergence test of the inverse crowd motion problem. Top to bottom: the snapshot of $\rhotilde$ at $t=0.5$, recovered $\obs$ with smallest relative error, upper-level objective value versus the number of iterations, relative error of $\obs$ versus the number of iterations. Left to right: $\weight_b=0.05,0.1,0.5,1$. }
\label{fig: obs gaussian}
\end{figure}

\begin{table}[htb]
\centering
\caption{Convergence test of the inverse crowd motion problem.}
\begin{tabular}{l|ccccc}
\toprule
    $\weight_b$ & $\rhotilde(0,0,0.5)$ &  \makecell[c]{upper-level\\ objective value} & \makecell[c]{relative error\\ (best)} & \makecell[c]{relative error \\(terminial)} & \makecell[c]{time elapsed \\(second)}\\ \midrule
    0.05 & 0.7831 & 1.1792 & 0.0139 & 0.0148 & 1570.1611\\
    0.1  & 0.0293 & 2.2504 & 0.0134 & 0.0161 & 1537.9703\\
    0.5  & 0.0079 & 6.2426 & 0.2186 & 0.3500 & 1565.9526\\
    1    & 0.0054 & 8.5889 & 0.3326 & 0.4354 & 1549.7152\\
\bottomrule
\end{tabular}
\label{tab: obs gaussian}
\end{table}

Figure \ref{fig: obs gaussian} and Table \ref{tab: obs gaussian} compare the results with $\weight_b=0.05,0.1,0.5,5$. 
For a fair comparison, we initialize the algorithm with obstacle $\obs^0=\zero$ so that the initial relative errors all start from 1 for different $\weight_b$. 
We run each inner loop for 5 iterations and run the outer loop for 6000 iterations.

The first row in Figure \ref{fig: obs gaussian}
plots training data $\rhotilde(\cdot,\cdot,0.5)$ and $\vmtilde(\cdot,\cdot,0.5)$. In the first column of table \ref{tab: obs gaussian}, we report the density value $\rhotilde(\cdot,\cdot,0.5)$ at the center, reflecting the value of $\min\rhotilde$.
It is clear to see that more agents avoid the center of the obstacle as $\weight_b$ grows larger, thus the density value in the center decreases.

The third row of Figure \ref{fig: obs gaussian} presents the progression of upper-level objective values across upper-level iterations, while Table \ref{tab: obs gaussian}, Column 2, details the final upper-level objective values. To enhance the precision of the upper-level objective calculation, we execute the forward solver to convergence every 10 upper-level iterations. This approach yields a refined approximation of $(\rho^*(\obs^{(k)}),\vm^*(\obs^{(k)}))$, thereby providing a more accurate estimation of the upper-level objective values. Theorem \ref{thm-conv} implies that convergence is achieved when $\min \rhotilde > 0$. Supporting this, Table \ref{tab: obs gaussian}, Column 1, indicates that $\min\rhotilde>0$ for all considered $\weight_b$ values. Furthermore, Figure \ref{fig: obs gaussian}, Row 3, demonstrates numerical convergence for each $\weight_b$ selection. This verifies the algorithm convergence Theorem \ref{thm-conv}.

We qualitatively show the numerical solutions of the obstacle in the second row of Figure \ref{fig: obs gaussian}, while we report the relative error in the fourth row of Figure \ref{fig: obs gaussian}  and  list the best relative error and terminal step relative error in the third and fourth column of Table \ref{tab: obs gaussian}, respectively. 
Given that $\min\rhotilde>0$, Theorem \ref{thm: unique identifiablity} suggests the possibility of uniquely recovering the ground truth obstacle, up to a constant, for all $\weight_b$ values of 0.05, 0.1, 0.5, and 1. Numerically, this unique recovery is observed for $\weight_b=0.05$ and $0.1$. However, for higher $\weight_b$ values of 0.5 and 1, the reconstructed $\obs$ does not align perfectly with the ground truth $\obstilde$, as one might expect. This deviation is accounted for by Remark \ref{rem: robustness reasoning}, which discusses the robustness of the reconstruction. Specifically, when $\weight_b$ is set to 0.5 or 1, the lower bound of the data $\rhotilde$ decreases. According to Remark \ref{rem: robustness reasoning}, a smaller $\rhotilde$ lower bound leads to less robust solutions, making them more susceptible to distortions from small perturbations in the ground truth. In our experiments, since the forward solver typically produces an approximation of the exact minimizer after a finite number of iterations, the data represents a slight deviation from the ground truth. Consequently, This causes the reconstructed obstacle to differ from the exact obstacle and the discrepancy is more obvious when $\weight_b=0.5, 1$.

\subsubsection{Improving results with multiple data}
\label{sssec: multidata}

We conduct an experiment to show that multiple training data help to enhance reconstruction results for the inverse metric problem. 

The example is defined on space domain $[-0.5,0.5]$ and time domain $[0,1]$.
We discretize the space domain $[-0.5,0.5]$ with $\nx=64$ and the time domain $[0,1]$ with $\nt=16$.
The ground truth metric is $\widetilde{\metric}(x) = 0.7-0.3\cos(2\pi x)$.
The parameters in the forward problem are $\weight_I=0.01,\weight_T=0.5$. Then we obtain the first pair of data with $\rhoinit(x)=1.25-0.25\cos(4\pi x),\rhoend=1.25+0.25\cos(2\pi x)$ and the second pair with $\rhoinit(x)=p_g(x;0,0.1),\rhoend=1$.

We solve the inverse problem with the first pair of data $(N=1)$ or both data $(N=2)$.
When solving the inverse problem, we take the information on the left end $\calG_k=\{\idx:\idx=1\}$ as known and fix it.
We choose $\calR(\metric):=\half\weight_{\calR}\int\|\nabla g(x)\|_2^2\dd x$ to regularize the smoothness of the metric. 
The discretization is therefore $\calR_{\calG}(\metric):=\half\weight_{\calR}\dx\sum_{\idx=1}^{\nx-1}((\metric)_{\idx+1}-(\metric)_{\idx})^2.$

\begin{figure}[h!]
\centering
\includegraphics[width=3.9cm]{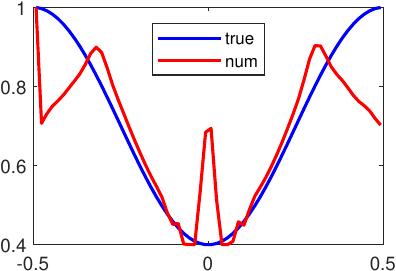}
\includegraphics[width=3.9cm]{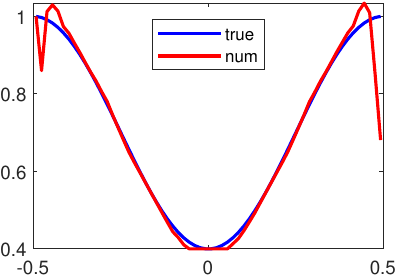}
\includegraphics[width=3.9cm]{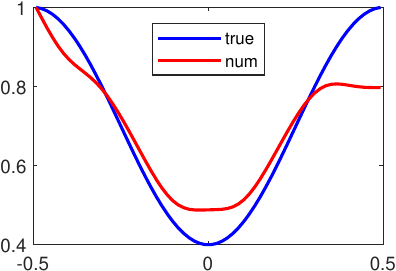}
\includegraphics[width=3.9cm]{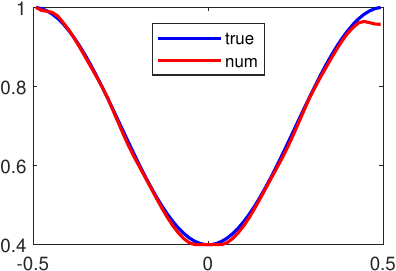}
\\
\includegraphics[width=3.9cm]{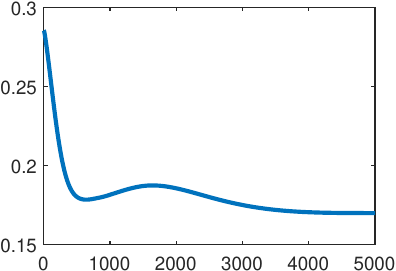}
\includegraphics[width=3.9cm]{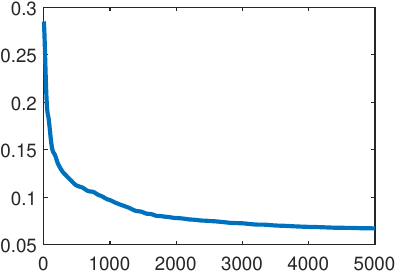}
\includegraphics[width=3.9cm]{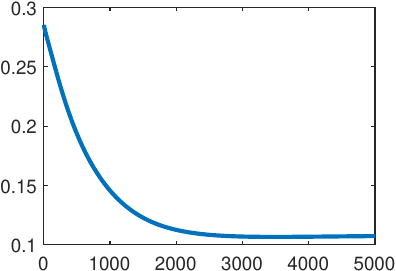}
\includegraphics[width=3.9cm]{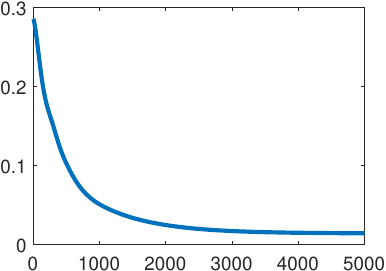}
\caption{Improving results with multiple data. 
Top to bottom: comparison of numerical $\metric$ and the ground truth $\mettilde$, the relative error of $\metric$ versus the number of iterations.
Left to right: $(N=1,\weight_{\calR}=0)$, $(N=2,\weight_{\calR}=0)$, $(N=1,\weight_{\calR}=10^{-5})$, $(N=2,\weight_{\calR}=10^{-4})$.}
\label{fig: metric 1d multidata}
\end{figure}

We run the Algorithm \ref{alg: AGM general} for 5000 iterations with 5 iterations per each inner loop. The initialization on $\idx=1$ is set as the true value and the initialization on other points is 0.7.
Figure \ref{fig: metric 1d multidata} shows the comparison of numerical results and ground truth (row 1) and the relative error of the metric versus the number of upper-level iterations.
Table \ref{tab: metric 1d multidata} reports the weight of regularization $\weight_{\calR}$, relative error, and running time of the algorithm. 
For one comparison, we choose no regularization $(\weight_{\calR}=0)$ in the model. The results with the first data $(N=1)$ are presented in row 1 and the results with both data $(N=2)$ are in row 2.
Then we tune the regularization parameter and report the best results with the first data in row 3 and with both data in row 4.
It is easy to see that when using both data, our model captures the ground truth metric better and achieves lower relative error.
It is worth noting that when using both data to solve the inverse problem, our model captures the shape of the ground truth metric even without smoothness regularization. However, when using the first data, the model fails to learn the information in the center and on both ends.

\begin{table}[htb]
\centering
\caption{Improving results with multiple data.}
\begin{tabular}{l|ccc}
\toprule
    $N$ & $\weight_{\calR}$ &  relative error & time elapsed (seconds)\\ \midrule
    1    & $0$              & 0.1700 &  60.0653\\
    2    & $0$              & 0.0673 & 128.5328\\
    1    & $1\times10^{-5}$ & 0.1073 &  67.3291\\
    2    & $1\times10^{-4}$ & 0.0145 & 118.9042\\
\bottomrule
\end{tabular}
\label{tab: metric 1d multidata}
\end{table}

\subsection{Robustness with respect to data}
\subsubsection{Unknown obstacles}
To test the robustness of our method for noisy input as discussed in Remark \ref{rem: robustness reasoning} , we design the following numerical experiment.

We discretize the space $[-0.5,0.5]^2$ with $\nx=\ny=64$ and choose $\nt=16$. 
We let the obstacle function be $\obs(x,y)=\begin{cases}
    0.5, & x<0,0.05<y<0.1, \text{ or } x>0,-0.1<y<-0.05,\\
    0, & \text{ otherwise.}
\end{cases}$.
Assume there is one pair of observations, with initial density $\rhoinit=p_g(\cdot,\cdot;-0.3,0.3,0.1,0.1)$, preferred terminal density $\rhoend=p_g(\cdot,\cdot;0.3,-0.3,0.1,0.1)$ and $\weight_I=0.1,\weight_T=1$.
We use the perturbed observation $\rhotilde+\weight_n n_{\rho},\vmtilde+\weight_n n_{\vm}$ to solve the inverse problem, where $\weight_n=0,0.25,0.5,0.75$ and noise $n_{\rho},n_{\vm}$ are generated by pointwise i.i.d sampling from the uniform distribution $U[-0.5,0.5]$.
To avoid numerical instability caused by zero value or negative density values, we threshold the perturbed density by 0.01. 
All experiments initialize with the same random choice of $\obs$. Every inner loop contains 5 iterations and 5000 outer iterations have been conducted.
In addition, we do not add any regularizer in this experiment. 
\begin{figure}[h!]
\centering
\includegraphics[width=4cm]{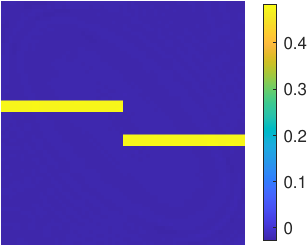}
\includegraphics[width=4cm]{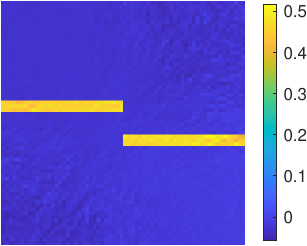}
\includegraphics[width=4cm]{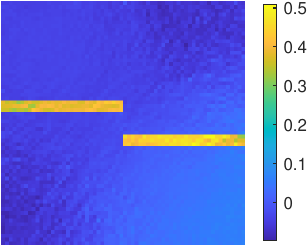}
\includegraphics[width=4cm]{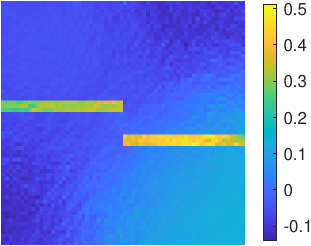}
\\
\includegraphics[width=4cm]{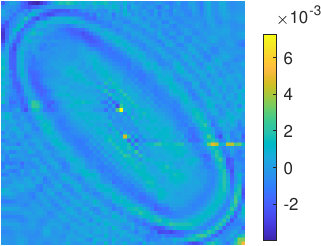}
\includegraphics[width=4cm]{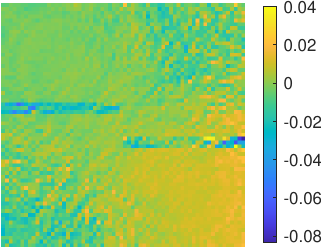}
\includegraphics[width=4cm]{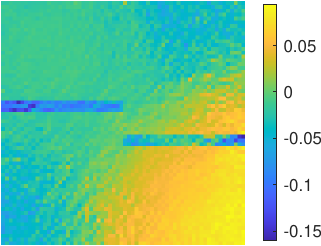}
\includegraphics[width=4cm]{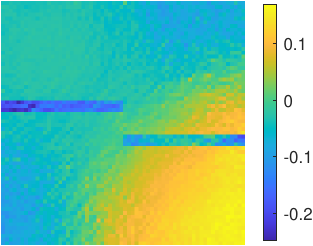}
\\
\includegraphics[width=4cm]{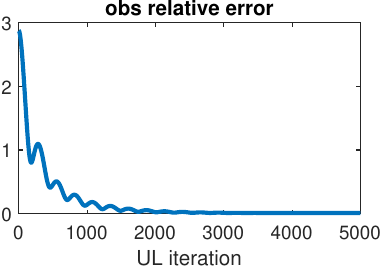}
\includegraphics[width=4cm]{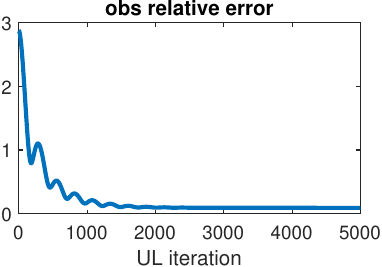}
\includegraphics[width=4cm]{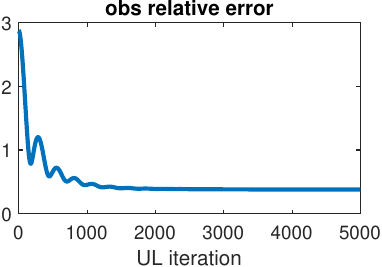}
\includegraphics[width=4cm]{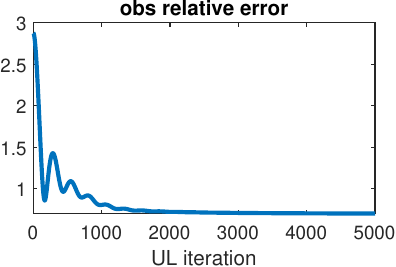}
\caption{Robustness test of the inverse crowd motion problem. 
Top to bottom: numerical $\obs$, the difference between the numerical results and the ground truth $\obs-\obstilde$, the relative error of $\obs$ versus the number of iterations. 
Left to right: noise level $\weight_n=0,0.25,0.5,0.75$. }
\label{fig: obs twobar}
\end{figure}
From Figure \ref{fig: obs twobar} and Table \ref{tab: obs robust}, we observe that with larger noise, the relative errors between numerical results and the ground truth are larger. Overall, the numerical results capture the shape of the ground truth and the algorithm converges to a close result to the ground truth $\obstilde$ with reasonably low relative errors.

\begin{table}[htb]
\centering
\caption{Robustness test of the inverse crowd motion problem.}
\begin{tabular}{l|cc}
\toprule
    $\weight_n$ &  relative error (last) & time elapsed (second)\\ \midrule
    0     & 0.0081 & 1437.0926\\
    0.25  & 0.0897 & 1343.8082\\
    0.5   & 0.3771 & 1397.4578\\
    0.75  & 0.7035 & 1379.7269\\
\bottomrule
\end{tabular}
\label{tab: obs robust}
\end{table}

\subsubsection{Unknown 1D metric}
\label{sssec: met robust}

\begin{figure}[h!]
\centering
\includegraphics[width=3.9cm]{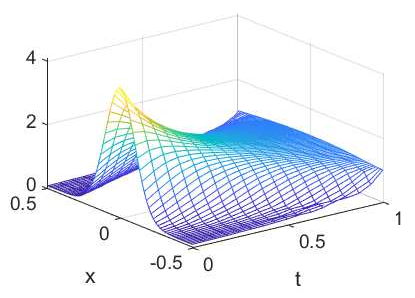}
\includegraphics[width=3.9cm]{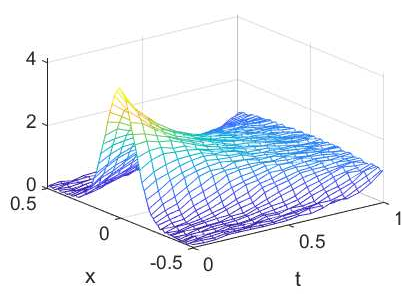}
\includegraphics[width=3.9cm]{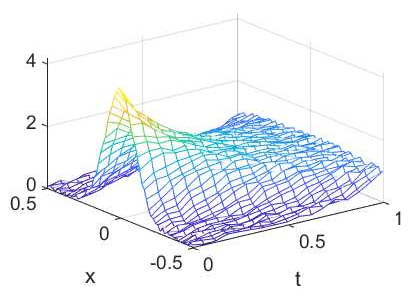}
\includegraphics[width=3.9cm]{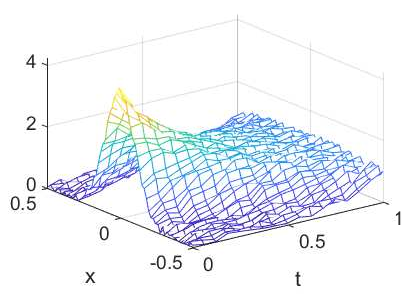}
\includegraphics[width=3.9cm]{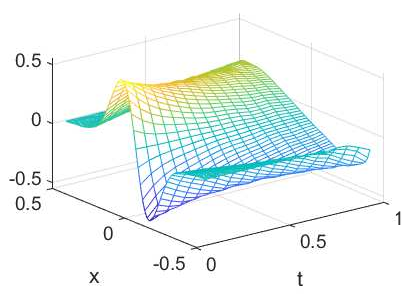}
\includegraphics[width=3.9cm]{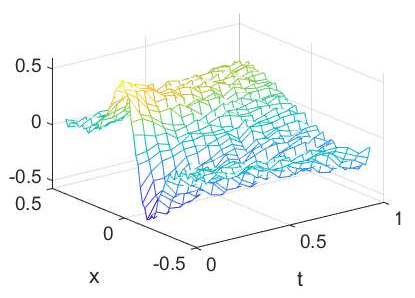}
\includegraphics[width=3.9cm]{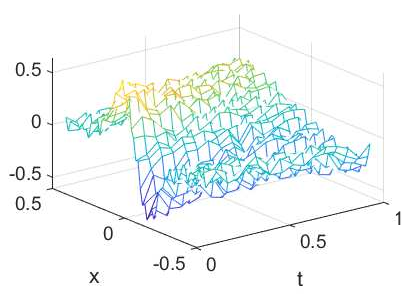}
\includegraphics[width=3.9cm]{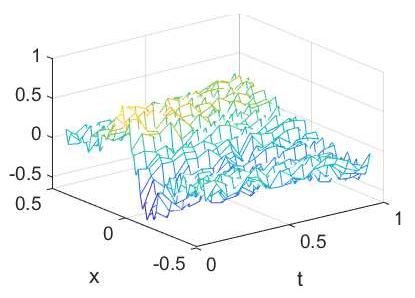}
\includegraphics[width=3.9cm]{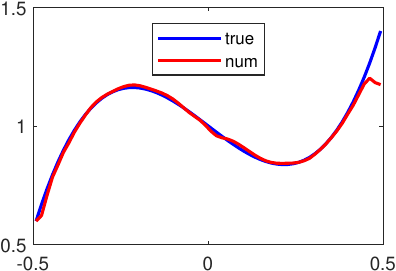}
\includegraphics[width=3.9cm]{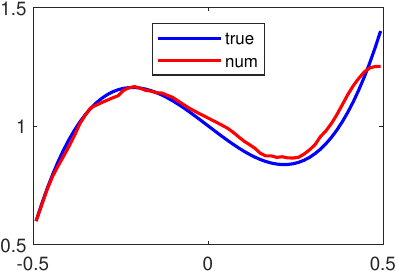}
\includegraphics[width=3.9cm]{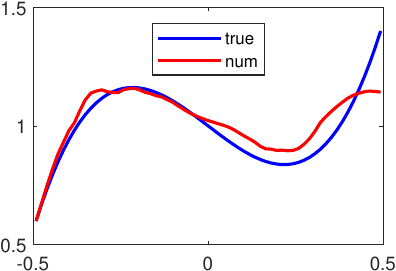}
\includegraphics[width=3.9cm]{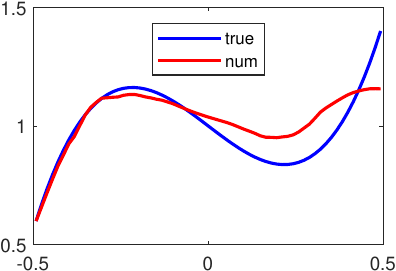}
\\
\includegraphics[width=3.9cm]{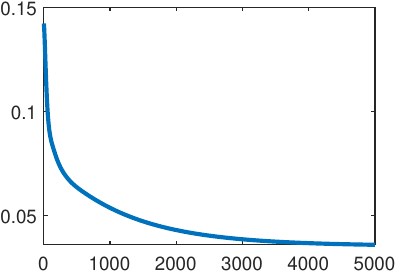}
\includegraphics[width=3.9cm]{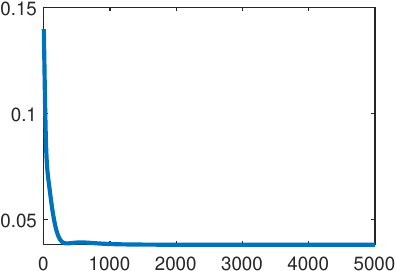}
\includegraphics[width=3.9cm]{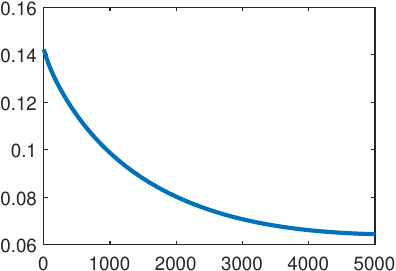}
\includegraphics[width=3.9cm]{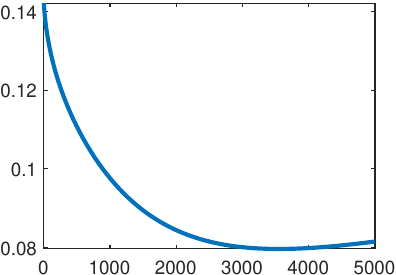}
\caption{Robustness test of the inverse metric problem. 
Left to right: $\weight_n=0,0.1,0.2,0.3$.
Top to bottom: perturbed data $\rhotilde+\weight_n n_{\rho}$, perturbed data $\vmtilde+\weight_n n_{\vm}$, comparison of numerical $\metric$ and the ground truth $\widetilde{\metric}$, the relative error of $\metric$ versus the number of iterations.}
\label{fig: metric 1d robust}
\end{figure}

This is a 1D example on $[-0.5,0.5]\times[0,1]$.
We discretize the space domain $[-0.5,0.5]$ with $\nx=64$ and the time domain $[0,1]$ with $\nt=16$.
The ground truth metric is $\widetilde{\metric}(x) = 8x(x-0.375)(x+0.375) + 1$.
The data is obtained by taking $\rhoinit(x)=p_g(x;0,0.1),\rhoend=1$ and $\weight_I=0.01,\weight_T=0.5$.
We test the robustness of the model by perturbing the observation $\rhotilde,\vmtilde$. 
The noises $n_{\rho},n_{\vm}$ share the same size with $\rhotilde,\vmtilde$ and are pointwise i.i.d samples from $U[-0.5,0.5]$.
We use the perturbed data $\rhotilde+\weight_n n_{\rho},\vmtilde+\weight_n n_{\vm}$ to solve the inverse problem, where $\weight_n=0,0.1,0.2,0.3$.
Row 1-2 of figure \ref{fig: metric 1d robust} illustrate the perturbed data. 

When solving the inverse problem, we take the information on the left end $\calG_k=\{\idx:\idx=1\}$ as known and fix it.
Same as section \ref{sssec: multidata}, we choose $\calR(\metric):=\half\weight_{\calR}\int\|\nabla g(x)\|_2^2\dd x$ to regularize the smoothness of the metric. 
The regularization weight $\weight_{\calR}$ takes different values for different $\weight_n$ and the values are in Table \ref{tab: met robust}.
We run Algorithm \ref{alg: AGM general} for 5000 iterations with 5 iterations per each inner loop. The initialization of $\metric$ takes value 1 everywhere.
Figure \ref{fig: metric 1d robust} and table \ref{tab: met robust} compare the result with different $\weight_n$.

\begin{table}[htb]
\centering
\caption{Robustness test of the inverse metric problem.}
\begin{tabular}{l|ccc}
\toprule
    $\weight_n$ & $\weight_{\calR}$ &  relative error (last) & time elapsed (second)\\ \midrule
    0    & $1\times10^{-5}$   & 0.0358 & 63.4809\\
    0.1  & $3\times10^{-4}$   & 0.0380 & 63.2121\\
    0.2  & $1\times10^{-3}$   & 0.0645 & 61.5193\\
    0.3  & $3\times10^{-3}$   & 0.0815 & 60.7215\\
\bottomrule
\end{tabular}
\label{tab: met robust}
\end{table}

From the comparison in Figure \ref{fig: metric 1d robust} and the relative error in Table \ref{tab: met robust}, we observe that as the noise level increases, the recovered metric deviates more from the ground truth. However, it is crucial to highlight that, on the whole, our model adeptly captures the underlying shape of the metric with reasonable fidelity, and the associated relative error remains consistently small. This robust performance underscores the resilience of our model in the presence of added noise to the data. 

\subsection{Robustness with respect to unknowns}
We present more numerical results to show that our method effectively recovers various types of obstacles and metrics.

\subsubsection{Unknown obstacles}

\begin{figure}[h]
    \centering
    \includegraphics[width=3.9cm]{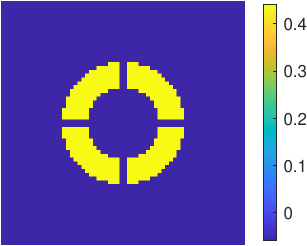}
    \includegraphics[width=3.9cm]{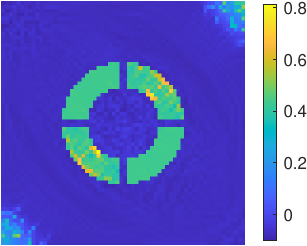}
    \includegraphics[width=3.9cm]{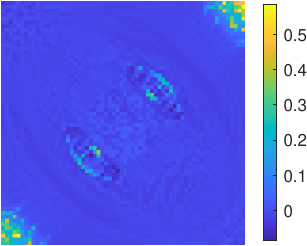}
    \includegraphics[width=3.9cm]{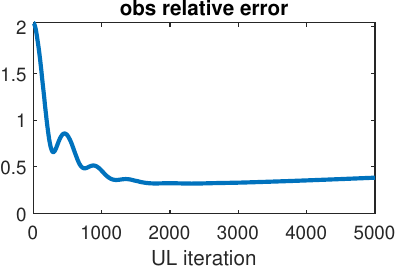}

    \includegraphics[width=3.9cm]{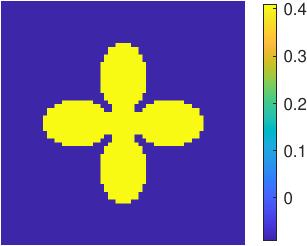}
    \includegraphics[width=3.9cm]{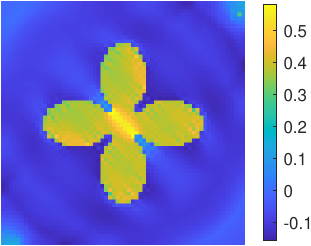}
    \includegraphics[width=3.9cm]{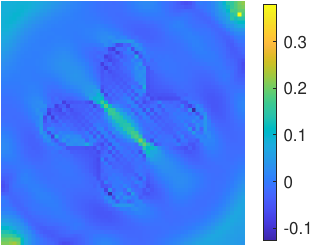}
    \includegraphics[width=3.9cm]{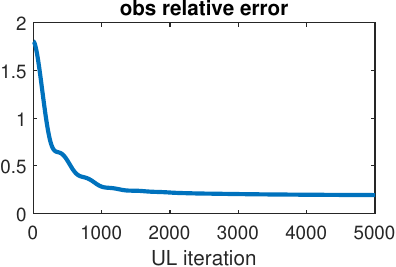}
    
    \caption{Robustness test of the inverse obstacle problem with respect to the obstacle. Mesh grid size: $\nt=16,\nx=\ny=64$. Left to right: ground truths, numerical results, the difference between ground truths and numerical results, the relative error of the obstacle versus the number of iterations. Top to bottom: relative error=0.3837, 0.1935, time elapsed=4103s, 3247s.}
\label{fig: obs var}
\end{figure}

Besides the obstacle of the Gaussian type and of a ``two-bar'' shape, we conduct experiments on obstacles with more irregular shapes.
We plot examples of ``the segmented ring'' and ``clover'' in figure \ref{fig: obs var}.
In both experiments, only one pair of data is used to recover the unknown obstacle.
The figure shows that our algorithm produces consistently good results when recovering various obstacles. Our model and algorithm recover the shape of the obstacle and achieve very low relative errors.

\subsubsection{Unknown 1D metric}

\begin{figure}[h]
    \centering
    \includegraphics[width=3.9cm]{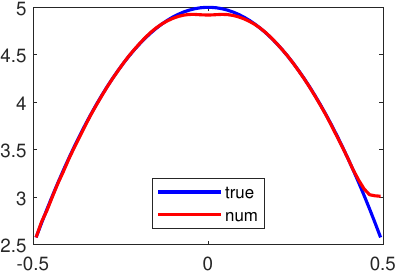}
    \includegraphics[width=3.9cm]{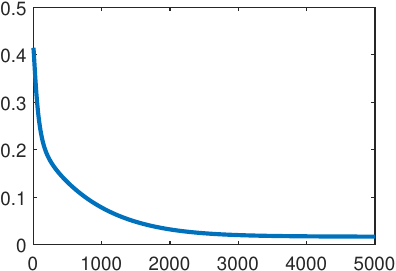}
    \includegraphics[width=3.9cm]{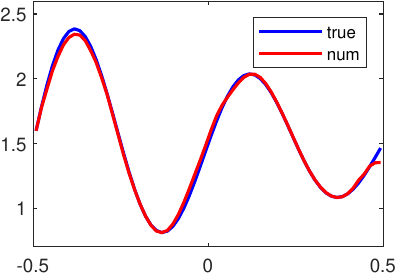}
    \includegraphics[width=3.9cm]{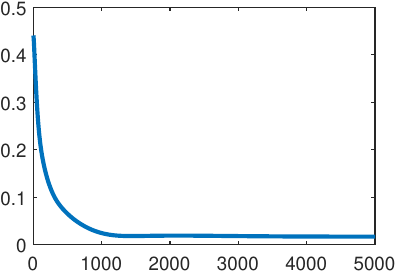}    
\caption{Robustness test of the inverse metric problem with respect to the metric. Mesh grid size: $\nt=16,\nx=\ny=64$. Columns 1,3: comparison of numerical $\metric$ and the ground truth $\mettilde$, columns 2,4: the relative error of $\metric$ versus the number of iterations. Column 1,2: $\lambda=10^{-5}$, relative error=0.0172, time elapsed=62.8513s, column 3,4: $\lambda=10^{-5}$, relative error=0.0172, time elapsed=63.0395s.}
\label{fig: met var}
\end{figure}

Apart from the experiments in sections \ref{sssec: multidata} and \ref{sssec: met robust},  we conduct experiments on more different metrics and plot the results in Figure \ref{fig: met var}.
In both experiments, we use only one pair of data and the ground truth information on the left end.
The figure shows that our model and algorithm consistently recover the ground truth metric and achieve low relative errors.

\subsection{Unknown 2D metric}
\label{subsec: num 2d met}
\begin{figure}[h]
    \centering
    \includegraphics[width=15cm]{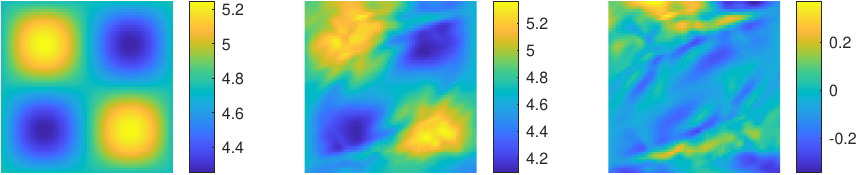}
    \includegraphics[width=15cm]{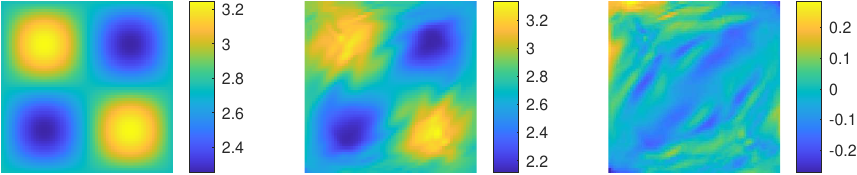}
    \includegraphics[width=15cm]{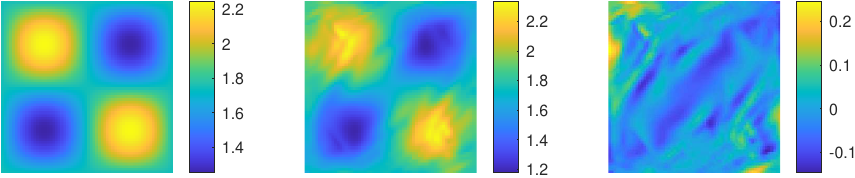}
\caption{Solving an inverse problem with an unknown metric in 2D. Mesh grid size: $\nt=16,\nx=\ny=64$. Left to right: ground truths, numerical results, the difference between ground truths and numerical results. Top to bottom: $\metric_{xx},\metric_{xy},\metric_{yy}$. Relative error=0.0260, time elapsed=4327.5671s.}
\label{fig: met 2d}
\end{figure}

The last example is a 2D inverse metric problem on $[-0.5,0.5]^2\times[0,1]$. 
We take $\nx=\ny=64$ and $\nt=16$.
The ground truth metric is $\widetilde{\metric}(x,y) = \left(\begin{matrix}
    \metric_0(x,y)+4 & \metric_0(x,y)+2\\
    \metric_0(x,y)+2 & \metric_0(x,y)+1
\end{matrix}\right)$ with $\metric_0(x,y)=0.75+0.5\sin(2\pi x)\cos(2\pi y-0.5\pi)$.
The data is obtained by taking $\weight_I=0.1,\weight_T=1$. 
We take $N=4$, i.e. 4 observations, in this example. The initial densities are $\rhoinit =p_g(\cdot,\cdot;a_x,a_y,0.1,0.1)$ with $(a_x,a_y)=(-0.3,-0.3),(-0.3,0),(-0.3,0.3),(0,0.3),$ and the terminal densities are $\rhoend(x,y)=p_g(\cdot,\cdot;a_x,a_y,0.1,0.1)$ with $(a_x,a_y)=(0.3,0.3),(0.3,0),(0.3,-0.3),(0,-0.3)$.
We solve the inverse problem with the weights of smoothness regularizers $\weight_{\calR}=10^{-4}$. 
The algorithm initiates from $\metric_{xx}=4,\metric_{xy}=2$ and $\metric_{yy}=1$. Each inner loop takes 5 iterations and each outer loop takes 5000 iterations.
Columns 1-3 of Figure \ref{fig: met 2d} shows the ground truth, the recovered metric, and the difference between the numerical result and ground truth. Our model and algorithm capture the symmetricity of the ground truth metric and achieve a relative error of value 0.0260.

\section{Conclusion}
\label{sec:con}
In conclusion, this paper introduces a novel bilevel optimization framework to tackle inverse mean-field games for learning metrics and obstacles. 
We also design an alternating gradient descent algorithm to solve the proposed bilevel problems. 
The primary advantage of our proposed formulation is its ability to retain the convexity of the objective function and the linearity of constraints in the forward problem. 
Focusing on the inverse mean-field games involving unknown obstacles and metrics, we have achieved numerical stability in these setups. 
A significant contribution of our research is establishing unique identifiability in the inverse crowd motion model with unknown obstacles based on one pair of input and revealing when the solution of the bilevel problem is stable to the noisy data. 
Employing an alternating gradient-based optimization algorithm within our bilevel approach, we ensure its convergence and illustrate its effectiveness through comprehensive numerical experiments. 
These experiments serve as robust validation, underscoring the practical applicability and reliability of our algorithm in resolving inverse problems. Our model and techniques offer a new approach to understanding and further explorations and application of inverse mean-field games.

\bibliography{bib_mfg,bib_bilevel}
\bibliographystyle{plain}

\end{document}